\documentclass[pdflatex,sn-mathphys-num]{sn-jnl}
\usepackage[english]{babel}
\usepackage[T1]{fontenc}
\usepackage[utf8]{inputenc}
\usepackage{amsmath}
\usepackage{amssymb}
\usepackage{amsfonts}
\usepackage{amsthm}
\usepackage{bm}
\usepackage{mathtools}
\usepackage{ascmac}
\usepackage{framed}
\usepackage{ulem}
\usepackage{graphicx}
\usepackage{tikz}
\usepackage{enumerate}
\usepackage{multirow}
\usetikzlibrary{intersections, calc, arrows.meta, patterns}

\newcommand{\eps}{\varepsilon}

\newcommand{\N}{\mathbb{N}}

\newcommand{\R}{\mathbb{R}}

\newcommand{\Z}{\mathbb{Z}}

\newcommand{\boB}{\mathcal{B}}

\newcommand{\boD}{\mathcal{D}}

\newcommand{\boH}{\mathcal{H}}

\newcommand{\boK}{\mathcal{K}}

\newcommand{\boM}{\mathcal{M}}

\newcommand{\boU}{\mathcal{U}}

\newcommand{\la}{\ensuremath{\lambda}}
\newcommand{\La}{\ensuremath{\Lambda}}
\newcommand{\te}{\ensuremath{\theta}}

\newcommand{\be}{\ensuremath{\beta}}
\newcommand{\gam}{\ensuremath{\gamma}}

\newcommand{\Om}{\ensuremath{\Omega}}

\newcommand{\vph}{\varphi}

\DeclareMathOperator{\Ker}{{\rm Ker}}

\DeclareMathOperator{\supp}{{\rm supp}}
\DeclareMathOperator{\Arccos}{{\rm Arccos}}

\newtheorem*{claim*}{Claim}

\newtheorem*{cor*}{Corollary}
\newtheorem{lem}{Lemma}
\newtheorem{prop}{Proposition}

\newtheorem{thm}{Theorem}
\newtheorem{alpthm}{Theorem}

\theoremstyle{definition}

\newtheorem{remark}{Remark}

\theoremstyle{remark}

\makeatletter												%
\renewcommand{\theequation}{
  \thesection.\arabic{equation}}							
  \@addtoreset{equation}{section}						%
\makeatother

\newcommand{\nor}[2]{\left\| {#1} \right\|_{#2}}		
\newcommand{\ovl}[1]{\overline{#1}}					
		
\newcommand{\inp}[2]{\left\langle {#1} , {#2} \right\rangle }	

\newcommand{\Del}{{\Delta}}							
\newcommand{\del}{{\delta}}								
\newcommand{\rd}{{\partial}}								
\newcommand{\nab}{{\nabla}}							

\newcommand{\bh}{{\mathbf{h}}}
\newcommand{\bn}{{\mathbf{n}}}
\newcommand{\boe}{{\mathbf{e}}}

\newcommand{\bom}{{\mathbf{m}}}
\newcommand{\bphi}{{\bm{\phi}}}

\newcommand{\bxi}{{\bm{\xi}}}
\newcommand{\bJ}{{\mathbf{J}}}

\newcommand{\wto}{\rightharpoonup}
\newcommand{\tif}{\tilde{f}_\be}

\begin{document}

\title{ 
Global perturbation of isolated equivariant chiral skyrmions from the Bogomol'nyi case}
\date{}
\author*[1]{\fnm{Slim} \sur{Ibrahim}}\email{ibrahims@uvic.ca}
\author[2]{\fnm{Ikkei} \sur{Shimizu}}\email{shimizu.ikkei.8s@kyoto-u.ac.jp}
\equalcont{These authors contributed equally to this work.}

\affil*[1]{\orgdiv{Department of Mathematics and Statistics}, \orgname{University of Victoria}, \orgaddress{\city{Victoria}, \state{BC}, \country{Canada}}}

\affil[2]{\orgdiv{Department of Mathematics}, \orgname{Kyoto University}, \orgaddress{\city{Kyoto}, \postcode{606-8502}, \country{Japan}}}

\abstract{
Isolated skyrmion solutions to the two-dimensional Landau–Lifshitz equation with Dzyaloshinskii-Moriya interaction, Zeeman term, and easy-plane anisotropy of various strengths are studied.
In the full range of parameter values for which the energy is a positive variation of the Bogomol'nyi case, we construct solutions to the corresponding Euler–Lagrange equation and analyze their qualitative properties, including monotonicity, exponential decay, and stability.
Our analysis is global and non-perturbative.
Moreover, we derive precise estimates quantifying the difference between these solutions and those in the Bogomol'nyi regime.
A key ingredient of our approach is a novel resolvent estimate for the linearized operator, which remains uniform with respect to additional implicit potentials arising in the problem.

}

\keywords{skyrmion, Landau-Lifshitz energy, Dzyaloshinskii-Moriya interaction, harmonic map\\
\textbf{Mathematics Subject Classification: 
} 35B25, 35C08, 49S05, 82D40
}

\maketitle


\section{Introduction}

\subsection{Introduction}


Micromagnetics is a continuum theory that describes the static and dynamic behavior of ferromagnets \cite{LanLif35, HubSch, BogYab89, BogHub94}. The state of a ferromagnet is characterized by its magnetization, a unit-length vector field
\[
\bn(x) =
\begin{pmatrix}
n_1(x) \\ n_2(x) \\ n_3(x)
\end{pmatrix},
\]
which is a map from $\R^2$ into the unit sphere $\mathbb{S}^2\subset\R^3$. The equilibrium configurations of a ferromagnet are the critical points of the Landau-Lifshitz energy functional, which for our purposes is given by
\begin{equation}\label{E1.1}
E_{r,h,k}[\bn] = D[\bn] + rH[\bn] + h\, Z[\bn] + k\, A[\bn],\qquad r>0,\quad h,k\in\R,
\end{equation}
where the individual terms are defined as
\begin{align*}
  D[\bn] &= \frac{1}{2} \int_{\R^2} |\nabla \bn|^2 \,dx, \quad\text{(Dirichlet energy)},\\[1mm]
  H[\bn] &= \int_{\R^2} (\bn - \boe_3) \cdot (\nabla \times \bn)\, dx, \quad\text{(Dzyaloshinskii-Moriya interaction)},\\[1mm]
  Z[\bn] &= \int_{\R^2} (1-n_3)\, dx, \quad\text{(Zeeman interaction)},\\[1mm]
  A[\bn] &= \frac{1}{2} \int_{\R^2} (1-n_3^2)\, dx, \quad\text{(easy-plane anisotropy)},
\end{align*}
with $\boe_3=(0,0,1)^T$. For clarity, on $\R^2$ we define
\[
\nabla \times \bn = 
\begin{pmatrix}
\partial_1 \\ \partial_2 \\ 0
\end{pmatrix}
\times 
\begin{pmatrix}
n_1 \\ n_2 \\ n_3
\end{pmatrix}.
\]
We work on the space
\[
H^1_{\boe_3} := \big\{ \bn:\R^2\to \R^3 \mid \bn - \boe_3 \in H^1(\R^2,\R^3),\; |\bn|=1 \text{ a.e.} \big\},
\]
which ensures that each term in \eqref{E1.1} is finite since
\[
Z[\bn] = \frac{1}{2}\int_{\R^2} |\bn - \boe_3|^2\,dx < \infty
\]
and
\[
A[\bn] = \int_{\R^2}\frac{1}{2}(1+n_3)(1-n_3)\,dx \le Z[\bn].
\]

The Euler-Lagrange equation associated with \eqref{E1.1} is
\begin{equation}\label{EF1.5}
-\Delta \bn + 2r\, \nabla \times \bn - (h+kn_3)\boe_3 - \Bigl(|\nabla \bn|^2 + 2r\, \bn\cdot (\nabla\times \bn) - (hn_3 + kn_3^2)\Bigr)\bn =0,
\end{equation}
which is a form of the Landau-Lifshitz equation. By classical regularity theory for harmonic maps, every $H^1_{\boe_3}$ solution to \eqref{EF1.5} is smooth (see the Appendix for an outline of the proof).

The variational problem for \eqref{E1.1} serves as a model for the equilibrium states of ferromagnetic materials. In particular, the Dzyaloshinskii-Moriya interaction term $H[\bn]$ is of great interest as it underlines the stabilization of spatially inhomogeneous states (see, e.g., \cite{BogYab89, BogHub94, BogHub99, Leo16}). From a mathematical standpoint, the variational analysis of \eqref{E1.1} has attracted significant attention \cite{Mel14, DorMel17, LiMel18, Sch19, BarSinRosSch20, 
BerMurSim20, BerMurSim21, IbrShi23, MurSimSla}. For instance, Melcher \cite{Mel14} established the existence of minimizers in the homotopy class of degree $-1$ for $k=0$ and $h\ge r^2$, while Li and Melcher \cite{LiMel18} constructed equivariant solutions (defined later) for $k=0$ and large $h$, showing that these solutions are local minimizers among symmetry-breaking perturbations. Further extensions to the equivariant setting and asymptotic analyses as $h,k\to\infty$ have been carried out in \cite{GusWan21, KomMelVen20, KomMelVen21, KomMelVen23}.\par
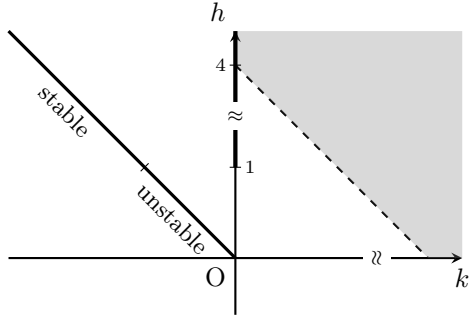
\begin{figure}
\centering
\begin{tikzpicture}[scale=1.5]
	\draw [thick, -stealth](-2,0)--(2,0) node [anchor=north]{$k$};
	\draw [thick, -stealth](0,-0.5)--(0,2) node [anchor=south east]{$h$};
	\node [anchor=north east] at (0,0) {O};        
	\draw [thick, dashed] (0,1.7)--(1.7,0);
	\fill [gray, opacity=0.3]
	 (0,2)	-- (0, 1.7)
	-- (1.7,0) -- (2,0) 
	-- (2,2) -- cycle;
	\draw [ultra thick] (0,0.8) -- (0,2);
	\draw (-0.05,0.8) -- (0.05,0.8);
        \draw (-0.05,1.7) -- (0.05,1.7);
 	\node [right] at (0,0.8) {\footnotesize $1$}; 
 	\node [left] at (0,1.7) {\footnotesize $4$}; 
	\node [fill=white] at (0,1.25) {\footnotesize $\approx$};
	\node [fill=white] at (1.25,0) {\footnotesize \rotatebox{90}{$\approx$}};
	\draw [very thick] (0,0)--(-2,2);
	\draw (-0.765,0.835)--(-0.835,0.765);
	\path (0,-0.2) -- (-1.1, 0.9) node[midway,sloped] {\scriptsize unstable} ;
	\path (-1,0.8) -- (-2, 1.8) node[midway,sloped] {\scriptsize stable} ;
\end{tikzpicture}
\caption{The $hk$ diagram summarizing the known mathematical results on $E_{1,h,k}$. 
\cite{Mel14} and \cite{LiMel18} investigate the bold line on $h$-axis. 
The gray region is covered by 
\cite{GusWan21}. 
The Bogomol'nyi case corresponds to the sloped line lying in the region $k\le 0$, examined in \cite{BarSinRosSch20, DorMel17, IbrShi23}.
}
\label{Fig1}
\end{figure}

An interesting factorization structure arises when $h+k=0$ 
leading to 
additional insights into the
structure of solutions. This is referred to as the Bogomol'nyi trick, as can be commonly seen in nonlinear field theories (see \cite{Felsager}). 
More specifically, we have
\begin{equation}\label{EG1.2}
\begin{aligned}
E_{r,h,-h} [\bn] = 
-4\pi \left(1-\frac{2r^2}{h}\right) Q[\bn]
+ \frac {r^2}{2h} \int_{\R^2} 
&|\boD_1^{\kappa} \bn + \bn\times \boD_2^{\kappa} \bn|^2 dx \\
& + 
\frac 12
\left(
1-\frac{r^2}{h} 
\right)\int_{\R^2} 
|\rd_1 \bn - \bn\times \rd_2 \bn|^2 dx,
\end{aligned}
\end{equation}
where $Q[\bn] = \int_{\R^2} \bn\cdot \rd_1 \bn\times \rd_2\bn dx$ is the topological degree of $\bn$, and $\boD_j^{\kappa} = \rd_j - \frac {1}{\kappa} \boe_j\times \cdot$ for $j=1,2$. 
and $\kappa= \frac{r}{h}$. 
In \cite{DorMel17}, D\"oring and Melcher realized from \eqref{EG1.2} that \eqref{EF1.5} has an explicit solution
\begin{equation}\label{EG1.4} 
\bh^{\kappa} (x) := 
\left(
\frac{-4\kappa x_2}{|x|^2 +4\kappa^2} ,
\frac{4\kappa x_1}{|x|^2 +4\kappa^2} ,
\frac{|x|^2 -4\kappa^2}{|x|^2 +4\kappa^2} 
\right)^T
,
\end{equation}
which is a common solution to 
\begin{equation}\label{EG1.3}
\boD_1^\kappa \bn + \bn\times \boD_2^\kappa \bn =0\quad \text{and}\quad 
\rd_1 \bn - \bn\times \rd_2 \bn=0.
\end{equation}
Note that $\bh^\kappa$ is a harmonic map with degree $Q=-1$. 
In \cite{Sch19} and \cite{BarSinRosSch20}, 
the solution family of the first equation in \eqref{EG1.3} has been completely given via the stereographic coordinates of the sphere. 
Here, 
one can obviously see from \eqref{EG1.2} that when $\frac{r}{\sqrt{h}}\le 1$, $\bh^{\kappa}$ is an energy minimizer in the homotopy class of degree $Q=-1$. 
On the other hand, if $\frac{r}{\sqrt{h}}>1$, 
the authors \cite{IbrShi23} showed that $\bh^{\kappa}$ is energetically unstable. 
We visualize the range of parameters of the known results 
in the $h-k$ plane, fixing $r=1$, see Figure \ref{Fig1}.
\par
%
%
%
%
The focus of this paper is 
to extend the result 
away from the Bogomol'nyi case. 
More specifically, 
we work on the parameter region as
$$
\frac{h-k}2 >0,\qquad \text{and}\qquad \frac{h+k}2>0.
$$
By rescaling, 
we may normalize $\frac{h-k}2 =1$, without any loss of generality. 
In fact, if we let $\bn (x) = \tilde \bn (y)$ with $y=\la x$ for $\la>0$, we have
$$
E_{r,h,k}[\bn] 
= D[\tilde\bn] 
+\frac{r}{\la} H[\tilde\bn] + \frac{h}{\la^2} Z[\tilde\bn] + \frac{k}{\la^2} A[\tilde\bn],
$$
and thus taking $\la= \sqrt{\frac{h-k}2}$ leads to the desired reduction. 
Henceforth, we reformulate  the variational problem as
\begin{equation}\label{EF1.6}
E_{r,\be} [\bn] =
D[\bn] + rH[\bn] + V_-[\bn] + \be^2 V_+[\bn],\qquad \be>0,
\end{equation}
where  $\be^2=\frac{h+k}{h-k}$ and
\begin{align*}
    V_-[\bn] &= 
    Z[\bn]-A[\bn] =\frac 12
    \int_{\R^2} (1-n_3)^2 dx,\\
    V_+[\bn] &= 
    Z[\bn] + A[\bn]
    =
    \int_{\R^2} 
    \left(2-\frac 12 (1+n_3)^2\right) dx.
\end{align*}
The advantage of this framework is that now both $V_\pm$ are positive functionals on $H_{\boe_3}$. 
Now the Bogomol'nyi case corresponds to $\be=0$. 
Hence in our parameter range, 
\eqref{EF1.6} 
is a positive 
variation of the Bogomol'nyi case. 
When $\beta$ is small, one might expect standard perturbation from $\bh^r$ to work. However, since
\begin{equation}\label{EG1.5}
V_+[\bh^r] =\infty,
\end{equation}
this becomes a singular perturbation problem: the reference minimizer lacks the integrability needed to define the perturbed energy. 
%
%
%
%

In this paper, we restrict ourselves to maps $\bn$ in a certain symmetry class of $\bn$ that also includes $\bh^r$. 
A map $\bn$ is called \textit{equivariant} if 
$\bn$ is invariant under the simultaneous rotation of the ambient and target spaces:
\begin{align*}
\bn(x) \mapsto e^{-\phi R}\bn(e^{\phi \ovl{R}}x), &&
\phi\in\R,&&
R= 
\begin{pmatrix}
0 & -1 & 0\\
1 & 0 & 0\\
0 & 0 & 0
\end{pmatrix}
,&&
\ovl{R}
=
\begin{pmatrix}
0 & -1 \\
1 & 0 
\end{pmatrix}
.
\end{align*}
Using polar coordinates $(\rho,\psi)$ of $\R^2$, 
one can rewrite 
such a map 
in the form
\begin{align}\label{E1.3}
  \bn =
  \begin{pmatrix}
    -\sin \psi \sin f(\rho) \\
    \cos \psi \sin f(\rho)  \\
    \cos f(\rho)
  \end{pmatrix}
  ,
\quad
f:(0,\infty)\to\R.
\end{align}
In particular, we focus on equivariant maps with a profile $f$ satisfying the boundary condition
\begin{align}\label{EZG1}
    f(0)=\pi,&& f(\infty)=0.
\end{align}
In physics, 
solutions to \eqref{EF1.5} with \eqref{E1.3} and \eqref{EZG1} are called chiral magnetic skyrmions, and known as isolated skyrmions of Bloch type. 
These are topological magnetic solitons that 
behave like particles in real materials. They are of potential important use in  storage devices. For more physical background and motivation, we refer interested readers to 
\cite{BogHub94, NagTok13, Leo16, FerCroSam13, ButLemBea18} and the references therein. 
If we insert \eqref{E1.3} into \eqref{E1.1}, we can rewrite
$$
E_{r,\be} [\bn] =
2\pi \left(
\tilde{D}[f] + r \tilde{H}[f]  
+ \tilde{V}_-[f] + \be^2 \tilde{V}_+[f]
\right)
,
$$
where
\begin{align*}
\tilde{D}[f] = \frac 12 \int_0^\infty \left( |f'|^2 + \frac{\sin^2 f}{\rho^2} 
\right) \rho d\rho,
,&&
\tilde{H}[f] = \int_0^\infty \left(f' - \frac{\sin f}{\rho} \right) (1-\cos f) \rho d\rho,
\end{align*}
\begin{align*}
\tilde{V}_-[f] = \frac 12 
 \int_0^\infty (1-\cos f)^2 \rho d\rho,
&&
\tilde{V}_+ [f] =
\int_0^\infty 
\left(2- \frac 12 (1+\cos f)^2\right) \rho d\rho.
\end{align*}
Since the multiple $2\pi$ is no longer our concern, 
we redefine the energy with respect to $f$ as
\begin{equation}\label{E1.4}
\tilde{E}_{r,\be}[f] := \tilde{D}[f] + r \tilde{H}[f] + \tilde{V}_-[f] + \be^2 \tilde{V}_+[f].
\end{equation}
The corresponding function space for $f$ is 
$$
\boM := \left\{
f:(0,\infty)\to \R\ |\ 
f, \rd_\rho f, \frac{\sin f}{\rho} \in L^2_{\rho d\rho}
,\quad
f(0)=\pi,\quad f(\infty)=0
\right\}.
$$
Note that $\boM$ is still not a normed space because of the boundary condition at the origin. However, 
the $H^1$-distance of $\bn$ is translated into the terminology of norm for $f-g$ as
$$
d_\boM (f,g) := \nor{f-g}{L^2_{\rho d\rho}} + 
\nor{\rd_\rho (f- g)}{L^2_{\rho d\rho}} + 
\nor{\frac {f-g}\rho}{L^2_{\rho d\rho}},\quad f,g\in \boM
$$
by noting that $(f-g)(0)=0$ for $f,g\in \boM$. 
Accordingly, we introduce the normed space 
$$
X:= \left\{ 
f:(0,\infty)\to \R\ |\ 
\nor{f}{X} <\infty
\right\}
$$
where
$$
\nor{f}{X}^2 := \nor{\rd_\rho f}{L^2_{\rho d\rho}}^2 + \nor{\frac{f}\rho}{L^2_{\rho d\rho}}^2.
$$
Since $\nor{f}{X} \sim \nor{\nab (e^{i\psi}f(\rho))}{L^2(\R^2)}$, $X$ can be regarded as $\dot{H}^1$ for equivariant functions on $\R^2$. 
Moreover, $f\in X$ enjoys the properties
\begin{align}\label{EH2}
\lim_{\rho\to 0} f(\rho) = \lim_{\rho\to \infty} f(\rho)= 0,&&
\nor{f}{L^\infty} \le C \nor{f}{X}.
\end{align}
Indeed, for $0<\rho_1<\rho_2<\infty$, we have
\begin{equation}\label{EH1}
\begin{aligned}
|f(\rho_1)^2 - f(\rho_2)^2|
&\le 2 \int_{\rho_1}^{\rho_2} |\rd_\rho f(\rho) f(\rho)| d\rho \\
&\le 2 \left(\int_{\rho_1}^{\rho_2}
\frac{f(\rho)^2}{\rho^2} \rho d\rho
\right)^{\frac 12}
\left(\int_{\rho_1}^{\rho_2}
(\rd_\rho f(\rho))^2
\rho d\rho
\right)^{\frac 12} .
\end{aligned}
\end{equation}
This shows that $\{f(\rho)\}_{\rho>0}$ is Cauchy as $\rho\to 0$ and $\rho\to\infty$, whose limits must be $0$. By letting \eqref{EH1} as $\rho_2\to\infty$, we obtain the last inequality in \eqref{EH2}. 
Now, the Euler-Lagrange equation for $f$ becomes
\begin{equation}\label{E1.5}
  \begin{aligned}
     & -f''(\rho) - \frac 1\rho f'(\rho) + \frac 1{\rho^2} \sin f(\rho) \cos f(\rho) \\
     & - 2r \frac{\sin^2 f(\rho)}{\rho} + \sin f(\rho) (1- \cos f(\rho))
    + \be^2 \sin f(\rho) (1+\cos f(\rho)) =0.
  \end{aligned}
\end{equation}
On the other hand, the harmonic map $\bh^{r}$ also satisfies the equivariant symmetry with a profile
$$
f(\rho) = \te(\rho),\qquad 
\te(\rho) := 
\Arccos \frac{\rho^2-4r^2}{\rho^2+4r^2},
$$
where the branch of $\Arccos$ is taken over $[0,\pi]$. 
%
One can directly see that $\te$ solves \eqref{E1.5} with $\be=0$. In addition, $\te$ satisfies
\begin{align}\label{1.135}
\sin \te = \frac{4r\rho}{\rho^2+4r^2}, &&
\te' = -\frac{\sin \te}{\rho},&&
\frac{2r\sin\te}{\rho} + \cos\te -1 =0.
\end{align}
As noted in \eqref{EG1.5}, we have $\tilde{V}_+ [\te] =\infty$. 
More precisely, $\te'$, $\frac{\sin{\te}}{\rho}
\in L^2_{\rho d\rho}$, $\te\in L^p_{\rho d\rho}$ with $p>2$, but $\te \not\in L^2_{\rho d\rho}$. 

\subsection{Main results}

%
%
Our first main result is the existence and shape property of the solution to \eqref{E1.5}. 
%
%
%
%
%
%
\begin{alpthm}
\label{TD}
For $\be>0$ and $r>0$, 
there exists a solution $f_{r,\be}\in \boM$ to \eqref{E1.5} such that the following conditions hold.
\begin{itemize}
\item (Pointwise bound and exponential decay) 
For $\rho\in (0,\infty)$, we have
\begin{equation}\label{EA0.1}
0< f_{r,\be}(\rho) < \te(\rho).
\end{equation}
Moreover, for any $0<\del<1$, there exist $C_\del=C_\del(r,\be)$, $\rho_\del=\rho_\del(r,\be)>0$ such that
\begin{align}\label{EA0.3}
C_\del
e^{-(1+\del)\be\rho} \le f_{r,\be} (\rho) \le 
C_\del e^{-(1-\del)\be\rho},
&&
-C_\del e^{-(1-\del)\be\rho} <f_{r,\be}'(\rho) < -C_\del e^{-(1+\del)\be\rho}
\end{align}
for $\rho\ge \rho_\del$.
\item (Difference estimate) 
There exists $\be_0(r)>0$ such that if $0<\be < \be_0(r)$, then 
for any $\del>0$, there exists $C_\del>0$ satisfying
\begin{equation}\label{EA0.2}
\nor{f_{r,\be} - \te}{X} \le C_\del \be^{1-\del}.
\end{equation}
%
%
%
%
\item (Monotonicity)  
If $r\le 1$ and $\be\le 1$, 
then $f_{r,\be}$ is monotonically decreasing on $(0,\infty)$. 
More precisely, we have
$$
f_{r,\be}'(\rho) + \frac 1\rho \sin f_{r,\be}(\rho) <0\qquad \text{on } (0,\infty).
$$
\end{itemize}
\end{alpthm}

The second part of our main result is the stability of corresponding equivariant map $\bn$. 

\begin{alpthm}[Stability]
\label{TE}
Let $\be>0$, $r>0$, and let $f_{r,\be}$ be a solution to \eqref{E1.5} as in Theorem \ref{TD}. 
Then, we have the following two statements.
\begin{itemize}
\item (Stability) 
Assume $r\le \frac 12$, then the corresponding equivariant map $\bn$ as in \eqref{E1.3} is a linearly stable critical point of $E_{r,\be} [\bn]$, i.e. $\nab^2 E_{r,\be}[\bn]$ is nonnegative definite. 
\item (Instability) 
For any $r>1$, there exists $\be_1(r) < \be_0(r)$ such that if $0<\be<\be_1(r)$,
the corresponding equivariant map $\bn$ as in \eqref{E1.3}
is an unstable critical point of $E_{r,\be} [\bn]$, i.e., for any neighborhood $\boU$ of $\bn$ in $\boM$, there exists $\bom\in \boU$ such that $E_{r,\be}[\bom]- E_{r,\be}[\bn]<0$.
\end{itemize}
\end{alpthm}

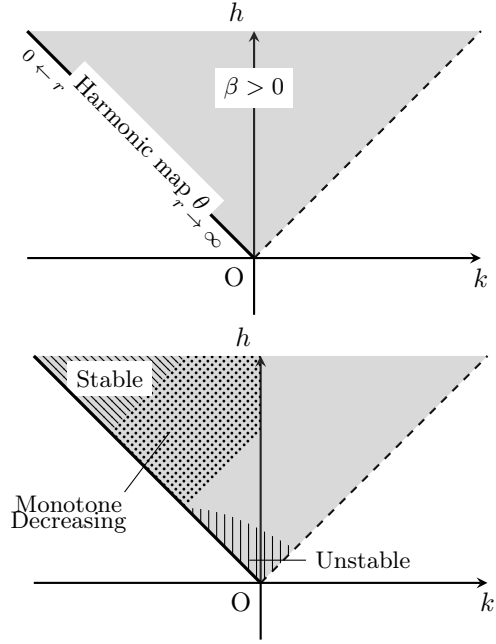
\begin{figure}
\centering
\begin{tikzpicture}[scale=1.5]
    \draw [thick, -stealth](-2,0)--(2,0) node [anchor=north]{$k$};
    \draw [thick, -stealth](0,-0.5)--(0,2) node [anchor=south east]{$h$};
    \node [anchor=north east] at (0,0) {O};
    \draw [thick, dashed](0,0)--(2,2);
    \fill [gray, opacity=0.30]
	 (-2, 2)
	-- (0, 0)
	-- (2,2)
     -- cycle;
    \draw [very thick](0,0)--(-2,2) node[midway, sloped, fill=white]{{\scriptsize Harmonic map $\te$}};
    \path (0,-0.2) -- (-1, 0.8) node[midway,sloped] {{\footnotesize $r \to \infty$}} ;
    \path (-1.7,1.5) -- (-2, 1.8) node[midway,sloped] {{\footnotesize $0 \gets r$}} ;
    \node [fill=white] at (0,1.5) {{\small $\be>0$}};
      \end{tikzpicture}
%
%
    \begin{tikzpicture}[scale=1.5]
    \draw [thick, -stealth](-2,0)--(2,0) node [anchor=north]{$k$};
    \draw [thick, -stealth](0,-0.5)--(0,2) node [anchor=south east]{$h$};
    \node [anchor=north east] at (0,0) {O};     
    \draw [thick, dashed](0,0)--(2,2);
    \fill [gray, opacity=0.3]
	 (-2, 2)
        -- (0, 0)
	-- (2,2)
        -- cycle;
    \fill [pattern=north west lines]
	 (-2, 2)
	-- (-1.33, 1.33)
	-- (-0.66,2)
        -- cycle;
    \fill [pattern=crosshatch dots]
	 (-1.33, 1.33)
	-- (-0.66, 0.66)
	-- (0,1.33)
	-- (0,2)
	-- (-0.66,2)
     -- cycle;
    \draw (-0.8, 1.33) -- (-1.2, 0.8); 
    \node[anchor=north east, align=center] at (-1.1,0.85)
    {{\scriptsize Monotone} \\[-6pt] {\scriptsize Decreasing}};
    \node[fill=white] at (-1.33,1.8) {{\scriptsize Stable}};
    \fill[pattern=vertical lines] (-0.66,0.66)
    cos (0.3,0.3) -- (0,0) -- cycle;
    \draw (-0.1, 0.2) -- (0.4,0.2) node[anchor=west] {{\scriptsize Unstable}};
    \draw [very thick](0,0)--(-2,2) ;
    \end{tikzpicture}
\caption{The 
upper diagram represents the corresponding $(h,k)$ region where the existence of solutions to \eqref{E1.5} is proved in Theorem \ref{TD}. The lower diagram indicates the $(h,k)$ region where monotonicity, stability, and instability are shown in Theorems \ref{TD} and \ref{TE}.
}\label{Fig2}
\end{figure}

%
%
%
\begin{remark}
We conjecture that solutions to \eqref{E1.5} in $\boM$ with \eqref{EA0.1} are unique, which will be addressed in our future work. 
We note that exponential decay, monotonicity, and stability stated in Theorems \ref{TD} and \ref{TE} hold for all such solutions. 
Moreover, the instability holds for any solution $f$ to \eqref{E1.5} with $0<\be\ll 1$, $0<f<\te$ and $\nor{f-\te}{X}=o_{\be\to 0}(1)$. 
See Theorem \ref{T3.2} below.
\end{remark}
\begin{remark}
Let us plot our results on $hk$-plane with fixed value of $r$. 
If we apply the rescaling $f(x) = \ovl{f}(y)$ with $y= r x$, 
we have
\begin{align*}
\tilde{E}_{r,\be}[f] 
&= 
\tilde{D}[\ovl{f}] + \tilde{H} [\ovl{f}] 
+ \tilde{V}_-[\ovl{f}] + \be^2 \tilde{V}_+[\ovl{f}] \\
&= 
\tilde{D}[\ovl{f}] + \tilde{H} [\ovl{f}] 
+ \frac{\be^2+1}{r^2} \tilde{Z}[\ovl{f}] 
+ \frac{\be^2-1}{r^2} \tilde{A}[\ovl{f}].
\end{align*}
Letting $h= \frac{\be^2+1}{r^2}$ and $k=\frac{\be^2-1}{r^2}$, 
one can visualize the range of $(h,k)$ where the isolated skyrmion solution is constructed in Theorem \ref{TD}, 
see Figure 1. 
On the other hand, the range where the monotonicity, stability, and instability have been proven is depicted in Figure \ref{Fig2}. 
\end{remark}
%
%

\subsection{Difficulty and Idea of the proof}

\subsubsection*{Construction of solutions} 
The basic idea to construct solutions to \eqref{E1.5} is based on the direct method for $\tilde{E}_{r,\be}$. 
However, if we simply 
discuss minimizing problem on $\boM$, we would encounter two difficulties. First, $\tilde{E}_{r,\be}$ may be unbounded from below if $r>1$. Second, even if a sequence in $\boM$ weakly converges in suitable spaces, 
the weak limit may satisfy different boundary conditions. 
To overcome the first issue, 
we introduce the pointwise restriction $0\le f(\rho) \le \te(\rho)$ on the function space. 
Thanks to this, the energy becomes bounded from below (see Lemma \ref{L1}).
Moreover, since $0$ and $\te$ are respectively a sub-, super-solution, minimizers still solve the Euler-Lagrange equation \eqref{E1.5} despite of the restriction of space, which is known as the so-called subsolution-supersolution method (see \cite{Struwe} for example). We will see this in Lemma \ref{L2.3}. 
Nevertheless, the topological change is still possible to occur under the weak limit of minimizing sequence. 
To cope with this, we separate our argument into two regimes of $\be$. When $\be$ is sufficiently large, then we may exclude the possibility of topological change by combining the upper estimate of the minimum (Lemma \ref{L2.2}) and the least energy of bubble (Lemma \ref{L2.1}). For the other values of $\be$, we instead discuss energy minimizing problem by restricting the function space into $f_0(\rho) \le f(\rho)\le \te(\rho)$, where $f_0(\rho)$ is a solution to \eqref{E1.5} with \eqref{EA0.1} for sufficiently large $\be_0$. 
Since $f_0$ is a subsolution to \eqref{E1.5} if $0<\be<\be_0$, the same methodology still works to conclude that the minimizer solves \eqref{E1.5}. 
\subsubsection*{Difference estimate}
%
To obtain \eqref{EA0.2}, we work with the perturbative analysis of the ODE \eqref{E1.5}. Writing $f_{r,\be}=\te +\xi$, one obtains the equation of $\xi$ from \eqref{E1.5} as
\begin{align}
\left(
F^*F + \frac{8r^2}{\rho^2+4r^2} +2\be^2
\right) \xi + \be^2 \sin \te (1+\cos \te) + N(\xi,\be)=0, \quad
F= \rd_\rho + \frac{\cos \te}{\rho},
\end{align}
where 
$N(\xi,\be)$ consists of 
higher order terms with respect to $\xi$ and $\be$. 
Since $F^*F + \frac{8r^2}{\rho^2+4r^2}$ is positive, 
one expects that the fixed point argument would provide the desired estimate by using resolvent estimates. 
In fact, in the relevant analysis in \cite{GusWan21}, 
this strategy actually works once the issue caused by the threshold resonance of the linearized operator around the harmonic map is removed by taking a specific rescaling. 
However, in our case, we encounter the essential issue that 
$N(\xi,\be)$ includes some bad nonlinear terms preventing to close the nonlinear estimate even if we take $\be$ any small. 
In detail, the resolvent of the source term $\be^2\sin\te (1+\cos \te)$ is estimated by $O(\be^{1-\del})$ for any $\del>0$ (see \eqref{EX2.4}). 
On the other hand, if we suppose to discuss fixed point argument in the space $\{\xi\in X\ |\ \nor{\xi}{X} \le O(\be^{1-\del})\}$, 
then the terms $\frac{4r}{(\rho^2+4r^2)^{1/2}}\xi^2$ and $\frac 12 \xi^3$ included in $N(\xi,\be)$ are only controlled by $O(\be^{1-2\del})$, $O(\be^{1-3\del})$, respectively, and thus the estimate could not be closed. Even if we could improve the estimate of source term into  $O(\be)$, as seen in the case of \cite{GusWan21}, we still confront the same issue since the order of these bad terms just coincides, and the smallness of $\be$ does not matter to close the estimate. 
To avoid this obstruction, we exploit the structure of these bad nonlinearities: Extracting the leading decay order, we realize that the potential and these bad nonlinearities have positivity, as in \eqref{EX2.1}, which can be incorporated into the principal operator without losing the invertibility. 
Even though this modified operator has an extra potential depending on the solution itself, we are able to obtain uniform (in all small $\xi$) resolvent estimates. This approach enables us to absorb the ``bad" nonlinear terms into the resolvent rather than estimating them using the ``canonical" resolvent. This is one of the main novelties of the present paper (see Propositions \ref{Pc2.1} and \ref{PZC1}).  \par
One drawback of this trick is that the fixed point does not still work, since if we attempt to obtain a difference estimate, the positivity structure in \eqref{EX2.1} breaks down. 
Hence, this idea can only be applied to obtain estimate for already constructed \textit{small} solutions. 
However, 
$\xi$ is not necessarily small a priori, which causes problem to control the nonlinearity. 
To compensate that, we construct \textit{another} solution $\tilde{f}_{r,\be}$ which is sufficiently close to $\te$, by discussing direct method in a small neighborhood of $\te$. 
In the direct method, we need to guarantee that the minimum is in the interior of the neighborhood; otherwise it does not necessary solve \eqref{E1.5}. 
A key observation is the
convexity estimate of $\tilde{E}_{r,0}$ around $\te$ (see Lemma \ref{L2.3.1}), which can be seen by restricting functions on $0\le f\le \te$. This can be viewed as a rigidity estimate for $\te$ in terms of $\tilde{E}_{r,0}$.
%
%
\subsubsection*{Shape analysis} 
The exponential decay of $f_{r,\be}$ follows from the fact that $g := \sqrt{\rho}f_{r,\be}$ satisfies 
$$
\frac{g''(\rho)}{g(\rho)} 
\xrightarrow{\rho\to\infty} 2\be^2, 
$$
which implies \eqref{EA0.3} by applying classical ODE argument as in \cite{BerLio83}. We remark that this argument works for any decaying solution as $\rho\to\infty$. \par
%
%
%
%
On the other hand, the proof of the monotonicity requires a couple of new ideas. 
For the exposition, we shall compare with known related studies. 
The monotonicity of the solutions to \eqref{E1.5} was obtained in two asymptotic regimes, namely when $k=0$, $h\gg 1$ in \cite{LiMel18}, and when $h=k\gg 1$ in \cite{GusWan21}. 
%
The basis of their argument is to observe 
$$
F(\rho):= \rho^2 (f')^2 - \sin^2 f,\qquad 
f=f_{r,\be}.
$$
This is an integral of the Euler-Lagrange ODE for $\tilde{D}[f]$, 
as it can naturally be seen as Newton's second law via change of variable: 
$$
\frac{d^2f}{dt^2} = \sin f\cos f,\qquad 
\rho = e^t.
$$
For solutions to \eqref{E1.5}, $F$ satisfies a structured ODE
\begin{equation}\label{EZ1}
\rd_\rho F = \rho \sin f \left( 
1-\cos f - \frac{2r\sin f}{\rho} + \be^2 (1+\cos f)
\right).
\end{equation}
In \cite{LiMel18} and \cite{GusWan21},  
it was possible to control the dynamics of \eqref{EZ1} thanks to the asymptotic regime of the parameters. However,  these estimates fail outside that range, namely when 
$\be\sim 0+$ (i.e. $h\sim -k$). 
Instead, we will give two new key ideas to examine, more closely, the dynamics of $F$.
First, we observe that the following identity holds:
$$
1+\cos f - \frac{\rho}{2r} \sin f = 
\frac{2\cos \frac f2 \sin \frac{\te -f}2 }{\sin \frac \te 2}.
$$
Then, the pointwise bound $0\le f \le \te$ of our solutions implies the positive-definiteness of this quantity, concluding that $F$ cannot take negative minimal value on where $f'<0$ (see Lemma \ref{L3.1} for the details). 
Second, 
when $r\le 1$, we observe that \eqref{E1.5} bears another sign-definite quantity 
\begin{equation}\label{EZ2}
\rd_\rho f - \frac 1\rho \sin f + r(1-\cos f) <0,\quad (0<\rho <\infty).
\end{equation}
If one supposes the existence of a point where $f'$ vanishes, then one can find a point where $\rd_\rho f - \frac 1\rho \sin f>0$ by further assuming $\be\le 1$. This leads to a contradiction to \eqref{EZ2} (see Lemmas \ref{L3.2} and \ref{L3.4}). 
%
%
%
\subsubsection*{Stability and instability} 
The framework to show the linear stability is the same as that in \cite{LiMel18}. Namely, we decompose the Hessian into each Fourier mode with respect to $\psi$, reducing the problem to one-dimensional ones. If $r\le \frac 12$, the monotonicity of the reduced Hessian in terms of the Fourier mode holds, thanks to which we only need to examine the zeroth and first Fourier modes. 
We remark that the sense of stability is weaker than that in \cite{LiMel18}, which comes from the technical reason that $B^{(0)}$ defined in \eqref{EX3.2} might have a nontrivial kernel, causing the failure of spectral gap of the Hessian even if we take orthogonal to the kernel from the translational invariance.\par 
For the instability, we make use of the fact that $\bh$ is unstable if $r>1$, as shown in \cite{IbrShi23}. By \eqref{EA0.2}, we are able to show that the Hessian of $E_{r,\be}$ at $\bn$ is very close to that of $E_{r,0}$ at $\bh$, and hence the instability is inherited for small $\be$. 

\subsubsection*{Comparison with the related results}
%
%
%
\begin{table}[h]
\begin{tabular}{|c||c|c|}
\hline
&Existence of minimizers & 
Difference estimate \\
\hline
%
%
\begin{tabular}{c}
Case\\ 
$\be \gg 1$
\end{tabular}
&
\begin{tabular}{c}
Lower bound $\inf E_\be \ge 0$;\\
Classical (Non-Restricted) 
Minimization;\\
Result similar to  
\cite{BerMurSim21,GusWan21}.
\end{tabular}
&
\begin{tabular}{c}
Reference energy is $D$;\\
Scope of \cite{BerMurSim21,GusWan21}.
\end{tabular}
\\
\hline
%
%
\begin{tabular}{c}
Case\\ 
$0<\be\lesssim$ 1
\end{tabular}
&
\begin{tabular}{c}
Lack of lower bound of $E_\be$;\\
The standard direct method fails;
\\
Use $f_{\be_0}$ ($\be_0\gg 1$) \\
to define and implement\\ 
a restricted minimization 
%
\end{tabular}
& 
\begin{tabular}{c}
Reference energy is $E_0$;\\
Lack of appropriate upper bound of $E_\be$;\\
Remedy: Use ODE \&\\
uniform resolvent estimates
\end{tabular}
\\
\hline
\end{tabular}
\end{table}


We conclude this section by summarizing the challenges of the problem and our proposed ideas.
Our parameter range splits naturally into two regimes: the asymptotic range ($\be\gg 1$) and the remaining range ($0<\be\lesssim 1$).
In the first case ($\be\gg 1$), the energy can be regarded as a perturbation of the Dirichlet energy, making the situation analogous to that in \cite{BerMurSim21, GusWan21}.
In contrast, the second regime presents new and significant difficulties. 
Regarding the construction of solutions, 
when $\beta\to0^+$, one cannot obtain a uniform lower bound $\inf \tilde{E}_{r,\be} >0$, so standard minimization methods fail.  
To address this, we introduce a constrained minimization problem with pointwise bounds $f_{\be_0}\le \bn\le \te$, where $\be_0$ is large and $f_{\be_0}$ is the solution from the first regime. This approach ensures both boundedness from below and the preservation of topology under weak limits. 
There, the fact that $f_{\be_0}$, $\te$ are sub-, super-solutions prevents the boundary from the boundary. 
%
Next, we discuss the difference estimate. In the case $\be\gg 1$, 
the question about the difference estimate from harmonic maps 
falls within the framework of
in \cite{BerMurSim21, GusWan21}, which we do not pursue in the present work. 
There, the rigidity of the harmonic maps \cite{BerMurSim21} or the perturbative analysis of the ODE \cite{GusWan21} played a key role. 
In contrast, 
when $\be\to 0+$, the rigidity argument as in \cite{BerMurSim21} does not work due to the lack of appropriate upper bound 
$\tilde{E}_{r,\be}[f_\be]\le (1+o_\be(1))\tilde{E}_{r,0}[f_\be]$. 
Although we follow the same spirit as in \cite{GusWan21}, 
the standard framework using the resolvent estimate is not sufficient due to the presence of uncontrollable nonlinearities. Therefore, we develop the uniform resolvent estimate with respect to the potential (coming from nonlinearities).

\subsection{Notations and organization of the paper}
Throughout the paper, 
$C$ is used to represent a universal constant which differs from line to line. 
The relation of two norms of functions  
$\nor{\cdot}{X}\sim \nor{\cdot}{Y}$ means that 
there exists $C>0$ such that 
$C^{-1}\nor{\cdot}{X}\le 
\nor{\cdot}{B}
\le C
\nor{\cdot}{A}
$. 
$\N$ denotes the set of all positive numbers. 
For $p\in [1,\infty]$ and $\Om\subset \R^n$ with $n\in\N$, $L^p(\Om)$ denotes the Lebesgue space on $\Om$ 
with exponent $p$, 
whose norm is denoted by $\nor{\cdot}{L^p(\Om)}$. 
For $f=f(\rho):(0,\infty)\to\R$ and $p\in[1,\infty]$, we also set 
$\nor{f}{L^p_{\rho d\rho}} = 
\left(\int_0^\infty |f|^p \rho d\rho\right)^{1/p}$, and define the space of $f$ with $\nor{f}{L^p_{\rho d^\rho}}$ by $L^p_{\rho d\rho}$. 
If there is no risk of confusion, we write $L^p_{\rho d\rho}$ as $L^p$ for simplicity. 
For a Hilbert space $X$, we denote the inner product on $X$ by $\inp{\cdot}{\cdot}_{X}$. 
For a Banach space $X$, 
let $B_X(R)$ be the ball in $X$ with radius $R$ whose center is at the origin of $X$.\par
The organization of the paper is as follows. 
Section \ref{S2} is devoted to the proof of the existence of solutions to \eqref{E1.5}. 
More specifically, 
we start with the case when $\be$ is sufficiently large (Section \ref{S2.1}), and then fill the remaining values of $\be$ (Section \ref{S2.2}). 
Next, when $\be$ is sufficiently small, we construct solutions with the difference estimate \eqref{EA0.2} in Section \ref{S2.3}. 
We give a supplementary argument concerning the pointwise bound of the above constructed solutions in Section \ref{S2.4}. 
In Section \ref{SX3}, 
we investigate the shape of the solutions. In detail, 
we show the exponential decay of the solutions in Section \ref{S3.1}, 
monotonicity in Section \ref{S3.2}, and the existence of the right derivative at the origin in Section \ref{S3.3}. 
In Section \ref{S3}, we turn to the analysis of stability. 
In Section \ref{S3.4}, we examine the stability of corresponding equivariant solutions to \eqref{E1.3}, while the instability in Section \ref{S3.5}. 
Section \ref{S4} is devoted to the proof of lemmas postponed in the main sections. In Section \ref{S4.1}, we show the regularity of weak solutions to \eqref{EF1.5}. 
In Section \ref{S4.2}, we give a proof of the uniform resolvent estimate.

\section{Existence: variational approach}
\label{S2}

In what follows, we drop the tildes from the functionals for $f$ for simplicity; namely, we write $\tilde{E}_{r,\be}[f]$ as $E_{r,\be}[f]$ for instance. 
Also, we abbreviate the index $r$ in $E_{r,\be}$ when there is no risk of confusion. 

\subsection{The regime of large values of $\be$}\label{S2.1}
We first address the 
construction of the solution to \eqref{E1.5} in the case when 
$\be$ is sufficiently large:
\begin{thm}\label{P1}
For any $r>0$, let $\be > 2r$, then 
the minimization problem
$$
\inf_{f\in\boM_*} E_\be[f]
$$
with
\begin{equation}\label{EFeb4}
\boM_* := \left\{
f:(0,\infty)\to\R\ \left| \ 
f\in \boM,\ 0\le f(\rho) \le \te(\rho)\quad(\forall \rho\in (0,\infty))
\right.
\right\}
\end{equation}
attains its minimum $f_* \in \boM_*$. Moreover, $f_*$ solves \eqref{E1.5}.
\end{thm}
The advantage when $\be>2r$ is that the energy becomes positive definite in $\boM_*$:
\begin{lem}\label{L1}
For $f\in L^2$ with $\rd_\rho f$, $\frac{\sin f}\rho \in L^2$, we have
\begin{equation}\label{2.1}
E_\be[f] \ge \left(1 - \frac{4r^2}{\be^2} \right) D[f] +  V_-[f] + \frac 12\be^2 V_+[f].
\end{equation}
\end{lem}
\begin{proof}
In fact, the Schwarz inequality (see also \cite{DorMel17}) yields
\begin{equation}\label{E2.2}
H^2[f] \le 8 D[f] V_-[f]  
\le 8 D[f] V_+[f],
\end{equation}
which implies
\begin{align*}
E_\be[f] 
&\ge 
D[f] - 2\sqrt 2 r \sqrt{D[f]} \sqrt{V_+[f]} 
+ V_-[f] + \be^2 V_+[f] \\
&\ge  
D[f] -2\sqrt 2 r \frac 12 \left( \frac {2\sqrt 2 r}{\be^2} D[f] 
+ \frac{\be^2}{2 \sqrt 2 r} V_+[f] \right) + V_-[f]+ \be^2 V_+[f]  \\
&= 
\left(1 - \frac{4r^2}{\be^2} \right) D[f]  + V_-[f] +\frac 12 \be^2 V_+[f] .
\end{align*}
\end{proof}
The positivity of the energy will enable us to exclude the possibility of topological change as a minimizing sequence converges to its weak limit.

\subsubsection{Direct method}
Set $E_{\be,*}= 
\inf_{f\in\boM_*} E_\be[f]$. 
Let $\{f_j\}_{j=1}^\infty$ be a minimizing sequence of $E_\be$ in $\boM_*$:
$$
\{f_j\}_{j=1}^\infty \subset \boM_*,\qquad 
\lim_{j\to\infty} E_\be[f_j] = E_{\be,*}.
$$
As a first step, we show that $f_j$ converges weakly to some $f_*\in L^2$ with $\rd_\rho f$, $\frac{\sin f_*}\rho \in L^2$, $0\le f_* \le \te$ and $E_\be[f_*] \le E_{\be,*}$. ($f_*$ will satisfy all the conditions for being in $\boM$ except the boundary conditions.) 
First of all, we observe that $\{f_j\}$, $\{\rd_\rho f_j\}$, and $\{\frac{\sin f_j}{\rho}\}$ are bounded in $L^2$. 
Indeed, 
by \eqref{2.1} and $\be> 2r$, 
$D[f_j]$ and $V_+[f_j]$ are bounded in $j$. Hence, 
we have
\begin{gather*}
\frac 12 \nor{\rd_\rho f_j}{L^2}^2
+ 
\frac 12 \nor{\frac{\sin f_j}{\rho}}{L^2}^2
\le D[f_j] \le C,\\
\int_0^\infty |f_j|^2 \rho d\rho 
\le 
C\int_0^\infty (1-\cos f_j) \rho d\rho 
\le C Z[f_j] 
\le C V_+[f_j] \le C
\end{gather*}
with constant $C$ independent of $j$. 
In particular, $f_j(\rho)$ converges to $0$ as $\rho \to\infty$ uniformly in $j$, thanks to inequality
\begin{align}\label{2.7}
f^2(\rho) = -2 \int_\rho^\infty f(s) \rd_\rho f (s) ds \le \frac 2 \rho \nor{f}{L^2} \nor{\rd_\rho f}{L^2}.  
\end{align}
Moreover, 
there exist a subsequence of $\{f_j\}$, which is also denoted by $\{ f_j\}$, and $f_*\in L^2$ with $\rd_\rho f_*, 
\frac{\sin f_*}{\rho} \in L^2$ such that
\begin{equation}\label{2.3}
f_j\xrightharpoonup{j\to\infty} f_*,\quad 
\rd_\rho f_j\xrightharpoonup{j\to\infty} \rd_\rho f_*,\quad 
\frac{\sin f_j}{\rho} \xrightharpoonup{j\to\infty}
\frac{\sin f_*}{\rho} 
\quad 
\text{weakly in }
L^2.
\end{equation}
Interpreting $f_j$ as radial functions in $H^1(\R^2)$, we apply the compact embedding $H^1(B_R(0)) \subset L^p(B_R(0))$ for any $p\in [1,\infty)$ and for any $R>0$.  
Then we have
\begin{equation}\label{2.4}
f_j \xrightarrow{j\to\infty} f_* \qquad
\text{strongly in } L^p ((0,R); \rho d\rho) \text{ for all } R>0, \ p\in [1,\infty).
\end{equation}
Also, 
using the compact embedding for one-dimensional functions $H^1 (\frac 1R , R) \subset C(\frac 1R, R)$, we have
\begin{equation}\label{2.5}
f_j \xrightarrow{j\to\infty} f_* \qquad
\text{uniformly in } \left(\frac 1R, R\right),\qquad \text{ for all } R>0.
\end{equation}
In particular, we have $0\le f_* \le \te$ on $(0,\infty)$. 
Therefore, by weak lower semi-continuity of norms and Fatou's lemma, we obtain
\begin{gather}
\label{2.8}
\nor{\rd_\rho f}{L^2}^2 \le \liminf_{j\to\infty} \nor{\rd_\rho f_j}{L^2}^2;
\\
\label{2.9}
\nor{\frac{\sin f}{\rho}}{L^2}^2 \le \liminf_{j\to\infty} \nor{\frac{\sin f_j}{\rho}}{L^2}^2;
\\
\label{2.10}
V_-[f] \le \liminf_{j\to\infty} V_- [f_j],\qquad 
V_+[f] \le \liminf_{j\to\infty} V_+[f_j].
\end{gather}
Moreover, 
we claim that
\begin{equation}\label{2.11}
\lim_{j\to\infty} H[f_j] = H[f_*].
\end{equation}
Indeed, the family of integrand $\left\{\left(f'_j -\frac{\sin f_J}{\rho}\right) (1-\cos f_j)\right\}_{j\ge 1}$ is tight since \eqref{2.7} implies
\begin{align*}
\left|
\int_R^\infty 
\left(f'_j -\frac{\sin f_J}{\rho}\right) (1-\cos f_j) \rho d\rho \right|
&\le 
\int_R^\infty 
\left(
|f'_j| + \left| \frac{\sin f_j}{\rho} \right|
\right) |f_j|^2 \rho d\rho \\
&\le 
\frac C{\sqrt R} 
\left( \nor{f_j'}{L^2} + \nor{\frac{\sin f_j}{\rho}}{L^2} \right) \nor{f_j}{L^2} \\
&\xrightarrow{R\to\infty} 0 \quad \text{uniformly in } j.
\end{align*}
On the other hand, 
on $(0,R)$, we have
\begin{align*}
&\int_0^R 
\left( f'_j -\frac{\sin f_j}{\rho} \right)
(1-\cos f_j) \rho d\rho - 
\int_0^R 
\left( f'_* -\frac{\sin f_*}{\rho} \right)
(1-\cos f_*) \rho d\rho \\
&= 
\int_0^R \left( f'_j -\frac{\sin f_j}{\rho} \right)
(\cos f_* -\cos f_j) \rho d\rho \\
&\hspace{40pt}
+ 
\int_0^R \left( f'_j - f_*' - 
\frac{\sin f_j - \sin f_*}{\rho} \right)
(1-\cos f_*) \rho d\rho \\
&\xrightarrow{j\to\infty} 0,
\end{align*}
by using \eqref{2.3}. 
Therefore, \eqref{2.8}, \eqref{2.9}, \eqref{2.10}, \eqref{2.11} imply
\begin{equation}\label{2.12}
E_\be [f_*] \le \liminf_{j\to\infty} E_\be [f_j] = E_{\be,*}. 
\end{equation}

\subsubsection{Topological change scenario}
By \eqref{2.12}, 
to show that $f_*$ is a minimizer in $\boM_*$, it suffices to check the boundary conditions: 
\begin{align}\label{2.13}
f_* (0)=\pi,&& f_*(\infty)=0.
\end{align}
Note that $f_*(\infty)=0$ immediately follows from \eqref{2.7}. 
On the other hand, the conditions that $\rd_\rho f, 
\frac{\sin f}\rho \in L^2_{\rho d\rho}$ imply that $f_*(0)=m\pi$ for 
some $m\in \Z$ (see \cite{GusWan21} for detailed proof). 
Hence, the pointwise bound $0\le f_* \le \te$ yields 
\begin{align}\label{2.135}
f_*(0) =0 \qquad \text{or} \qquad  f_*(0)=\pi.
\end{align}
To show that $f_*\in\boM_*$, it suffices to exclude the former scenario.
\begin{lem}[Topological change costs at least energy  $2$]\label{L2.1} 
If $f_*(0)=0$, then
$$
\liminf_{j\to\infty} D[f_j] \ge 2 +D[f_*].
$$
\end{lem}
See Lemma 11 in \cite{GusWan21} for the detailed proof. Let us briefly mention the mechanism underlying Lemma \ref{L2.1}. 
If the boundary condition is changed under the weak limit, at least one rescaling family of harmonic maps should be bubbled out as $j\to\infty$. 
Accordingly, the energy of $f_j$ is decoupled into the Dirichlet energy of each bubble, the energy for $f_*$, and an error $o_j(1)$. 
Since the least Dirichlet energy of harmonic maps is $2$, we arrive at the conclusion in Lemma \ref{L2.1}.
\begin{lem}[Estimate of $E_{\be,*}$]
\label{L2.2}
$$
E_{\be,*} <2.
$$
\end{lem}
\begin{proof}
Let $\chi\in C_0^\infty (\R)$ be a cut-off function satisfying
$$
0\le \chi \le 1 \text{ on } \R,\qquad 
\chi(\rho) = 
\left\{
\begin{aligned}
&1 \quad &\text{if} \quad \rho \le 1,\\
&0 \quad &\text{if} \quad \rho \ge 2.
\end{aligned}
\right.
$$
Then for $\la>1$ and $R>0$, define
$$
\te_{\la,R} (\rho) 
= \te (\la\rho) \chi \left(\frac{\la\rho}R \right). 
$$
Obviously, $\te_{\la,R}\in \boM$. Moreover, 
since $\te$ is strictly decreasing,  
we have $0\le \te_{\la,R} \le \te$ if $\la>1$, and thus $\te_{\la,R}\in \boM_*$. 
Moreover, direct computation yields
\begin{gather*}
D[\te_{\la,R}] = D[\te] + O_{R\to\infty}\left(\frac 1{R^2}\right) = 2 + O_{R\to\infty}\left(\frac 1{R^2}\right) ,
\\
H[\te_{\la,R}] = \frac 1\la \left(H[\te] + O_{R\to\infty}\left(\frac 1{R^2}\right) \right)
= \frac 1\la \left(-8r+ O_{R\to\infty}\left(\frac 1{R^2}\right) \right),
\\
V_-[\te_{\la,R}] = \frac 1{\la^2} \left(V_-[\te] + O_{R\to\infty}\left(\frac 1{R^2}\right)\right) 
= \frac 1{\la^2} \left(16r^2 + O_{R\to\infty}\left(\frac 1{R^2}\right)\right) ,
\\
V_+[\te_{\la,R}] = 
\frac{1}{\la^2} O_{R\to\infty}(\log R).
\end{gather*}
Thus, for $\la>1$, we have
$$
E_\be [\te_{\la,R}] = 
2 - \frac{8r^2}{\la} + \frac{16r^2}{\la^2} 
+ \frac{\be^2}{\la^2} O_{R\to\infty}(\log R)
+ O_{R\to\infty}\left(\frac 1{R^2}\right).
$$
If we take $R= \la$, then
$$
E_\be [\te_{\la,R}] = 2 - \frac{8r^2}\la + O_{\la\to\infty}\left(\frac{1}{\la^2}\log \la\right).
$$
Hence, if we take $\la$ sufficiently large, we have $E_\be [\te_{\la,R}]<2$.
\end{proof}
Now we show $f_*(0)=\pi$, which concludes that $f_*$ is a minimizer of $E_\be$ in $\boM_*$. 
Suppose to the contrary, then 
by \eqref{2.135}, 
we may suppose $f_*(0)=0$. 
By Lemma \ref{L2.2}, there is $\eps>0$ such that $E_{\be,*} < 2-\eps$. 
Then for large $j$, we have
$$
E_\be[f_j] < 2-\eps.
$$
Taking $j\to\infty$ and using Lemma \ref{L2.1}, we have
\begin{align*}
2-\eps &\ge \liminf_{j\to\infty} E_\be [f_j] \\
&\ge
\liminf_{j\to\infty} D[f_j]
+r\liminf_{j\to\infty} H[f_j]
+\liminf_{j\to\infty} V_-[f_j]
+\be^2\liminf_{j\to\infty} V_+[f_j]\\
&\ge 2+ D[f_*] + rH[f_*] + V_-[f_*] + \be^2 V_+[f_*] = 2+ E_\be[f_*] ,
\end{align*}
which implies
$$
E_\be [f_*] < -\eps <0.
$$
If $\be>2r$, however, then this leads to a contradiction since Lemma  \ref{L1} implies $E_\be[f_*]\ge 0$. Hence we conclude $f_*(0)=\pi$.

\subsubsection{Euler-Langrange equation}
Next, we show the following:
\begin{lem}
\label{L2.3}
Let $f_*$ be a minimizer of $E_\be$ in $\boM_*$. 
Then $f_*$ solves \eqref{E1.5}.
\end{lem}
The only issue is that $f_*$ might touch the boundary of $\boM_*$, which could prevent us to take arbitrary variation from $f_*$. However, this obstacle can be overcome by using the fact that $0$ and $\te$ are sub-solution and super-solution of \eqref{E1.5}, respectively. This trick is now classical, known as the sub/supersolution method, see for example \cite{Struwe}. 
We reproduce the proof here, since a little subtlety exists since $\te$ is not in $\boM$.
\begin{proof}
Let $\vph \in C_0^\infty (0,\infty)$, and $\eps>0$. Define
$$
g_\eps := \min \{ \te ,  \max \{ 0, f_* +\eps \vph\} \} = 
f_* + \eps \vph - \vph^\eps + \vph_\eps,
$$
where
\begin{align*}
&\vph^\eps = \max\{0,f_* + \eps \vph -\te \} \ge 0,\\
&\vph_\eps = -\min\{0,f_* + \eps \vph \} \ge 0.
\end{align*}
Note that both are weakly differentiable on $(0,\infty)$, and $\supp \vph^\eps$, $\supp \vph_\eps
\subset \supp \vph$ because of the pointwise bound $0\le f_* \le \te$. 
By definition, we have $0\le g_\eps \le \te$. 
Therefore, we have
$
(1-t) f_* + t g_\eps \in \boM_*$ for all $0\le t\le 1$. 
Since $f_*$ is a minimizer on $\boM_*$, it follows that
\begin{align*}
0 &\le \lim_{t\to 0+} \frac{E_\be[(1-t) f_* + t g_\eps] - E_\be [f_*]}{t} \\
&= 
\inp{E_\be'[f_*]}{ g_\eps -f_*}_{L^2} 
= \eps \inp{E_\be'[f_*]}{ \vph}_{L^2} - \inp{E_\be'[f_*]}{ \vph^\eps}_{L^2} + \inp{E_\be'[f_*]}{ \vph_\eps}_{L^2}
\end{align*}
with 
$$
E_\be'[f] = -\rd_\rho^2f - \frac 1\rho  \rd_\rho f + \frac 1{\rho^2} 
\sin f \cos f
- 2r\frac{\sin^2 f}\rho +
\sin f(1-\cos f) + \be^2 \sin f(1+\cos f)
$$
for $f\in \boM$, where the derivative is in distributional sense. 
This implies
\begin{equation}\label{eb1}
\inp{E_\be'[f_*]}{\vph}_{L^2} \ge \frac 1\eps \left( 
\inp{E_\be'[f_*]}{ \vph^\eps} - \inp{E_\be'[f_*]}{ \vph_\eps}_{L^2}
\right).
\end{equation}
Here we observe that 
even though $\te\not\in \boM$, the coupling $\inp{E_\be'[\te]}{\vph^\eps}_{L^2}$ is well-defined 
since $\supp \vph^\eps$ is compact. Moreover,
\begin{equation}\label{2.16}
\inp{E_\be'[\te]}{\vph^\eps}_{L^2} 
= \int_0^\infty \be^2 \sin \te (1+\cos \te)\vph^\eps \rho d\rho
\ge 0
\end{equation}
where we used the positivity of $\vph^\eps$. 
Hence
\begin{align*}
&\inp{E_\be'[f_*]}{ \vph^\eps}_{L^2} = \inp{E_\be'[f_*] - E_\be'[\te]}{ \vph^\eps}_{L^2} + 
\inp{E_\be'[\te]}{\vph^\eps}_{L^2} 
\ge \inp{E_\be'[f_*] - E_\be'[\te]}{ \vph^\eps}_{L^2} \\
&= \inp{\rd_\rho (f_*-\te)}{\rd_\rho \vph^\eps }_{L^2}\\
&\quad +  \inp{\frac{\sin f_*\cos f_*}{\rho^2} - 2r
\frac{\sin^2 f_*}{\rho} +\sin f_* (1-\cos f_*) + \be^2 \sin f_* (1+\cos f_*)
 }{ \vph^\eps }_{L^2} \\
&\quad - 
\inp{\frac{\sin \te\cos \te}{\rho^2} - 2r
\frac{\sin^2 \te}{\rho} +\sin \te (1-\cos \te) + \be^2 \sin \te (1+\cos \te)
 }{ \vph^\eps }_{L^2}. 
\end{align*}
Here, we define
\[
\Om^\eps \coloneqq 
\{\rho\in (0,\infty) \ |\ 
f_*<\te \le f_*+\eps \vph
\}.
\]
We note that $\supp \vph^\eps \subset \Om^\eps \subset \supp \vph$, and that $|\Om^\eps|\to 0$ as $\eps\to 0$ since $\cap_{\eps>0} \Om^\eps =\emptyset$. Moreover, on $\Om^\eps$, we have
\[
\vph^\eps \le |f_*-\te| +\eps|\vph| \le 2\eps |\vph|.
\]
%
Then we obtain
\begin{align*}
\inp{\rd_\rho (f_*-\te)}{\rd_\rho \vph^\eps }_{L^2} 
&= \int_{\Om^\eps} \rd_\rho (f_*-\te) \rd_\rho (f_*+\eps \vph - \te) \rho d\rho \\
&= \int_{\Om^\eps} |\rd_\rho (f_*-\te) |^2 \rho d\rho 
+\eps \int_{\Om^\eps} \rd_\rho (f_*-\te) \rd_\rho \vph \rho d\rho \\
&\ge \eps \int_{\Om^\eps} \rd_\rho (f_*-\te) \rd_\rho \vph \rho d\rho = o(\eps). 
\end{align*}
On the other hand, the compactness of the support gives
\begin{align*}
& \left| \inp{\frac{\sin f_*\cos f_*}{\rho^2} - 2r
\frac{\sin^2 f_*}{\rho} +\sin f_* (1-\cos f_*) + \be^2 \sin f_* (1+\cos f_*)
 }{ \vph^\eps }_{L^2} \right| \\
& + \left| \inp{\frac{\sin \te\cos \te}{\rho^2} - 2r
\frac{\sin^2 \te}{\rho} +\sin \te (1-\cos \te) + \be^2 \sin \te (1+ \cos \te)
 }{ \vph^\eps }_{L^2}  \right| \\
&\le C \eps \int_{\Om^\eps} |\vph(\rho)| \rho d\rho
= o(\eps).
\end{align*}
Therefore, we have
$$
\frac 1\eps \inp{E_\be'[f_*]}{ \vph^\eps }_{L^2} \ge o_{\eps\to 0} (1).
$$
Using the same argument as above, we have
$$
-\frac 1\eps \inp{E_\be'[f_*]}{ \vph_\eps }_{L^2} \ge o_{\eps\to 0} (1).
$$
Hence, by \eqref{eb1}, we have
$$
\inp{E_\be'[f_*]}{\vph}_{L^2} \ge o_{\eps\to 0}(1).
$$
By taking $\eps\to 0$, we have
$$
\inp{E_\be'[f_*]}{\vph}_{L^2}  \ge 0.
$$
If we repeat the above argument by replacing $\vph$ by $-\vph$, we obtain
$$
\inp{E_\be'[f_*]}{\vph}_{L^2}  \le 0.
$$
Hence we have 
$$
\inp{E_\be'[f_*]}{\vph}_{L^2}  = 0.
$$
Since $\vph\in C_0^\infty$ is arbitrary, we conclude that $f_*$ solves $E_\be'[f_*]=0$.
\end{proof}


\subsection{The regime of the remaining values of $\beta$}
\label{S2.2}

Fix $\be_0 >4r$, and let $f_0= f_{\be_0}$ be the solution to \eqref{E1.5} as in Theorem \ref{P1}. 
Then we claim the following:
%
%
%
\begin{thm}\label{P2}
For $r>0$ and $0<\be<\be_0(r)$, 
the minimization problem
$$
\inf_{f\in \boM_*} E_\be[f]
$$
with
$$
\boM_* := \left\{
f:(0,\infty)\to\R\ \left| \ 
f\in \boM,\  f_0(\rho) \le f(\rho) \le \te(\rho)\ \text{on }  (0,\infty)
\right.
\right\}
$$
attains its minimum $f_* \in \boM_*$. Moreover, $f_*$ solves \eqref{E1.5}.
\end{thm}
%
%
%
\begin{proof}
First of all, we observe that
$$
E_{\be,*} := \inf_{\boM_*} E_\be[f] >-\infty.
$$
Indeed, \eqref{E2.2} and $0\le f\le \te$ give
\begin{align*}
E_\be[f] 
&= D[f] + rH[f] + V_-[f] + \be^2 V_+[f] \\
&\ge 
\frac 12 D[f] + (1-4r^2) V_-[f] + \be^2 V_+[f] \\
&\ge \frac 12 D[f] + \min\{0, (1-4r^2) V_-[\te] \} + \be^2 V_+[f] > -\infty.
\end{align*}
Now let $\{f_j\}_{j=1}^\infty$ be a minimizing sequence of $E_\be$:
$$
\{f_j\}_{j=1}^\infty \subset \boM_*,\qquad 
\lim_{j\to\infty} E_\be[f_j] = E_*.
$$
Then the same argument as in Section \ref{S2.1} implies that $f_j$ has a weak limit $f_*$ in $L^2$ if we take subsequence if necessary, and that $f_*$ satisfies
\begin{itemize}
\item $f_*, \rd_\rho f_*, \frac{\sin f_*}{\rho} \in L^2$;
\item $f_j \to f_*$ as $j\to\infty$ uniformly in $(\frac 1R,R)$ for all $R>0$; 
\item $E_{\be} [f_*] \le \liminf_{j\to\infty} E_\be [f_j] = E_{\be,*}$.
\end{itemize}
Here, 
since $f_0 \le f_j \le \te$ for all $j$, the pointwise convergence implies that 
$
f_0 (\rho) \le f_* (\rho) \le \te (\rho)
$ on $(0,\infty)$. 
Since $f_0 (0) = \te(0)=\pi$ and $f_0(\infty)=\te(\infty)=0$, we necessarily have $f_*(0)=\pi$, $f_*(\infty)=0$. Therefore, $f_*\in \boM_*$, and hence $f_*$ is a minimizer of $E_\be$ in $\boM_*$. \par
That $f_*$ solves \eqref{E1.5} follows similarly to the proof of Lemma \ref{L2.3}, by modifying the definition of $\vph_\eps$ by
$$
\vph_\eps = -\min\{0, f_* +\eps \vph -f_0 \},
$$
and by using the fact that $f_0=f_{\be_0}$ is a subsolution to \eqref{E1.5} when $0<\be<\be_0$. Hence the proof is complete.
\end{proof}


\subsection{The regime of small values of $\be$}\label{S2.3}
The main theorem of this section is the following:
%
%
%
\begin{thm}\label{P3}(Existence of slightly deviated solutions)
For any $r>0$, there exists $\eps_0 = \eps_0(r)>0$, $\be_0 = \be_0(r)>0$ such that if $0<\be<\be_0$, then 
the minimization problem
$$
\inf_{f\in \boM_{\eps_0}} E_\be[f]
$$
with
$$
\boM_{\eps_0} := \left\{
f:(0,\infty)\to\R\ \left| \ 
\begin{aligned}
&f\in \boM,\ 0\le f(\rho) \le \te(\rho)\quad(\forall \rho\in (0,\infty)), \\
&\nor{f-\te}{X} \le \eps_0
\end{aligned}
\right.
\right\}
$$
attains its minimum at $\tilde{f}_\be \in \boM_{\eps_0}$. Moreover, $\tilde{f}_\be$ solves \eqref{E1.5}. 
Furthermore, 
for any $0<\del<1$, there exists $C=C_\del>0$ such that
\begin{equation}\label{EY1}
\nor{\tilde{f}_\be-\te}{X} < C\be^{1-\del}.
\end{equation}
\end{thm}
%
%
%
Using the same argument as in Section \ref{S2.1}, 
we can verify that $E_\be$ still attains its minimum 
in $\boM_{\eps_0}$. 
However, the problem is that due to the additional restriction in $\boM_{\eps_0}$, a minimizer may touch the boundary of $\boM_{\eps_0}$ and therefore does not have to satisfy \eqref{E1.5}. 
We rule out this possibility by 
observing that $E_0$ enjoys a convexity property around $\te$ via the Taylor expansion. 
To obtain \eqref{EY1}, we  linearize the Euler-Lagrange equation around $\te$. The desired 
difference estimate 
then follows 
from  
improved resolvent estimates given in Section \ref{S4.2} (see Section \ref{S2.3.3}).

\subsubsection{Verification of the Euler-Langrange equation}

The key lemma 
to verify that $\tilde{f}_\be$ solves \eqref{E1.5} is the following.
\begin{lem}[Key lemma]\label{L2.3.1}
There exist $C_0 = C_0(r)>0$, $\eps_1 = \eps_1 (r)>0$ such that the following holds: 
If $f\in \boM$ satisfies $0\le f\le \te$ and $\nor{f-\te}{X} \le \eps_0$, 
then
\begin{equation}\label{EX2.3}
E_0[f] - E_0[\te] \ge C_0 \nor{f-\te}{X}^2.
\end{equation}
\end{lem}
We postpone the proof of Lemma \ref{L2.3.1} and first show \eqref{E1.5}. 
\begin{lem}[Bound of minimum]\label{L2.3.2}
Let $C_0$, $\eps_1$ be as in Lemma \ref{L2.3.1}, and let $0<\eps_0\le \eps_1$. Then, 
there exists $0<\be_0 = \be_0(r,\eps_0) \ll 1$ such that if $0<\be<\be_0$, we have
$$
\inf_{f\in \boM_{\eps_0}} E_{\be}[f] < E_0[\te] + \frac 12 C_0 \eps_0^2.
$$
\end{lem}
\begin{proof}
For $R>0$, let
$$
\te_R =  \chi_R \te,\qquad \chi_R (\rho)= \chi \left(\frac{\rho}{R}\right),\qquad 
\chi (\rho) =
\left\{
\begin{aligned}
&1 \quad &\text{if} \quad \rho \le 1,\\
&0 \quad &\text{if} \quad \rho \ge 2.
\end{aligned}
\right.
$$
Then 
$\nor{\te_R-\te}{X} = O(\frac 1{R^2})$, and
\begin{align*}
D[\te_R] = D[\te] + O(\frac 1{R^2})
,&&
H[\te_R] = H[\te] + O(\frac 1{R^2}),
\end{align*}
\begin{align*}
V_-[\te_R] = V_-[\te] + O(\frac 1{R^2}),
&&
V_+[\te_R] = 
O(\log(1+R^2)).
\end{align*}
We take $R$ sufficiently large such that $\nor{\te_R-\te}{X} \le \eps_0$, 
which 
in particular implies $\te_R\in \boM_{\eps_0}$. 
On the other hand, we get
$$
E_\be [\te_R] = E_0[\te] + O\left(\frac 1{R^2}\right) + \be^2 N_R 
,
\qquad 
N_R 
= O(\log(1+R^2)).
$$
If we retake $R$ sufficiently large, we have
$$
O(\frac 1{R^2}) < \frac 14 C_0 \eps_0^2.
$$
Once fixing  such $R$, we take $\be_0 = \be_0(r)$ sufficiently small, such that 
$\be^2 N_R <\frac 14 C_0 \eps_0^2$ 
if $0<\be<\be_0$. 
Then we obtain
$$
\inf_{f\in \boM_{\eps_0} } E_{\be}[f] \le E_{\be} [\te_R] < E_0[\te] + \frac 12 C_0 \eps_0^2.
$$
\end{proof}


\begin{proof}[Proof of \eqref{E1.5} for $\tilde{f}_\be$.]  
Let $\eps_0$, $\eps_1$, $\be_0$, $C_0$ be the numbers as in Lemmas \ref{L2.3.1} and \ref{L2.3.2}.  
Then, by Lemma \ref{L2.3.2}, we have
\begin{equation}\label{EX2.2}
E_0[\tif] \le 
E_\be[\tif] < 
E_0[\te] + \frac 12 C_0\eps_0^2
\end{equation}
for $0<\be<\be_0$. 
If $\nor{\tif -\te}{X}=\eps_0$, then \eqref{EX2.2} contradicts \eqref{EX2.3}. 
Hence we have $\nor{\tif-\te}{X}<\eps_0$. 
In particular, 
any small perturbation of $\tif$ in $X$ is included in this ball, and hence 
we can apply the same argument as in the proof of Lemma \ref{L2.3}. 
Therefore, we conclude that $\tif$ satisfies \eqref{E1.5}.
\end{proof}

\subsubsection{Proof of Lemma \ref{L2.3.1}}

Now we prove Lemma \ref{L2.3.1}. 
We start by showing the following coercivity estimate. 
\begin{lem}\label{L2.3.3}
Let $f:(0,\infty)\to \R$ be a function satisfying $0\le f\le \te$. Then, $\xi:= f-\te$ satisfies 
$$
E_0[f] -E_0[\te] \ge 
\frac 12 \nor{F \xi}{L^2}^2 
+ \int_0^\infty \frac{1}{18} \frac{4r^2}{\rho^2 + 4r^2} \xi^2 \rho d\rho  
- C\nor{\xi}{X}^3, 
$$
where $C>0$ is a constant  independent of $f$, and
\begin{equation}\label{EZB1}
F= \rd_\rho + \frac{\cos \te}\rho.
\end{equation}
\end{lem}
\begin{proof}
We consider the Taylor expansion of $E_0$ around $\te$. 
If we compute each functional in $E_0$, we obtain
\begin{align*}
D[f]
= \int_0^\infty 
&\left[\frac 12 (\te')^2 +  \te' \xi' + \frac 12 (\xi')^2 \right.\\
&
+\frac 1{\rho^2} \left(
\frac 12 \sin^2 \te + \sin \te \cos \te \xi 
+ \frac 12  (\cos^2\te - \sin^2\te) \xi^2
- \frac 23 \sin \te \cos \te \xi^3
\right.\\
& \qquad 
\left.\left.
+ \frac 16 (\sin^2\te  
-\cos^2\te
) \xi^4 
+ R_1 \xi^5
\right)
\right] \rho d\rho,
\end{align*}
where 
$R_1$ is the fifth coefficient of the Taylor expansion in $\xi$ of $\frac 12 \sin^2 (\te+\xi)$ around $\te$. Notice that 
$R_1=R_1(\rho)$ is 
bounded on $(0,\infty)$ whose bound is independent of $\xi$. Similarly,
%
\begin{align*}
H[f]
= 
\int_0^\infty 
&
\left[ (\te' - \frac 1\rho \sin\te )(1- \cos \te)\right. \\
& + (1- \cos \te)(\xi' - \frac{1}\rho \cos\te \xi ) 
+ (\te' - \frac 1\rho \sin\te ) \sin \te \xi \\
& + \frac 1{2\rho} \sin \te (1- \cos \te) \xi^2  
+ (\xi' - \frac{1}\rho \cos\te \xi ) \sin \te\xi 
+ \frac 12 (\te' - \frac 1\rho \sin\te )  \cos\te \xi^2\\
&
+\frac 1{6\rho} \cos \te  (1-\cos\te) \xi^3
+ \frac 1{2\rho} \sin^2 \te \xi^3 
+ (\xi' - \frac{1}\rho \cos\te \xi ) \frac 12 \cos\te \xi^2 
\\
&-\frac 16 \sin \te
 (\te' - \frac 1\rho \sin\te )  \xi^3 
- \frac 1{24\rho} \sin \te  (1- \cos \te) \xi^4
+ \frac 5{12\rho}  \sin \te \cos \te \xi^4 \\
&
\left.
+ (\xi' - \frac{1}\rho \cos\te \xi ) (-\frac 16 \sin\te \xi^3) 
+ (\te' - \frac 1\rho \sin\te ) (-\cos \te \frac{\xi^4}{24}) 
+ \frac 1{\rho} R_2 \xi^5  + R_3 \xi' \xi^4\right] \rho d\rho,
\end{align*}
\begin{align*}
V_-[f] = 
\int_0^\infty 
&
\left[\frac 12(1- \cos \te)^2  + \sin \te (1-\cos \te) \xi 
+\frac 12 [(1-\cos\te)\cos\te + \sin^2\te] \xi^2 \right. \\
&
+ \left[-\frac 16 (1-\cos \te) \sin\te + \frac 12 \sin \te\cos\te \right] \xi^3 \\
&
\left.+ \left[
-\frac 1{24} \cos\te (1-\cos \te) - \frac 16 \sin^2\te 
+ \frac 18 \cos^2\te\right] \xi^4 + R_4 \xi^5\right] \rho d\rho.
\end{align*}
Here $R_j$'s are the coefficients  defined 
by the Taylor expansion, satisfying 
similar properties to $R_1$. 
Based on this, we expand $E_0$ as
$$
E_0[f] -E_0[\te] = 
\int_0^\infty 
\sum_{j=1}^5 \boK_j (\xi,\xi') \rho d\rho, 
$$
where $\boK_j$ is the collection of $j$-order terms with respect to $\xi$ and $\xi'$. 
We note that $\int_0^\infty \boK_1 \rho d\rho =0$ since $E'_0[\te]=0$. 
For $\boK_2$, we have
\begin{align*}
\boK_2
&= \frac 12 (\xi')^2 
+ \frac 1{2\rho^2} (\cos^2 \te -\sin^2 \te)\xi^2 
+ r\left[
\frac 1{2\rho} \sin\te (1-\cos\te) \xi^2 \right] \\
&\quad 
+r\left[
(\xi' - \frac 1\rho \cos\te \xi) \sin\te \xi 
+ \frac 12 (\te' - \frac 1\rho \sin\te) \cos\te \xi^2
\right] 
+\frac 12 [(1-\cos\te)\cos\te + \sin^2\te] \xi^2 \rho d\rho \\
&=
\frac 12  (\xi')^2 + \frac{1}{2\rho^2} (1 -2\sin^2\te) \xi^2
+  g(\rho) ,
\end{align*}
where
\begin{align*}
g(\rho) &= r \sin\te \xi' \xi \\
&\quad +
\left[
\frac r{2\rho} \sin\te (1-\cos\te) -\frac r\rho \sin\te\cos\te 
+\frac r2 \te' \cos \te 
+\frac 12 (1-\cos\te)\cos\te + \frac 12 \sin^2\te \right] \xi^2.
\end{align*}
For the first term of $g$, integration by part, together with \eqref{1.135}, gives
\begin{align*}
\int_0^\infty r\sin\te \xi'\xi \rho d\rho 
= -\int_0^\infty \frac r2 \frac {\sin\te}\rho (1 - \cos\te )  \xi^2 \rho d\rho.
\end{align*}
Hence
\begin{align*}
\int_0^\infty g(\rho) \rho d\rho 
&= \int_0^\infty \left[
-\frac{2r}{\rho} \sin\te\cos\te 
+\frac 12 (1-\cos\te)\cos\te + \frac 12 \sin^2\te \right] \xi^2 \rho d\rho \\
&= 
\int_0^\infty
\frac{4r^2}{\rho^2+4r^2}\xi^2
\rho d\rho .
\end{align*}
Similarly, we have
\[
\int_0^\infty 
\left[
\frac 12 (\xi')^2 + \frac{1}{2\rho^2} (1-2\sin^2 \te) \xi^2 
\right]
\rho d\rho = \frac 12 
\nor{F\xi}{L^2}^2.
\]
Therefore, we obtain
$$
\int_0^\infty
\boK_2 \rho d\rho = 
\frac 12 \nor{F \xi}{L^2}^2 + 
\int_0^\infty 
\frac{4r^2}{\rho^2+4r^2}\xi^2
\rho d\rho,
\qquad 
F= \rd_\rho +\frac{\cos \te}\rho.
$$
For the higher order terms, we write as
$$
\boK_3 = 
\frac r2 \cos \te
\left(
\xi' -\frac 1\rho \cos \te \xi
\right) \xi^2 
+\frac 12 \sin \te \cos\te \xi^3 + 
\frac 1{\rho^2} R_5(\rho)\xi^3,
$$
$$
\boK_4 = 
\frac 18 \cos^2 \te \xi^4 + \frac 1{\rho^2} R_6(\rho) \xi^4 + 
\frac 1\rho
R_7(\rho) \xi'\xi^3,
$$
where 
$R_5,R_6,R_7$ are bounded functions independent of $\xi$. 
Then by noting $0\le |\xi|\le \te\le \frac{C}{\rho}$, we have
\begin{equation}\label{Ec2.1}
\begin{aligned}
\int_0^\infty
\left|\frac1{\rho^2}
R_1\xi^5 \right| + 
\left|\frac 1\rho R_2 \xi^5 \right|+ 
\left| 
R_3 \xi'\xi^4
\right|
+
\left| 
R_4\xi^5
\right|
+
\left|
\frac1{\rho^2} R_5\xi^3
\right|
+
\left|
\frac1{\rho^2} R_6\xi^4
\right|
+
\left|
\frac1{\rho} R_7\xi'\xi^3
\right| \rho d\rho
\\
\le 
C\int_0^\infty 
\frac1{\rho^2} |\xi|^3 + 
\frac1\rho |\xi'|\xi^2 \rho d\rho \le 
C \nor{\xi}{X}^3 .
\end{aligned}
\end{equation}
Also, by integration by part,
\begin{equation}\label{Ec2.2}
\begin{aligned}
\int_0^\infty \cos\te \xi' \xi^2 \rho d\rho
&=
\left[ \cos\te \frac 13 \xi^3 \rho \right]_{0}^{\infty} 
-\frac 13\int_0^\infty  \frac 1\rho (\cos \te \rho)' \xi^3  \rho d\rho \\
&=
\int_0^\infty 
\left(
-\frac 13 \frac 1\rho \sin^2 \te  \xi^3 -\frac 13 \frac 1\rho \cos \te \xi^3 \right) \rho d\rho,
\end{aligned}
\end{equation}
where we used $|\xi|\le \te$ to show that the boundary term vanishes. 
In the right hand side in \eqref{Ec2.2}, the first term is also bounded by $C\nor{\xi}{X}^3$ in the same way as \eqref{Ec2.1}. 
Therefore, we obtained
\begin{equation*}
\begin{aligned}
& E_0[f] -E_0[\te]\\
&\ge 
\frac 12 \nor{F\xi}{L^2}^2 
+ \int_0^\infty 
\frac{4r^2}{\rho^2+4r^2} \xi^2 - \frac r{6\rho} \cos \te 
\xi^3
- \frac r{2\rho} \cos^2 \te
 \xi^3 
+ \frac 12 \sin\te\cos\te \xi^3 
+ \frac 18 \cos^2\te \xi^4 \rho d\rho \\
&\qquad - C\nor{\xi}{X}^3.
\end{aligned}
\end{equation*}
Focusing on the leading order as $\rho\to \infty$, we split the integral terms as
\begin{align*}
&\int_0^\infty 
-\frac r6 \frac 1\rho \cos\te \xi^3
-\frac r2 \frac 1\rho \cos^2\te \xi^3 
+\frac 12 \sin\te\cos\te \xi^3 + \frac 18 \cos^2\te \xi^4 \rho d\rho\\
&= \int_{8r}^\infty 
\left[
-\frac r{6\rho} \xi^3
-\frac r{2\rho} \xi^3 
+\frac{2r}{\rho} \xi^3 + \frac 18 \xi^4 
\right] \rho d\rho \\
&\quad 
+ \int_0^{8r} \left[-\frac r6 \frac 1\rho \cos\te \xi^3
-\frac r2 \frac 1\rho \cos^2\te \xi^3 
+\frac 12 \sin\te\cos\te \xi^3 + \frac 18 \cos^2\te \xi^4 \right] \rho d\rho \\
&\quad
+
\int_{8r}^\infty \left[-\frac r6 \frac 1\rho (\cos\te -1) \xi^3
-\frac r2 \frac 1\rho (\cos^2\te -1) \xi^3 \right.\\
&\left. \hspace{60pt}
+\frac 12 (\sin\te\cos\te - \frac{4r}{\rho}) \xi^3 
+ \frac 18 (\cos^2\te -1) \xi^4 \right] \rho d\rho.
\end{align*}
The second integral in the right hand side is bounded by
$$
C \int_0^{8r} \frac 1\rho \xi^3 \rho d\rho 
\le C \int_0^{8r} \frac 1{\rho^2} \xi^3 \rho d\rho \le C \nor{\xi}{X}^3.
$$
The third integral in the right hand side is controlled by
$$
C \int_{8r}^\infty \frac 1{\rho^2} \xi^3 \rho d\rho \le C \nor{\xi}{X}^3.
$$
The integrand of the first term is rewritten as
$$
\frac{4r}{3\rho} \xi^3 + \frac 18\xi^4 
= \frac 18 \xi^2 \left( \xi^2 + \frac{32r}{3\rho} \xi \right)
= \frac 18 \xi^2 \left( \xi + \frac{16r}{3\rho} \right)^2 
- \frac{32 r^2}{9} \frac{\xi^2}{\rho^2}.
$$
Therefore, 
\begin{equation}\label{Ec2.3}
\begin{aligned}
&E_0[f] -E_0[\te] \\
&\ge 
\frac 12 \nor{F \xi}{L^2}^2 
+ \int_0^{8r} \frac{4r^2}{\rho^2 + 4r^2} \xi^2 \rho d\rho  
 + \int_{8r}^\infty \left( \frac{4r^2}{\rho^2+4r^2} - \frac{32}9 \frac {r^2}{\rho^2} \right) \xi^2 \rho d\rho  -C\nor{\xi}{X}^3
 \\
 &\ge 
 \frac 12 \nor{F\xi}{L^2}^2 + 
 \int_0^\infty \frac{1}{18} \frac{4r^2}{\rho^2 + 4r^2} \xi^2 \rho d\rho 
 -C\nor{\xi}{X}^3,
\end{aligned}
\end{equation}
where we used
$$
\frac{32}{9}\frac {r^2}{\rho^2} 
= \frac{17}{18} \frac{4r^2}{\rho^2 + \frac 1{16}\rho^2} 
\le \frac{17}{18} \frac{4r^2}{\rho^2 + 4r^2}\qquad 
\text{on}\quad [8r,\infty).
$$
\end{proof}
\begin{proof}[Proof of Lemma \ref{L2.3.1}]
We shall prove
$$
\frac 12 \nor{F f}{L^2}^2 
+ \int_0^\infty \frac{1}{18} \frac{4r^2}{\rho^2 + 4r^2} f^2 \rho d\rho 
\ge C \nor{f}{X}^2\qquad (f\in X),
$$
which, together with Lemma \ref{L2.3.3}, implies Lemma \ref{L2.3.1} by choosing $\eps_1$ sufficiently small. 
Suppose the contrary. 
We write the left hand side by $G(f)$.  Then there exists a sequence $\{f_n\}_{n=1}^\infty$ such that
$$
\nor{f_n}{X} =1,\qquad G(f_n) \xrightarrow{n\to\infty} 0.
$$
Since $f_n$ is bounded in $X$, there exists $f_*$ such that
$
f_n \rightharpoonup f_*$
in $X$ by taking subsequence if necessary. 
Now we claim that 
\begin{equation}\label{Ec2.4}
G(f_*) \le \liminf_{n\to\infty} G(f_n) =0,
\end{equation}
which implies $f_*=0$. Indeed, noting that
\begin{equation}\label{Ec2.6}
\begin{aligned}
G(f_n) 
&= \left\langle \frac 12 
\left(-\rd_\rho^2 -\frac 1\rho\rd_\rho +\frac 1{\rho^2} - \frac{32r^2}{(\rho^2 +4r^2)^2} 
+ \frac{1}{9} \frac{4r^2}{\rho^2 +4r^2} \right)
f_n,f_n\right\rangle_{L^2} \\
&= 
\frac 12 \nor{f_n}{X}^2 
-\int_0^\infty \frac{16r^2}{(\rho^2 +4r^2)^2} f_n^2 \rho d\rho
+ \int_0^\infty \frac{1}{18} \frac{4r^2}{\rho^2 +4r^2} f_n^2 \rho d\rho.
\end{aligned}
\end{equation}
The lower semi-continuity of norms yields
\begin{align*}
\frac 12 \nor{f_*}{X}^2 + \int_0^\infty \frac{1}{18} \frac{4r^2}{\rho^2 +4r^2} f_*^2 \rho d\rho
\le 
\liminf_{n\to\infty} 
\left(\frac 12 \nor{f_n}{X}^2 + \int_0^\infty \frac{1}{18} \frac{4r^2}{\rho^2 +4r^2} f_n^2 \rho d\rho
\right).
\end{align*}
On the other hand, compactness yields
\begin{equation}\label{Ec2.5}
\lim_{n\to\infty} \int_0^\infty \frac{16r^2}{(\rho^2 +4r^2)^2} f_n^2 \rho d\rho
=
\lim_{n\to\infty} \int_0^\infty \frac{16r^2}{(\rho^2 +4r^2)^2} f_*^2 \rho d\rho.
\end{equation}
Hence \eqref{Ec2.4} follows. 
Using \eqref{Ec2.5} again, we have
$
\lim_{n\to\infty} \int_0^\infty \frac{16r^2}{(\rho^2 +4r^2)^2} f_n^2 \rho d\rho =0$. 
Therefore, \eqref{Ec2.6} implies
\begin{align*}
G(f_n) 
&\ge
\frac 12 -\int_0^\infty \frac{16r^2}{(\rho^2 +4r^2)^2} f_n^2 \rho d\rho
\xrightarrow{n\to\infty} 
\frac 12,
\end{align*}
which contradicts \eqref{Ec2.4}. 
Consequently, the proof is complete.
\end{proof}

\subsubsection{Difference estimate}
\label{S2.3.3}

So far, we constructed a solution $\tilde{f}_\be$ to \eqref{E1.5} as a minimizer of $\boM_{\eps_0}$. 
In the last step, we upgrade the estimate of $\tilde{f}_\be$ as 
\begin{equation}\label{Ec2.9}
\nor{\tilde{f}_\be-\te}{X} \le C_\del \be^{1-\del}
\end{equation}
for $0<\del<1$. 
For this purpose, we make use of the Euler-Lagrange equation. 
Write $\xi_* = \tilde{f}_\be-\te$, and insert it into \eqref{E1.5}. 
Then the Taylor expansion yields
\begin{equation}\label{Ec2.7}
\left(
F^*F + \frac{8r^2}{\rho^2+4r^2}
+2\be^2\right) \xi_* + \be^2 \sin \te (1+\cos \te) 
+
\frac{4r}{(\rho^2+4r^2)^{1/2}} 
\xi_*^2 
+ \frac 12 \xi_*^3 + 
R_1 +R_2 =0,
\end{equation}
where 
\begin{align*}
R_1
&= 
\frac{R_{11}(\rho)}{\rho^2} \xi_*^2
+ 
\frac{R_{12}(\rho,\xi_*)}\rho 
\xi_*^3
+\frac{R_{13}(\rho,\xi_*)}\rho \xi_*^4,
\\
R_2
&= 
\be^2 \frac{R_{21}(\rho)}{\rho^2}\xi_* 
+ \be^2 \frac{R_{22}(\rho)}\rho  \xi_*^2
+ \be^2 R_{23}(\rho) \xi_*^3 
+ \be^2 R_{24}(\rho,\xi_*) \xi_*^4.
\end{align*}
Here $R_{1j}$, $R_{2j}$, $j=1,...,4$ are 
bounded functions uniformly in $\xi_*$. 
We note that the potential term and the nonlinearity are factorized as
\begin{equation}\label{EX2.1}
\frac{8r^2}{\rho^2+4r^2} 
+ \frac{4r}{(\rho^2+4r^2)^{1/2}} \xi_* + \frac 12 \xi_*^2 = 
\frac 12 \left(
\xi_*+ \frac{4r}{(\rho^2+4r^2)^{1/2}} 
\right)^2.
\end{equation}
Thanks to the positivity of \eqref{EX2.1}, we can rewrite \eqref{Ec2.7} as
\begin{equation}\label{Ec2.8}
\xi_* = \left(
F^* F + \frac 12 \left( \xi_*+ \frac{4r}{(\rho^2+4r^2)^{1/2}} \right)^2 + \be^2
\right)^{-1} (-\be^2 \sin \te (1+\cos \te)+ R_1 +R_2).
\end{equation}
The key ingredient for \eqref{Ec2.7} is the resolvent estimate 
for the operator 
\begin{equation}\label{EZB2}
F^* F + \frac12 \ovl{\xi}_*^2
\quad 
\text{with}\quad 
\ovl{\xi}_* = \xi_* + 4r (\rho^2+4r^2)^{1/2}.
\end{equation}
We claim the following:
\begin{prop}[Uniform resolvent estimate]
\label{Pc2.1}
There exist $C>0$, $\del_0>0$, $\be_*>0$ such that the following holds: Given $\xi_*\in X$, let $\ovl{\xi_*}$ be as in \eqref{EZB2}. Then we have
\begin{gather*}
\nor{(F^*F + \frac 12 \ovl{\xi}_*^2 +\be^2)^{-1} f}{X} \le 
\frac{C}{\be^{1-s}} \nor{\rho^s f}{L^2},\qquad 0\le s\le 1,
\\
\nor{(F^*F + \frac 12 \ovl{\xi}_*^2 +\be^2)^{-1} f}{X} \le \frac{C}{\be^{1+s}} \left(\nor{\rho^{1-s} \rd_\rho f}{L^2} 
+ \nor{\rho^{-s}f}{L^2}
\right),\qquad 0<s\le 1
\end{gather*}
for $f:(0,\infty)\to\R$ and $0<\be<\be_*$. 
\end{prop}
We will give a proof in the Appendix. 
Assuming this, we shall prove \eqref{Ec2.9}. 
We may take $\eps_0$ sufficiently small such that $\nor{\xi_*}{X} = \nor{\tilde{f}_\be-\te}{X} \le \del_0$, and thus Proposition \ref{Pc2.1} is available. 
By replacing $\be_0$ by $\min\{\be_0,\be_*\}$ if necessary, the estimates in Proposition \ref{Pc2.1} hold for all $0<\be<\be_0$. 
Since $\rho^{1-\del} \rd_\rho (\sin \te (1+\cos \te))$, $\rho^{-\del} \sin\te (1+\cos \te) \in L^2$, we have
\begin{equation}\label{EX2.4}
\nor{
(F^*F + \frac 12 \ovl{\xi}_*^2 +\be^2)^{-1} (\sin \te (1+\cos \te))}{X}
\le C_\del \be^{-1-\del}.
\end{equation}
On the other hand, combining the pointwise estimate $|\xi_*|\le \te \le \frac{C}\rho$, 
we estimate the remainders as 
\begin{align*}
\nor{(F^*F + \frac 12 \ovl{\xi}_*^2 + \be^2)^{-1} R_1}{X}
&\le 
C \nor{\rho R_1}{L^2} \\
&\le 
C\left(
\nor{\frac{1}{\rho}\xi_*^2}{L^2} + 
\nor{\xi_*^3}{L^2} +
\nor{\xi_*^4}{L^2}
\right)\\
&\le 
C \nor{\frac{1}{\rho} \xi_*}{L^2} 
\left(
\nor{\xi_*}{L^\infty} + 
\nor{\xi_*}{L^\infty}^2 
\right) \\
&\le 
C 
\nor{\xi_*}{X} 
\left(
\nor{\xi_*}{X} + 
\nor{\xi_*}{X}^2 
\right)
\le 
C \eps_0 \nor{\xi_*}{X},
\end{align*}
\begin{align*}
\nor{(F^*F + \frac 12 \ovl{\xi}_*^2 + \be^2)^{-1} R_2}{X}
&\le 
C \nor{\rho R_2}{L^2} \\
&\le 
C\be^2 \left(
\nor{\frac{1}{\rho}\xi_*}{L^2} + 
\nor{\xi_*^2}{L^2} +
\nor{\rho \xi_*^3}{L^2}
+
\nor{\rho \xi_*^4}{L^2}
\right)\\
&\le 
C \be^2 \nor{\frac{1}{\rho} \xi_*}{L^2} 
\left( 1+ 
\nor{\xi_*}{L^\infty} 
\right) \\
&\le 
C \be^2 
\nor{\xi_*}{X} 
\left(
1+ \nor{\xi_*}{X}
\right)
\le 
C \be_0^2 \nor{\xi_*}{X},
\end{align*}
where 
the implicit constants $C$ are independent of $\be$. 
Therefore, if $\be_0$ is sufficiently small, \eqref{Ec2.8} yields 
$$
\nor{\xi_*}{X} 
\le C \be^{1-\del} 
+ C \eps_0 \nor{\xi_*}{X}.
$$
Choosing $\eps_0$ sufficiently small, we obtain \eqref{Ec2.9}.


\subsection{Strict pointwise bound}\label{S2.4}
%
%
%
\begin{prop}\label{P2.4}
Let $f$ be a solution to \eqref{E1.5} with $0\le f \le \te$ on $(0,\infty)$ and $f(0)=\pi$. Then 
$0<f <\te$ on $(0,\infty)$. 
In particular, $f_\be$, $\tilde{f}_\be$ in Theorems \ref{P1}, \ref{P2} and \ref{P3} satisfy $0<f_\be, \tilde{f}_\be <\te$ on $(0,\infty)$. 
\end{prop}
%
%
%
\begin{proof}
If there exists $\rho=\rho_0\in (0,\infty)$ such that $f(\rho_0)=0$, then 
the condition $0\le f \le \te$ yields $f'(\rho_0)=0$. 
Then the uniqueness of ODE implies $f\equiv 0$, which is a contradiction.\par
Next suppose there exists $\rho=\rho_1\in (0,\infty)$ such that $f(\rho_1)=\te(\rho_1)$. 
Then the condition $f \le \te$ yields
$$
f'(\rho_1)=\te'(\rho_1)\qquad\text{and}\qquad f''(\rho_1)\le  \te''(\rho_1).
$$
Then 
\eqref{E1.5} yields
\begin{align*}
&f''(\rho_1) + \frac{1}{\rho_1} f'(\rho_1) \le \te''(\rho_1) + \frac{1}{\rho_1} \te'(\rho_1) \\
&\Longleftrightarrow \quad
\sin f(\rho_1) \left(\frac 1{\rho_1^2} \cos f(\rho_1) + 1 - \cos f(\rho_1) - \frac{2r}{\rho_1} \sin f(\rho_1) + \be^2 (1+\cos f(\rho_1))\right) \\
&\qquad\qquad \le 
\sin \te(\rho_1) \left(\frac 1{\rho_1^2} \cos \te(\rho_1) + 1 - \cos \te(\rho_1) - \frac{2r}{\rho_1} \sin \te(\rho_1)\right) \\
&\Longleftrightarrow 
\quad \be^2 \sin f(\rho_1) (1+\cos f(\rho_1)) \le 0.
\end{align*}
This contradicts $0<f \le \te <\pi$ at $\rho=\rho_1$. 
\end{proof}



\section{Shape analysis}
\label{SX3}


\subsection{Exponential decay}
\label{S3.1}

\begin{prop}\label{P3.1}
Let 
$\be>0$, $r>0$ and $f\in \boM$ be a solution to \eqref{E1.5} with \eqref{EA0.1}. 
Then, for any $0<\del<1$, 
there exist $C_\del=C_\del(\be,r)$ and $\rho_\del=\rho_\del(\be,r)$ such that
\begin{align}\label{E3.1}
C_\del e^{-(1+\del)\be\rho} \le 
f(\rho) \le C_\del e^{-(1-\del))\be \rho},
&&
-C_\del e^{-(1-\del)\be\rho} \le 
f'(\rho) \le -C_\del e^{-(1+\del))\be \rho}
\end{align}
for $\rho\ge  \rho_0$.
\end{prop}
\begin{proof}
The proof is based on \cite{LiMel18} and \cite{BerLio83}. 
We devide the proof into three steps.\par
%
%
\textit{Step 1.} 
Let $g(\rho):= \sqrt{\rho} f(\rho)$. Then computation yields
$$
g'' (\rho)
= -\frac 1{4\rho^2} g(\rho) + \frac{F(\rho)}{f(\rho)} g (\rho),
$$
where
$$
F(\rho) := \frac{1}{\rho^2} \sin f\cos f - 2r \frac{\sin^2 f}{\rho} + \sin f(1-\cos f) +\be^2 \sin f (1+\cos f).
$$
Here we have
$$
\frac{F(\rho)}{f(\rho)} = \frac{\sin f}{f} \left(
\frac{\cos f}{\rho^2} - \frac{2r \sin f}{\rho} + 1-\cos f + \be^2 (1+\cos f)
\right)
\xrightarrow{\rho\to\infty} 2\be^2 .
$$
Therefore, there exists $\rho_\del= \rho_\del(\be,r)>0$ such that 
$$
(1-\del)^2 \be^2 g \le g'' \le (1+\del)^2 \be^2 g  \qquad (\rho\ge \rho_0).
$$
%
%
\textit{Step 2.} 
Define 
$$
w(\rho) := e^{-(1-\del) \be \rho} (g' + (1-\del) \be g).
$$
Then
$$
w' = e^{-(1-\del) \be \rho} (g'' - (1-\del)^2 \be^2 g) \ge 0\qquad (\rho\ge \rho_\del).
$$
Here, we claim that
\begin{equation}\label{E3.2}
w(\rho) \le 0 \qquad (\rho\ge \rho_\del).
\end{equation}
Arguing by contradiction assuming that  there is $\rho_1\ge \rho_\del$ such that $w(\rho_1)>0$. Then 
$w(\rho) \ge w(\rho_1)>0$ if $\rho>\rho_1$. 
Then
\begin{align*}
 \left( e^{(1-\del) \be\rho} g \right)' \ge w(\rho_1) e^{2(1-\del)\be \rho}.
\end{align*}
Integrating over $[\rho_1,\rho]$, we have
$$
e^{(1-\del)\be\rho} g(\rho) - e^{(1-\del)\be\rho_1} g(\rho_1) \ge 
\frac{w(\rho_1)}{2(1-\del)\be} (e^{2(1-\del)\be\rho} - e^{2(1-\del)\be \rho_1}),
$$
and thus
\begin{align*}
g(\rho) 
&\ge \frac{w(\rho_1)}{2(1-\del)\be} e^{(1-\del)\be\rho} 
+ \left( e^{(1-\del) \be\rho_1} g(\rho_1) - \frac{w(\rho_1)}{2(1-\del)\be} e^{2(1-\del)\be\rho_1} \right) e^{-(1-\del)\be\rho}.  
\end{align*}
Recalling $g (\rho) = \sqrt\rho f(\rho)$, we get $\lim_{\rho\to\infty} f(\rho) =\infty$, 
which contradicts $f(\infty)=0$. 
From \eqref{E3.2}, we obtain
\begin{equation}\label{EFeb1}
g' + (1-\del) \be g \le 0
\quad \Longleftrightarrow 
\quad
\left( e^{(1-\del) \be \rho} g \right)' \le 0
\end{equation}
for $\rho\ge \rho_\del$. Integrating this over $[\rho_\del,\rho]$, we obtain
$$
g(\rho) \le e^{(1-\del) \be\rho_\del} g(\rho_\del) e^{-(1-\del) \be\rho} \qquad ( \rho \ge \rho_\del),
$$
yielding the upper bound of $f$ in \eqref{E3.1}. \par
%
%
\textit{Step 3.}
Define 
$$
z(\rho):= e^{-(1+\del) \be \rho} (g' + (1+\del) \be g).
$$
Then the same argument as in Step 2 gives $z(\rho)\ge 0$ for $\rho\ge\rho_\del$. 
Hence we have
\begin{equation}\label{EFeb2}
g'+(1+\del)\be g\ge 0\quad 
\Longleftrightarrow \quad 
\left( e^{(1+\del)\be \rho} g\right)' \ge 0
\end{equation}
if $\rho\ge \rho_\del$. 
Integrating this over $[\rho_\del,\rho]$, we obtain
$$
g(\rho) \ge e^{(1+\del)\be\rho_\del} g(\rho_\del) e^{-(1+\del) \be\rho} \qquad (\rho \ge \rho_\del),
$$
giving the lower bound of $f$ in \eqref{E3.1} by incorporating $\frac{1}{\sqrt{\rho}}$ into the exponential part. Finally, \eqref{EFeb1} and \eqref{EFeb2} yield
\begin{gather*}
- (1+\del) \be g \le g' \le - (1-\del) \be g \\
\Longleftrightarrow\qquad 
- (1+\del) \be f - \frac{1}{2\rho} f \le f' \le -(1+\del) \be f- \frac{1}{2\rho}f, 
\end{gather*}
which concludes the bound for $f'$ in \eqref{E3.1}, completing the proof.
\end{proof}


\subsection{Monotonicity}\label{S3.2}

\begin{prop}\label{P3.2}
Let $f\in \boM$ be a solution to \eqref{E1.5} with \eqref{EA0.1}. 
Suppose $r\le 1$ and $\be\le 1$. 
Then,
$$
f'+\frac{\sin f}\rho <0 \qquad \text{on } (0,\infty).
$$
In particular, we have $f'<0$ on $(0,\infty)$. 
\end{prop}

For the proof, we introduce several notations and giving a few preparatory lemmas. 
First, let
\begin{align*}
    Q(\rho) := f' + \frac 1\rho \sin f,&&
    \ovl{Q} (\rho) := f' - \frac 1\rho \sin f.
\end{align*}
Then define
$$
F (\rho) := 
\rho^2 Q(\rho)\ovl{Q} (\rho) =
\rho^2 (f')^2 - \sin^2 f.
$$
Note that it suffices to show 
$F>0$ on $\rho\in(0,\infty)$. 
The reason for focusing on $F$ rather than $Q$ is that it satisfies a structured ODE as
\begin{equation}\label{E3.3}
F'(\rho) = 
2\rho^2 f' \sin f 
P(\rho)
\end{equation}
where
\begin{equation*}
P(\rho)=
1- \cos f - \frac{2r}\rho \sin f
+\be^2 (1+\cos f).
\end{equation*}


%
%
%

\subsubsection{Reduction}

We first show that the monotonicity of $f$ implies the conclusion.
\begin{lem}\label{L3.1}
If $f'<0$ on $(0,\infty)$, then $F > 0 $. 
In particular, $Q<0$ holds.
\end{lem}
\begin{proof}
Suppose the contrary. 
Since $F(0)=F(\infty) =0$, $F$ has a non-positive minimum at some $\rho_1\in (0,\infty)$. 
Then by \eqref{E3.3}, $\rho_1$ satisfies
\begin{align}\label{E3.4}
F(\rho_1) \le 0,&& P(\rho_1) =0,&& 
\left[\rho P\right]' (\rho_1) \le 0.
\end{align}
The third condition is explicitly written as 
\begin{equation}\label{E3.5}
\left(-2r \cos f(\rho_1) + (1-\be^2) \rho_1 \sin f(\rho_1)
\right)
 f'(\rho_1) + 1 - \cos f(\rho_1) +\be^2 (1+\cos f(\rho_1))\le 0.
\end{equation}
Since $f'<0$, this particularly implies
\begin{equation}\label{E3.6}
(1-\be^2)\rho_1 \sin f (\rho_1) -2r \cos f(\rho_1)  \ge 0.
\end{equation}
Next, the first and second conditions yield
\begin{align}\label{E3.7}
Q(\rho_1)+ \frac1{2r} P(\rho_1)\ge 0 
&&
\Longleftrightarrow
&&
f'(\rho_1) \ge
- \frac{1}{2r} (1 - \cos f(\rho_1) +\be^2 [1+\cos f(\rho_1)]).
\end{align}
Therefore, from \eqref{E3.5}, \eqref{E3.6}, \eqref{E3.7}, we have
\begin{align*}
0 
&\ge 1 - \cos f(\rho_1) +\be^2 (1+\cos f(\rho_1))  +  f'(\rho_1) (-2r \cos f(\rho_1) + 
(1-\be^2) \rho_1 \sin f (\rho_1)) \\
&\ge 1 - \cos f(\rho_1) +\be^2 (1+\cos f(\rho_1)) \\
&\qquad - \frac{1}{2r} (
1 - \cos f(\rho_1) +\be^2 (1+\cos f(\rho_1))
 ) (-2r \cos f(\rho_1) + (1-\be^2)\rho_1 \sin f (\rho_1)) \\
&= (1 - \cos f(\rho_1) +\be^2 (1+\cos f(\rho_1)) ) \left(
1 + \cos f(\rho_1) - (1-\be^2) \frac {\rho_1}{2r} \sin f (\rho_1)
\right)  ,
\end{align*}
which implies
\begin{equation}\label{E3.8}
1 + \cos f(\rho_1) - (1-\be^2) \frac {\rho_1}{2r} \sin f (\rho_1)\le 0.
\end{equation}
Now we use the identity
$$
1+\cos f - \frac{\rho}{2r} \sin f
=
\frac{2 \cos \frac f2 \sin \frac{\te-f}2 }{\sin \frac \te 2} .
$$
Then the condition $0<f<\te$ yields
\begin{equation}\label{EX3.1}
1 + \cos f(\rho_1) - (1-\be^2) \frac {\rho_1}{2r} \sin f (\rho_1)
= 
\frac{2 \cos \frac f2 \sin \frac{\te-f}2 }{\sin \frac \te 2}
+ \frac{\be^2}{2r} \rho_1 \sin f(\rho_1)\ge 0,
\end{equation}
which must be $0$ by \eqref{E3.8}. 
This happens only if $f(\rho_1)=0$, 
which contradicts $f>0$ by Proposition \ref{P2.4}. Hence the conclusion follows.
\end{proof}
\subsubsection{Sign property}

Next, we observe the sign of $\ovl{Q}$.
\begin{lem}\label{L3.2}
If $r\le 1$, then $\ovl{Q}<0$ on $(0,\infty)$.
\end{lem}
\begin{proof}
We shall prove stronger bound
$$
N(\rho) := \ovl{Q} + r(1-\cos f) <0,\quad (0<\rho<\infty).
$$
We first compute
$$
N'= \left( -\frac{1+\cos f}{\rho} + r\sin f \right) N
+ \sin f \left[ (1-r^2) (1-\cos f) + \be^2 (1+\cos f) \right].
$$
Let
$$
A(\rho) := \frac{1-\cos f}{\rho} + r\sin f.
$$
Since $f$ decays exponentially, we can define 
$
\int_\rho^\infty A(s) ds. 
$
Then
\begin{equation}\label{E3.9}
\left( \rho^2 e^{ \int_\rho^\infty A(s) ds} N \right)'
=
\rho^2 e^{\int_\rho^\infty A(s) ds} \sin f \left[ (1-r^2) (1-\cos f) + \be^2 (1+\cos f) \right] >0.
\end{equation}
Hence, integrating \eqref{E3.9} over $(\rho,\infty)$ with $\rho>0$ implies the conclusion, 
by noting that $N$ decays exponentially.
\end{proof}

\subsubsection{Auxiliary lemma}

We will prove $f'<0$ via contradiction argument. 
For this, we prepare auxiliary lemma for the conclusive step.
\begin{lem}\label{L3.3}
Suppose 
$f'(\rho)\ge 0$ at some $\rho\in (0,\infty)$. Then, there exists $0<\rho_*<\rho_{**}<\infty$ such that 
\begin{equation}\label{EFeb3}
\begin{gathered}
f'(\rho_*) = f'(\rho_{**})=0,\qquad 
f'(\rho) >0 \quad \text{on } (\rho_*,\rho_{**}),\\
P(\rho_*)  >0,\qquad
P(\rho_{**}) \le 0.
\end{gathered}
\end{equation}
\end{lem}
\begin{proof}
We are going to prove this by contradiction argument. 
Suppose that there is no such $\rho_*,\rho_{**}$. We first claim that the assumption would imply the following. 
\begin{claim*}
There is a sequence $\{\rho_n\}_{n=1}^\infty$ with $\rho_1<\rho_2<\cdots$ such that .
\begin{enumerate}[(i)]
\item $f'(\rho_n)=0$;
\item $P(\rho_n)>0$;
\item There exists $\eps_n>0$ such that 
$f'(\rho)>0$ on $(\rho_n,\rho_n+\eps_n)$;
\item 
$\{\rho_n\}_{n=1}^\infty \subset \{\rho \ |\ f(\rho) \ge \frac\pi 2 \} $. 
In particular, $\{\rho_n\}_{n=1}^\infty$ is bounded.
\end{enumerate}
\end{claim*}
We construct such a sequence by induction. \par
\textit{Step 1}. 
Define
$$
\rho_1 := \inf \{ \rho\in (0,\infty) \ |\ f'(\rho) =0\}.
$$
By assumption, we have $\rho_1<\infty$. 
Moreover, Proposition \ref{PA3.5} shown later and the condition $0\le f\le \te$ implies $f'(\rho)<0$ for $0<\rho\ll1$, and hence $\rho_1>0$. 
We begin by showing that $\rho_1$ satisfies 
(i)--(iii). 
By definition, we have $f'(\rho_1)=0$ and $f'(\rho)<0$ on $(0,\rho_1)$. 
If $P(\rho_1) <0$, then $F'(\rho)>0$ 
on $(\rho_1-\eps , \rho_1)$ for sufficiently small $\eps$. 
Since $F(0)=0$ and $F(\rho_1) = -\sin^2 f(\rho_1)<0$, it follows that $F$ have a local minimum $\tilde{\rho}_1 \in (0,\rho_1)$. 
Since $f'(\tilde{\rho}_1)<0$, the same argument as Lemma \ref{L3.1} can be applied, leading to a contradiction. 
If $P(\rho_1) =0$, then 
\begin{align*}
(\rho P)' &= (-2r \cos f + (1-\be^2)\rho \sin f) f' + 1 -\cos f +\be^2(1+\cos f)\\
&= 1 -\cos f + \be^2 (1+\cos f) >0 \quad \text{at } \rho=\rho_1.
\end{align*}
Hence, for sufficiently small $\eps>0$, we have
\begin{center}
$\rho_1-\eps <\rho < \rho_1$\quad $\Longrightarrow$\quad 
$F'(\rho) = 2\rho f' \sin f \cdot \rho P >0$,
\end{center}
which leads to the same contradiction as above. Hence $P(\rho_1)>0$ follows. \par
Next we prove that there is $\eps_1>0$ such that
\begin{equation}\label{E3.10}
\rho\in (\rho_1,\rho_1+\eps_1)\quad \Longrightarrow\quad f'(\rho)>0.
\end{equation}
We first note that since $f'(\rho)<0$ for $\rho<\rho_1$, 
we have $f''(\rho_1)\ge 0$.  
We claim $f''(\rho_1)>0$ in fact. 
If $f''(\rho_1)=0$, then the Euler-Lagrange ODE \eqref{E1.5} gives
$$
\sin f \left( \frac 1{\rho^2} \cos f + P\right) =0\quad \text{at } \rho=\rho_1.
$$
Since $P(\rho_1)>0$, we necessarily have 
$\cos f(\rho_1)<0$. 
On the other hand, differentiating \eqref{E1.5} yields
\begin{equation}\label{E3.11}
\begin{aligned}
&-f''' - \frac 1\rho f'' + \frac 1{\rho^2} f' \\
&+ \cos f f' \left( \frac 1{\rho^2} \cos f - \frac {2r}\rho \sin f + 1-\cos f +\be^2 (1+\cos f) \right) \\
&+\sin f \left(
-\frac 1{\rho^3} \cos f- \frac 1{\rho^2} \sin f f' + \frac{2r}{\rho^2} \sin f 
-\frac{2r}\rho \cos f f' + (1-\be^2)\sin f f'
\right) =0.
\end{aligned} 
\end{equation}
At $\rho=\rho_1$, we have
$$
-f'''(\rho_1) + \sin f(\rho_1)\left(
-\frac 1{\rho_1^3} \cos f(\rho_1) + \frac{2r}{\rho_1^2} \sin f(\rho_1)
\right) =0
$$
using $f'(\rho_1)=f''(\rho_1)=0$. In particular, we have $f'''(\rho_1)>0$. 
This implies $f'(\rho) >0$ for $\rho\in (\rho_1-\eps,\rho_1+\eps)$ with sufficiently small $\eps$, contradicting $f'(\rho)<0$ on $\rho<\rho_1$. 
Hence $f''(\rho_1)>0$, and thus \eqref{E3.10} follows. 
Therefore, $\rho_1$ satisfies the desired conditions.\par
\textit{Step 2}. 
Suppose we have constructed $\rho_1<\cdots <\rho_n$ as in the claim. We shall construct $\rho_{n+1}$ satisfying (i)--(iii). 
Since $f'(\rho)>0$ on $(\rho_n,\rho_n+\eps_n)$ by assumption, and since $f'(\rho)<0$ for sufficiently large $\rho$ by \eqref{E3.1}, 
there exists $\rho>\rho_n$ with $f'(\rho)=0$. 
Then define 
$$
\tilde{\rho}_{n+1} := \inf\{\rho>\rho_n \ |\ f'(\rho)=0 \}.
$$
By the property of $\rho_n$, we have $f'(\rho)>0$ on $(\rho_n,\tilde{\rho}_{n+1})$, and
$P(\rho_n)>0$. In particular, we obtain 
$f''(\tilde{\rho}_{n+1}) \le 0$. 
Moreover, 
by the nonexistence of $\rho_*$, $\rho_{**}$ as in \eqref{EFeb3}, we have 
$P(\tilde{\rho}_{n+1}) > 0$. 
If $f''(\tilde{\rho}_{n+1})=0$, then the same argument as in Step 1 gives 
$f'''(\tilde{\rho}_{n+1})>0$. Hence $f'(\rho)>0$ if $0< \rho-\tilde{\rho}_{n+1} \ll 1$. 
In this case, $\rho_{n+1} := \tilde{\rho}_{n+1}$ satisfies all the three conditions. 
Now suppose $f''(\tilde{\rho}_{n+1})<0$. Then 
$f'(\rho)<0$ for $0 <\rho -\tilde{\rho}_{n+1} \ll 1$, and thus 
$\rho= \tilde{\rho}_{n+1}$ is a strict local maximum of $F$. Since $F(\infty)=0$ and $F(\tilde{\rho}_{n+1}) = -\sin^2 f(\tilde{\rho}_{n+1})<0$, there exists $\rho>\tilde{\rho}_{n+1}$ with $F'(\rho)=0$. 
Then we define
$$
\rho_{n+1} := \inf\{\rho>\tilde{\rho}_{n+1}\ |\ F'(\rho)=0 \}.
$$
Let us show that $\rho_{n+1}$ satisfies the desired properties. 
Since $F'(\rho_{n+1})=0$, we have either
$$
f'(\rho_{n+1})=0\qquad \text{or}\qquad P(\rho_{n+1})=0.
$$
Suppose $f'(\rho_{n+1})\neq 0$. Then by the definition of $\rho_{n+1}$, 
we have $f'(\rho)<0$ on $(\tilde{\rho}_{n+1} ,\rho_{n+1}+\eps)$ for some $\eps>0$, and $P>0$ on $(\tilde{\rho}_{n+1} ,\rho_{n+1})$. 
In particular, $F$ is decreasing on $(\tilde{\rho}_{n+1} ,\rho_{n+1})$, and thus 
$F(\rho_{n+1})< F(\tilde{\rho}_{n+1}) <0$. Therefore, $\rho=\rho_{n+1}$ satisfies the same condition as \eqref{E3.4}, which  
leads to the contradiction in the same way as Lemma \ref{L3.1}. 
Thus $f'(\rho_{n+1})=0$ follows. \par
Next we show $P(\rho_{n+1})>0$. 
If $P(\rho_{n+1})<0$, then by $P(\tilde{\rho}_{n+1})>0$, 
there is $\rho\in (\tilde{\rho}_{n+1},\rho_{n+1})$ such that $P=0$, and thus $F'(\rho)=0$. This contradicts the definition of $\rho_{n+1}$. 
If $P(\rho_{n+1})=0$, then 
\begin{align*}
(\rho P)' 
&= 1 -\cos f + \be^2 (1+\cos f) >0 \quad\text{at } \rho=\rho_{n+1},
\end{align*}
and thus $P<0$ for $\rho \in (\rho_{n+1}-\eps,\rho_{n+1})$ with sufficiently small $\eps>0$. 
Hence we arrive at the same contradiction as above. 
Hence $P(\rho_{n+1})>0$ follows. 
The inequality 
$f'(\rho) >0$ for $0<\rho-\rho_{n+1}\ll 1$ follows from the same argument as in the proof of \eqref{E3.10}.\par
\textit{Step 3.} 
We finally show the boundedness of $\{\rho_n\}_{n=1}^\infty$. 
For $\rho_n$, 
let $\tilde{\rho}_{n+1}$ be as defined in the last step. 
Since 
$f'(\tilde{\rho}_{n+1})=0$ and
$f''(\tilde{\rho}_{n+1})\le 0$, 
the Euler-Langrange ODE \eqref{E1.5} gives 
$$
\sin f \left( \frac1{\rho^2} \cos f + P  \right) \le 0
\qquad \text{at } \rho=\tilde{\rho}_{n+1}.
$$
Since $\sin f>0$ and $P(\tilde{\rho}_{n+1})>0$, we have $\cos f(\tilde{\rho}_{n+1})<0$, and thus 
$$
\rho_{n}\le \tilde{\rho}_{n+1} \le  
\sup\{\rho \ |\ f(\rho) \ge \frac\pi 2 \} <\infty,
$$
by noting $\lim_{\rho\to\infty} f(\rho) =0$. 
Therefore, we conclude the claim.\par
Now we conclude Lemma \ref{L3.3}. 
Since $\{\rho_n\}$ is bounded and increasing, there exists $\rho_\infty := \lim_{n\to\infty} \rho_n <\infty$. 
Since $f'(\rho_n)=0$ for all $n$, we have $f'(\rho_\infty)=0$. 
Also,
$$
f''(\rho_\infty) 
=\lim_{n\to\infty} \frac{f'(\rho_\infty) - f'(\rho_n)}{\rho_\infty -\rho_n} =0.
$$
Therefore, by \eqref{E3.11} and $\cos f(\rho_n)< 0$, 
we have $f'''(\rho_\infty) >0$, which implies that 
$f'(\rho)>0$ on $(\rho_\infty -\eps ,\rho_\infty +\eps) \setminus \rho_\infty$ for sufficiently small $\eps>0$. This contradicts that $\rho_n\to \rho_\infty$ and $f'(\rho_n)=0$. 
\end{proof}

\subsubsection{Conclusive step}

Now we conclude Proposition \ref{P3.2}, by showing the following lemma.
\begin{lem}\label{L3.4}
Suppose $f'(\rho)\ge 0$ at some $\rho\in (0,\infty)$. 
If $\be\le 1$, then there exists $\rho_{\#} \in (0,\infty)$ such that $\ovl{Q}(\rho_{\#})>0$. In particular, it contradicts Lemma \ref{L3.2}, concluding Proposition \ref{P3.2}. 
\end{lem}
\begin{proof}
Let $\rho_*$, $\rho_{**}$ be as in the claim in Lemma \ref{L3.3}. 
Since $(\rho P)'(\rho_{**}) = 1 -\cos f(\rho_{**}) +\be^2 (1+\cos f(\rho_{**}))>0$, 
we have $P(\rho)<0$ on $\rho_{**}-\eps <\rho < \rho_{**}$ with sufficiently small $\eps$. 
Therefore, if we define
$$
\rho_{\#} := \sup \{\rho \in (\rho_*, \rho_{**}) \ |\  P\ge 0 \},
$$
we necessarily have $\rho_{\#} \in (\rho_*, \rho_{**})$. 
We finally claim that $\ovl{Q}(\rho_{\#})>0$. 
If we 
suppose contrary, we have
\begin{align}\label{E3.12}
    \ovl{Q}(\rho_\#)\le 0,&&
    P(\rho_{\#})=0,&& 
    \left[\rho P\right]'(\rho_{\#})\le 0,
\end{align}
where the last two conditions follows from the definition. 
The first two conditions in \eqref{E3.12} imply 
$$
f'(\rho_{\#}) \le \frac1{2r} (1 -\cos f(\rho_{\#}) + \be^2(1+\cos f(\rho_{\#}))).
$$
On the other hand, the third inequality in \eqref{E3.12} is explicitly rewritten as
$$
(-2r \cos f(\rho_{\#}) + (1-\be^2)\rho_{\#} \sin(\rho_{\#})) f'(\rho_{\#}) 
+1-\cos f(\rho_{\#}) +\be^2 (1+\cos f(\rho_{\#}))\le 0.
$$
Since $f'(\rho_{\#})>0$ by $\rho_{\#}\in (\rho_*,\rho_{**})$, we necessarily have 
$$
-2r \cos f(\rho_{\#}) + (1-\be^2)\rho_{\#} \sin f(\rho_{\#}) <0.
$$
Therefore, if $\be\le 1$, then
\begin{align*}
0
&\ge (-2r \cos f(\rho_{\#}) + (1-\be^2)\rho_{\#} \sin f(\rho_{\#})) f'(\rho_{\#}) 
+1-\cos f(\rho_{\#}) +\be^2 (1+\cos f(\rho_{\#}))  \\
&\ge (-2r \cos f(\rho_{\#}) + (1-\be^2)\rho_{\#} \sin f(\rho_{\#})) \frac 1{2r} (1 - \cos f(\rho_{\#}) + \be^2(1+\cos f(\rho_{\#})))\\
&\qquad +1-\cos f(\rho_{\#}) +\be^2 (1+\cos f(\rho_{\#}))\\
&= \left(1-\cos f(\rho_{\#}) + \be^2 (1+\cos f(\rho_{\#}))\right) \left(1-\cos f(\rho_{\#}) + (1-\be^2)\frac{\rho_{\#}}{2r} \sin f(\rho_{\#})\right) > 0,
\end{align*}
which is a contradiction.

\end{proof}

\subsection{Derivative at the origin}\label{S3.3}


\begin{prop}\label{PA3.5}
Let $f$ be a solution to \eqref{E1.5} with $0\le f\le \te$ on $(0,\infty)$ and $f(0)=\pi$. Then the right derivative at the origin $f'(0+)$ exists, and $\lim_{\rho\to+0}f'(\rho)=f'(0)$. 
\end{prop}
\begin{proof}
We claim that for any $0<\del<1$, there exists $\rho_\del$ such that
\begin{equation}\label{EE3.2}
\pi - f(\rho) \le C_\del \rho^{1-\del}\qquad (0<\rho<\rho_\del)
\end{equation}
with some constant $C_\del>0$. 
Define $\tilde{f}:\R\to\R$ by
$$
\tilde{f}(y) = \pi - f(e^{-y}).
$$
Then $0\le \tilde{f} \le \pi$, and \eqref{E1.5} implies that
$$
\rd_{yy} \tilde{f} = F(y),\qquad 
F(y) := \sin \tilde{f} 
\left(
\cos \tilde f - e^{-2y} (1+\cos \tilde f -2r e^{y} \sin \tilde f + \be^2 (1-\cos \tilde{f}))
\right).
$$
Since $\tilde f(y) \to 0$ as $y\to\infty$, there exists $y_\del$ such that 
$$
y> y_\del \quad \Longrightarrow\quad
\frac{F(y)}{\tilde{f}(y)} > 1-\del.
$$
Hence the same argument as in Section \ref{S3.1} leads to \eqref{EE3.2}. \par
%
Next, we rewrite \eqref{E1.5} as
\begin{equation}\label{EE3.3}
\rd_\rho^2 f + \frac 1\rho \rd_\rho f - \frac{1}{\rho^2} f = G(\rho),
\end{equation}
$$
G(\rho) = -\frac{1}{\rho^2} (f- \frac 12 \sin 2f) 
+\sin f(1-\cos f) - \frac{2r}\rho \sin^2 f + \be^2 \sin f (1+\cos f).
$$
Since $\rho$, $\frac{1}{\rho}$ are the fundamental solutions to the homogeneous ODE for \eqref{EE3.3}, we can write $f$ as 
\begin{equation}\label{EE3.5}
f(\rho) = \frac{1}{2\rho} \int_\rho^\infty s^2 G(s) ds - \frac{\rho}{2} \int_\rho^\infty G(s) ds.
\end{equation}
By \eqref{E3.1} and \eqref{EE3.2}, the asymptotics of $G(\rho)$ as $\rho\to 0+$ and $\rho\to\infty$ are
\begin{equation}\label{EE3.4}
G(\rho) = -\frac{\pi}{\rho^2} + O_{\rho\to 0+}(\rho^{1-\del}),\qquad\text{and}\qquad 
G(\rho) = O_{\rho\to\infty} 
(e^{-\frac 14 \be\rho}).
\end{equation}
\eqref{EE3.4} implies
\begin{align*}
\int_\rho^\infty G(s) ds 
= -\frac{\pi}{\rho} + \ovl{G}(\rho) ,\qquad \ovl{G}(\rho)=O_{\rho\to 0}(1).
\end{align*}
Combining this with \eqref{EE3.5}, we have
$$
\frac{1}{2\rho} \int_\rho^\infty s^2 G(s) ds
= f(\rho) - \frac{\pi}2 + \frac{\rho}2 \ovl{G}(\rho) = O_{\rho\to 0+} (1), 
$$
which implies
\begin{equation}\label{EE3.6}
\int_0^\infty s^2 G(s) ds =0.
\end{equation}
In particular, it follows from \eqref{EE3.4} that
\begin{align}
\int_\rho^\infty s^2 G(s) ds 
&= -\int_0^\rho s^2 G(s) ds = \pi \rho +O(\rho^{4-\del}),
\end{align}
which, combined with \eqref{EE3.5}, yields
$$
f(\rho) = 
\frac{1}{2\rho} \left(
\pi \rho + O_{\rho\to 0}(\rho^{4-\del}) 
\right)
-\frac \rho2 \left(
-\frac \pi\rho + \ovl{G}(\rho)
\right)
= \pi 
- \frac 12 \rho \ovl{G}(\rho)
+ O(\rho^{3-\del}).
$$
Therefore
$$
\lim_{\rho\to 0+}
\frac{f(\rho)-\pi}{\rho}
= -\frac 12 \ovl{G}(0) .
$$
Finally, a direct computation yields
\begin{align*}
f'(\rho) = -\frac{1}{2\rho^2} \int_\rho^\infty s^2G(s) ds
- \frac 12 \int_\rho^\infty G(s) ds
= -\frac 12 \ovl{G} (\rho) + O(\rho^{2-\del}) \xrightarrow{\rho \to +0} -\frac12 \ovl{G}(0),
\end{align*}
which completes the proof. 
\end{proof}

\section{Stability}
\label{S3}

\subsection{Linear stability}\label{S3.4}

This section is devoted to the stability analysis of $\bn_\be$. 
To state the claim, we first derive the Hessian of $E_\be$ at $\bn_\be$. 
Let $\{\bn_{\be,\eps}\}_{\eps\in\R}$ be any smooth deformation of $\bn_\be$ with constraint $|\bn_{\be,\eps}|=1$. 
If we set the Taylor expansion 
with respect to $\eps$ as 
$\bn_{\be,\eps} = \bn_\be + \eps \bphi + \frac{\eps^2}2 \bxi + O(\eps^3)$, the unit length constraint imposes the compatibility conditions
\begin{align*}
\bphi \cdot \bn_\be=0,&&
|\bphi|^2 + \bxi \cdot \bn_\be =0.
\end{align*}
Using this, the energy can be algebraically expanded as
$$
E_\be[\bn_{\be,\eps}] - E_\be[\bn_\be]
= \frac {\eps^2}2 \boH_\be[\bphi] + O(\eps^3),
$$
where
\begin{equation}\label{EF3.4}
\boH_\be [\bphi] \equiv 
\boH_{\be.\bn_\be} [\bphi] := 
\int_{\R^2} |\nab \bphi|^2 + 2r \bphi \cdot \nab \times \bphi
+ (1-\be^2) \bphi_3^2 
- \La_\be (\bn_\be) |\bphi|^2 dx
\end{equation}
with
$$
\La_{\be} (\bn_\be) := |\nab \bn_\be|^2 + 2r \bn_\be \cdot \nab \times \bn_\be 
-(\be^2+1) n_{3,\be} - (\be^2-1) n_{3,\be}^2.
$$
The main purpose of this section is to show the 
positivity of $\boH_\be$ for the tangent field:
\begin{thm}\label{T3.1}
Let $f_\be$ be a local minimizer of $\tilde{E}_\be [f]$ in $\boM_*$, where $\boM_*$ is the space defined in \eqref{EFeb4}. 
Then if $0<r\le \frac 12$, the corresponding equivariant map $\bn_\be$ is linearly stable in the sense that 
$$
\boH_\be [\bphi] \ge 0 \quad \text{for } 
\bphi\in H^1(\R^2:\R^3) \text{ with } 
\bphi \cdot \bn_\be =0.
$$
\end{thm}
First,
we reduce the proof to showing 
\begin{equation}\label{EF3.3}
\boH_\be[\bphi] \ge 0\quad \text{for } 
\bphi\in C_0^\infty(\R^2\setminus \{0\}:\R^3) \text{ with } 
\bphi \cdot \bn_\be =0.
\end{equation}
Indeed, we may suppose $\bphi\in H^\infty(\R^2:\R^3)$. 
One may also assume that $\bphi(0)=0$ by the following procedure: 
Suppose $\bphi(0)\neq 0$. 
Since $\boH_\be[\rd_j\bn_\be]=0$ for $j=1,2$, we have
\begin{equation}\label{EF3.2}
\boH_\be [\bphi ] = 
\boH_\be [\bphi + c_1 \rd_1 \bn_\be
+c_2 \rd_2 \bn_\be]
\end{equation}
for all $c_1,c_2\in \R$. 
Moreover, using 
Proposition \ref{PA3.5}, we observe that
\begin{align}\label{EE3.1}
\rd_1 \bn_\be (0) =
\begin{pmatrix}
0 \\ -f'(0+) \\ 0
\end{pmatrix}
,&&
\rd_2 \bn_\be (0) =
\begin{pmatrix}
f'(0+) \\ 0 \\ 0
\end{pmatrix}
.
\end{align}
Since $0<f < \te$, we have 
$f'(0+)\le \te'(0+) <0$. 
Since $\bphi (0) \cdot \bn_\be(0)= \bphi (0)\cdot \boe_3 =0$, \eqref{EE3.1} implies that there exists $(c_1,c_2)\in \R^2\setminus 0$ such that 
$\tilde{\bphi}= \bphi +c_1 \rd_1 \bn_\be + c_2 \bn_\be$ satisfies $\tilde\bphi(0)=0$. 
Therefore, it suffices to replace $\bphi$ by $\tilde{\bphi}$. 
Now, Lemma 2 in \cite{LiMel18} can be applied, implying that $\bphi$ can be approximated in $H^1$ by a sequence $\{\bphi_j\}_{j=1}^\infty\subset C_0^\infty (\R^2\setminus \{ 0\}:\R^3)$. Hence we can reduce Theorem \ref{T3.1} to \eqref{EF3.3}.\par 
%
%
%
The proof of \eqref{EF3.3} is divided into several steps. 
In the proof, we abbreviate $f_\be$ by $f$ for simplicity. \par
\textit{Step 1.} Moving frame reduction. 
For $\bphi\in T_{\bn}\mathbb{S}^2$, we write
$$
\bphi = u_1 \bJ_1 + u_2 \bJ_2
$$
with
\begin{equation}\label{E3.16}
\bJ_1 =
\begin{pmatrix}
\cos\psi \\
\sin\psi\\
0
\end{pmatrix}
,
\qquad 
\bJ_2 =
\begin{pmatrix}
-\sin\psi \cos f(\rho) \\
\cos\psi \cos f(\rho)\\
-\sin f(\rho)
\end{pmatrix}
.
\end{equation}
Then, computation yields
\begin{align*}
&|\nab \bphi|^2 = |\nab u|^2 + \frac{2\cos f}{\rho^2} u\times \rd_\psi u
+ \frac 1{\rho^2} u_1^2 + \left(f'^2 + \frac{\cos^2 f}{\rho^2} \right)u_2^2;
\\
&\bphi \cdot (\nab\times \bphi) = 
-\frac 1\rho \sin f (u\times \rd_\psi u) 
+\left(f'- \frac 1\rho \sin f\cos f\right) u_2^2;
\\
&\phi_3^2 = u_2^2 \sin^2 f;
\\
&- \La_\be (\bn)
=
-\left(f'^2 + \frac{\sin^2 f}{\rho^2}\right) 
-2r \left(f' + \frac 1\rho \sin f\cos f\right) \\
&\hspace{170pt}
+(\be^2+1)\cos f + 
(\be^2-1) \cos^2 f.
\end{align*}
Hence
$$
\boH_\be [\bphi] =
\int_{\R^2}
\left[|\nab u|^2 + \left(
\frac{2\cos f}{\rho^2} - \frac{2r \sin f}{\rho}
\right) (u\times \rd_\psi u)
+ G_{1,\be}(f) u_1^2 + G_{2,\be}(f) u_2^2 \right] dx,
$$
where
\begin{equation}\label{EZ11}
\begin{aligned}
G_{1,\be}(f) &= \frac{\cos^2 f}{\rho^2} - f'^2 -2r f' 
-\frac{2r}{\rho} \sin f \cos f \\
&\hspace{80pt} + \cos f (1-\cos f) + \be^2 \cos f (1+\cos f),
\end{aligned}
\end{equation}
\begin{equation}\label{EZ12}
\begin{aligned}
G_{2,\be}(f) &= 
\frac{\cos^2 f - \sin^2 f}{\rho^2} - \frac{4r}{\rho} \sin f \cos f 
+ (1-\be^2) \sin^2 f \\
&\hspace{80pt} + \cos f (1-\cos f) + \be^2 \cos f (1+\cos f)
.
\end{aligned}
\end{equation}
\textit{Step 2.} Fourier splitting. 
Apply the Fourier expansion to $u_j$, $j=1,2$ with respect to $\psi$:
\begin{equation}\label{EE3.9}
u_j(\rho,\psi) =
a_j^{(0)} (\rho) + \sum_{n=1}^\infty \left( a_j^{(n)} (\rho) \cos n\psi 
+ b_j^{(n)} (\rho) \sin n\psi \right).
\end{equation}
Then computation similar to \cite{LiMel18, IbrShi23} gives
\begin{equation}\label{EE3.11}
\boH_\be [\bphi] = 2\pi \boH^{(0)} [a_1^{(0)} ,a_2^{(0)}]
+
\pi \sum_{n=1}^\infty \left( \boH^{(n)} [a_1^{(n)}, b_2^{(n)}] +
\boH^{(n)} [b_1^{(n)}, -a_2^{(n)}]  \right),
\end{equation}
with
\begin{align*}
\boH_\be^{(n)} [a,b] =
\int_0^\infty &
\left[ (a')^2 + (b')^2 
+ \frac{n^2}{\rho^2} (a^2+b^2)\right. \\
&\quad 
\left.+4n \left( \frac{\cos f}{\rho^2} - \frac{r\sin f}{\rho} \right) ab 
+ G_{1,\be}(f) a^2 +G_{2,\be}(f) b^2 
\right]
\ \rho d\rho.
\end{align*}
\textit{Step 3.} Monotonicity in $n$. We claim the following:
\begin{lem}\label{L3.5}
If $r\le \frac 12$ and $k\ge 1$, then $\boH^{(k)}$ is non-decreasing with respect to $k$:
$$
\boH^{(1)} \le \boH^{(2)} \le \boH^{(3)} \le \boH^{(4)} \le \cdots.
$$
\end{lem}
Therefore, to show the positivity of $\boH_\be$, we only need to check the cases $k=0,1$. 
\begin{proof}
$$
\boH_\be^{(k+1)} [a,b] - \boH_\be^{(k)} [a,b]
=
\int_0^\infty \left[\frac{2k+1}{\rho^2} \left( a^2+b^2 \right) + 4
\left( \frac{\cos f}{\rho^2} - \frac{r\sin f}{\rho} \right) ab\right] \rho d\rho.
$$
We first claim that if $r \le \frac 12$, then
\begin{equation}\label{E3.13}
\left| r\rho \sin f(\rho) - \cos f(\rho) \right| \le \frac 32.
\end{equation}
Indeed, 
if $\rho\ge 2r = \te^{-1} (\frac \pi2)$, then 
$0\le f\le \te \le \frac \pi 2$, and thus
$$
r\rho \sin f - \cos f \le r\rho \sin \te - \cos \te 
= \frac{4r^2\rho^2}{\rho^2+4r^2} - \frac{\rho^2-4r^2}{\rho^2+4r^2} 
= 4r^2-1 + \frac{8r^2 (1-2r^2)}{\rho^2+4r^2}\le 1.
$$
If $\rho \le 2r$, then
$$
r\rho \sin f - \cos f \le 2r^2 \sin f- \cos f \le 2r^2 +1 \le \frac 32.
$$
On the other hand, since $\sin f\ge 0$, we have
$$
r\rho \sin f - \cos f \ge -\cos f \ge -1.
$$
Hence \eqref{E3.13} follows. Therefore
$$
\boH_\be^{(k+1)} [a,b] - \boH_\be^{(k)} [a,b] 
\ge 
\int_0^\infty \frac{1}{\rho^2} \left\{
(2k+1) (a^2+b^2) - 6 |ab|
\right\} \rho d\rho.
$$
If $k\ge 1$, then
\begin{align*}
\boH_\be^{(k+1)} [a,b] - \boH_\be^{(k)} [a,b]  &\ge 
\int_0^\infty \frac{1}{\rho^2} \left\{
3 (a^2+b^2) - 6 ab
\right\} \rho d\rho \\
&=
\int_0^\infty \frac{1}{\rho^2} (|a|-|b|)^2 \rho d\rho \ge 0.
\end{align*}
\end{proof}
\textit{Step 4.} Positivity of $\boH_\be^{(0)}$. 
We claim that if $r\le 1$, then
$$
\boH^{(0)}_\be [a,b] 
= 
\int_0^\infty \left[ (a')^2 + (b')^2 + G_{1,\be}(f) a^2 + G_{2,\be} (f) b^2\right] \rho d\rho \ge 0 
$$
for $(a,b) \in C_0^\infty (0,\infty)^2$. 
First of all, we separate it as
$$
\boH^{(0)}_\be [a,b] = A^{(0)} [a] + B^{(0)}[b],
$$
where
\begin{align}\label{EX3.2}
A^{(0)} [a] = \int_0^\infty \left[ (a')^2 + G_{1,\be} (f) a^2 \right]\rho d\rho,\qquad 
B^{(0)} [b] = \int_0^\infty \left[(b')^2 + G_{2,\be} (f) b^2 \right]\rho d\rho,
\end{align}
and $G_{j,\be}$ given by \eqref{EZ11}, \eqref{EZ12}. 
Here we observe that the positivity of $B^{(0)}$ follows from the local minimality of $f\in \boM_*$. 
Indeed, let $g\in C_0^\infty(0,\infty)$. 
Since Proposition \ref{P2.4} implies $0<f<\te$ on $(0,\infty)$, 
we have
$f+\eps g\in \boM_*$ for sufficiently small $\eps>0$. 
Expanding $\tilde{E}_\be$ around $f$, we obtain
\begin{align*}
0\le \tilde{E}_\be [f+\eps g] - \tilde{E}_\be [f] 
&= \eps^2 B^{(0)} [g] + O(\eps^3), 
\end{align*}
which implies that $B^{(0)}[g]$ is nonnegative. 
For $A^{(0)}$, we apply the transformation $a = (\sin f) \xi$. Then, integration by parts (see Lemma 4 in \cite{LiMel18}) yields
$$
A^{(0)}[a] = \int_0^\infty \sin^2f (\xi')^2 \rho d\rho 
+ \int_0^\infty \xi^2 L_0 \sin f \cdot \sin f d\rho,\qquad 
L_0= - \frac d{d\rho} \rho \frac d{d\rho} + G_{1,\be} (f)
\rho.
$$
Computation using \eqref{E1.5} yields $L_0 \sin f = -2r f' \sin f \rho$, which gives
\begin{equation}\label{EE3.12}
A_0[a] = \int_0^\infty \left[\sin^2 f (\xi')^2 - 2r \sin^2 f f' \xi^2 \right]\rho d\rho > 0
\end{equation}
for $a\in C_0^\infty (0,\infty) \setminus \{0\}$, 
where we used $f'<0$. \par
\textit{Step 5.} 
Recall 
\begin{align*}
\boH_\be^{(1)} [a,b] =
\int_0^\infty &\left[(a')^2 + (b')^2 
+ \frac{1}{\rho^2} (a^2+b^2)\right. \\
&\quad \left.+4 \left( \frac{\cos f}{\rho^2} - \frac{r\sin f}{\rho} \right) ab 
+ G_{1,\be}(f) a^2 +G_{2,\be}(f) b^2\right] \rho d\rho.
\end{align*}
We change variables as
$$
a = \frac{\sin f}{\rho} \xi,\qquad
b = - (f') \eta,
$$
and get
\begin{equation}\label{EE3.10}
\begin{aligned}
&\boH^{(1)}_\be \left[ \frac{\sin f}{\rho} \xi, -f' \eta \right] \\
&=
\int_0^\infty \left[\frac{\sin^2 f}{\rho} (\xi')^2 + \rho (f')^2 (\eta')^2 
+ f' \left( \frac{2}{\rho^2} \sin f \cos f  - \frac{2r \sin^2 f}{\rho}  \right)
(\xi - \eta)^2\right] d\rho \\
&= 
\int_0^\infty \left[\frac{\sin^2 f}{\rho} \left( \xi' - \frac{\xi-\eta}{\rho} \right)^2 
+\frac{\sin^2 f (1-2r\rho^2 f')}{\rho^3} 
\left\{\xi-\eta + \frac{\rho}{1-2r\rho^2 f'} \eta' \right\}^2\right.\\
&\hspace{40pt} 
\left.+ \frac{\rho^2 (f')^2 - \sin^2 f - 2r\rho^4 (f')^3 }{\rho (1-2r\rho^2 f')} (\eta')^2\right] d\rho. 
\end{aligned}
\end{equation}
By Proposition \ref{P3.2}, 
we have 
$f'<0$ and $\rho^2 (f')^2 - \sin^2 f>0$, implying that all the coefficients in \eqref{EE3.10} are positive. 
Hence, we have
$$
\boH^{(1)}_\be [a,b] =
\boH^{(1)}_\be \left[ \frac{\sin f}{\rho} \xi, -f' \eta \right] >0
$$
for $(a,b) \in C_0^\infty (0,\infty)^2 \setminus \{(0,0) \}$. 
Therefore, we conclude \eqref{EF3.3}, as desired.

\subsection{Instability}\label{S3.5}

\begin{thm}\label{T3.2}
Let $\tilde{f}_\be$ be the critical point as in Theorem \ref{P3}. 
%
%
Let $\tilde{\bn}_\be$ be the map in \eqref{E1.3} with $f=\tilde{f}_\be$. 
Then for $r>1$, 
there exists $\be_1=\be_1(r)>0$ 
such that if $0<\be< \be_1$, 
$\tilde{\bn}_\be$ is an unstable critical point of $E_\be$. 
More precisely, for any neighborhood $\boU$ of $\tilde{\bn}_\be$ in $\boM$, 
there exists $\bn\in \boU$ such that $E_\be[\bn] -E_\be[\tilde{\bn}_\be] <0$.
\end{thm}

For the proof, we again use the energy expansion
$$
E_\be [\tilde{\bn}_\be+\bphi] - E_\be [\tilde{\bn}_\be] =
\frac 12 \boH_{\be,\tilde{\bn}_\be}[\bphi ],
$$
where $\boH_{\be,\tilde{\bn}_\be}[\bphi ]$ is the Hessian as in \eqref{EF3.4} with $\bn_\be$ replaced by $\tilde{\bn}_\be$. 
The idea is that the Hessian of $E_\be$ at $\tilde{f}_\be$ is close to that of $E_0$ at 
$\te$, which is known to be unstable if $r>1$ as established in \cite{IbrShi23}. 
\begin{proof}
To prove the instability, it suffices to show that there exists $\bphi\in H^1(\R^2)$ with $\bphi\cdot \tilde{\bn}_\be =0$ 
such that $\boH_\be [\bphi] <0$. (See \cite{IbrShi23} for the detail.) 
We use the moving frame
$$
\bphi = u_1 \bJ_1 + u_2 \bJ_2
$$
with $\bJ_1,\bJ_2$ defined in \eqref{E3.16}. 
Then 
\begin{align*}
\tilde{\boH}_{\be,\tilde{f}_\be} [u]
&:= \boH_{\be,\tilde{\bn}_\be} [\bphi] \\
&=
\int_{\R^2} \left[|\nab u|^2 + \left(
\frac{2\cos \tilde{f}_\be}{\rho^2} - \frac{2r \sin \tilde{f}_\be}{\rho}
\right) (u\times \rd_\psi u)
+ G_{1,\be}(\tilde{f}_\be) u_1^2 + G_{2,\be}(\tilde{f}_\be) u_2^2 \right]dx.
\end{align*}
Here we recall the result in \cite{IbrShi23} that if $r>1$, then there exists $u_* = (u_{1,*}, u_{2,*})\in C_0^\infty(\R^2)^2$ such that
$$
\tilde{\boH}_{0,\te} [u_*] <0.
$$
Now we decompose
\begin{align*}
\tilde{\boH}_{\be,\tilde{f}_\be} [u_*]
=
\tilde{\boH}_{0,\te} [u_*] +& \int_{\R^2} \left[2 \left( \frac{\cos \tilde{f}_\be -\cos \te}{\rho} 
- \frac{r(\sin \tilde{f}_\be - \sin \te)}{\rho}  \right) (u_* \times \rd_\psi u_*)\right. \\
&\qquad \left.+ (G_{1,\be} (\tilde{f}_\be) - G_{1,0}(\te)) u_{1,*}^2 
+ (G_{2,\be} (\tilde{f}_\be) - G_{2,0}(\te)) u_{2,*}^2 \right] dx. 
\end{align*}
Now, we take $0<\del<1$. 
Then \eqref{Ec2.9} implies 
\[
|\cos \tilde{f}_\be - \cos \te | +|\sin \tilde{f}_\be - \sin \te | \le C |\tilde{f}_\be-\te| \le C_\del \be^{1-\del},
\]
and hence
\begin{align*}
&\left|\int_{\R^2} 
2 \left( \frac{\cos \tilde{f}_\be -\cos \te}{\rho} 
- \frac{r(\sin \tilde{f}_\be - \sin \te)}{\rho}  \right) (u_* \times \rd_\psi u_*) dx
\right|\\
&\le C (|\cos \tilde{f}_\be -\cos \te| +|\sin \tilde{f} -\sin \te|) 
\nor{u_*}{L^2} \nor{\nab u_*}{L^2}
\le
C_\del \be^{1-\del}.
\end{align*}
Similarly, we have
\begin{align*}
&|G_{1,\be} (\tilde{f}_\be) - G_{1,0}(\te)| + |G_{2,\be} (\tilde{f}_\be) - G_{2,0}(\te)|\\
&\le |G_{1,0} (\tilde{f}_\be) - G_{1,0} (\te)| 
+|G_{2,0} (\tilde{f}_\be) - G_{2,0} (\te)|
+ C\be^2 \cos \tilde{f}_\be (1+\cos \tilde{f}_\be)\\
&\le C|\tilde{f}_\be-\te| 
+ C|\tilde{f}_\be' -\te'|+ C \be^2 
\le C_\del \be^{1-\del} + 
C |\tilde{f}_\be' -\te'|.
\end{align*}
Therefore,
\begin{align*}
&\left|\int_{\R^2} \left[(G_{1,\be}(\tilde{f}_\be) - G_{1,0}(\te)) u_{1,*}^2 + (G_{2,0} (\tilde{f}_\be) - G_{2,0} (\te)) u_{2,*}^2 \right]dx \right|\\
&\le C_\del \be^{1-\del} \nor{u_*}{L^2}^2 + C \nor{\tilde{f}_\be'-\te'}{L^2} \nor{u_*}{L^4}^2 \le C_\del \be^{1-\del}.
\end{align*}
Consequently, we obtain
\begin{equation}\label{EF3.5}
\left| 
\tilde{\boH}_{\be,\tilde{f}_\be} [u_*] - \tilde{\boH}_{0,\te} [u_*]
\right| \le C_\del \be^{1-\del}.
\end{equation}
Now we take $\be_1=\be_1(r)$ sufficiently small such that $C_\del \be_1^{1-\del} \le -\frac 12 \tilde{\boH}_{0,\te}[u_*]$. Then for $0<\be<\be_1$, \eqref{EF3.5} implies $\tilde{\boH}_{\be,\tilde{f}_\be} [u_*] < \frac 12 \tilde{\boH}_{0,\te}[u_*] <0$. As a consequence, $\bphi_* = u_{1,*} \bJ_1 + u_{2,*} \bJ_2$ is the desired tangent field.
\end{proof}


\section{Appendix}
\label{S4}

\subsection{Regularity of $\bn$}\label{S4.1}

\begin{prop}
Let $\bn \in H^1_{\boe_3}(\R^2)$ be a solution to \eqref{EF1.5} 
with $|\bn|=1$. 
Then $\bn\in C^\infty$. 
\end{prop}

The proof is essentially due to \cite{Helein}; see also \cite{LiMel18}. For reader's convenience, we give a short outline of the proof.

\begin{proof}
We may restrict the domain to an arbitrary open ball $B\subset \R^2$. 
Let $\bom := \bn- \boe_3$. Then $\bom$ solves
\begin{equation}\label{EF3.6}
\begin{aligned}
-\Del \bom 
&=-2r\nab\times \bom + k m_3\boe_3 
- (h+k) \bom
\\
&\qquad+
\left( 
|\nab \bom|^2 + 2r (\bom+\boe_3) \cdot \nab \times \bom - h m_3 
- k m_3 (m_3+2)
\right) (\bom + \boe_3).
\end{aligned}
\end{equation}
Observe that the RHS of \eqref{EF3.6} consists of $L^2$-terms except $|\nab \bom|^2 (\bom+\boe_3)
$. Thus we rewrite it in the form 
$$
-\Del\bom = 
|\nab \bom|^2 (\bom+\boe_3) 
+ \ovl{\bom},\qquad 
\ovl{\bom} \in L^2(B:\R^3).
$$
Letting $\bm{w}_0$ be the solution to the Dirichlet problem
\begin{align*}
-\Del \bm{w}_0 = \ovl{\bom}\quad \text{in } B,&& 
\bm{w}_0 = \bom 
\quad \text{on } \rd B,
\end{align*}
we have $\bm{w}_0 \in H^1(B)\cap C (B)$ (see Sections 8 in \cite{GilTru}), and 
%
\begin{align}\label{EF3.7}
-\Del (\bom -\bm{w}_0) 
= |\nab \bom|^2 (\bom +\boe_3)
\quad \text{in } B,&& 
\bm{\bom - \bm{w}_0} = 0 
\quad \text{on } \rd B.
\end{align}
Here, we observe that $|\bn|=1$ implies
$$
|\nab \bom|^2 (\bom+\boe_3)
=|\nab \bn|^2 \bn = 
\rd_l \bn \times (\bn \times \rd_l \bn),
$$
where repeated indices indicates the summation in these indices. Thus for $i=1,2,3$, we can write
$$
|\nab \bn|^2 n_i = \eps_{ijk}\rd_l n_j \left(  \eps_{kab} n_a\rd_l n_b \right),
$$
where $\eps_{ijk}$ ($i,j,k=1,2,3$) is the Levi-Civita symbol. 
For $k=1,2,3$, 
the divergence of the later components are computed as
\begin{align*}
\rd_l (\eps_{kab} n_a \rd_l n_b)
&=
 \eps_{kab} \left(
\nab n_a \cdot \nab n_b + n_a \Del n_b
\right) \\
&=
\eps_{kab} n_a \left( 2r\nab \times \bn 
- (h+kn_3)\boe_3
\right)_b \in L^2(\R^2),
\end{align*}
where we used \eqref{EF1.5} in the last equality. 
Therefore, for $k,l=1,2,3$, we can write
$$
\eps_{kab} n_a \rd_l n_b
= F^{(1)}_{kl} + F^{(2)}_{kl}
$$
with 
\begin{align}\label{EZA1}
F^{(1)}_{kl}\in H^1_0(B), && 
F^{(2)}_{kl} \in L^2(B),&& \rd_l F^{(2)}_{kl} =0
\end{align}
(see Lemma 2.1.1 in \cite{Sohr} for example). 
For $i=1,2,3$, 
since $\eps_{ijk} \rd_l n_j F_{kl}^{(1)}\in L^p(B)$ for $1<p<2$, 
the elliptic regularity theory implies that the solution $w_{1i}$ to the Dirichlet problem
\begin{align*}
-\Del w_{1i} = \eps_{kab} n_a F_{kl}^{(1)} \quad \text{in } B,&& 
w_{1i} = 0 
\quad \text{on } \rd B
\end{align*}
satisfies $w_{1i}\in W^{2,p}_0(B) \subset C (B)$ (see Sections 9 in \cite{GilTru}). 
Letting $\bm{w}_1 = (w_{11},w_{12},w_{13})$,  we have
\begin{align*}
-\Del \left( \bom - \bm{w}_0 - \bm{w}_1 \right)_i = 
\eps_{ijk} \rd_l n_j F^{(2)}_{kl}
\quad \text{in } B\quad (i=1,2,3),&& 
\bom - \bm{w}_0 -\bm{w}_1 = 0 
\quad \text{on } \rd B.
\end{align*}
On the other hand, since 
$(\eps_{ll'} F_{kl'}^{(2)})_{l=1,2}$ is a curl-free vector field on $B$ for $k=1,2,3$, a distributional Poincar\'e lemma (see Lemma 2.2.2 in \cite{Sohr} for example) implies that there exists $A_k\in H^1(B)$ such that
$$
F^{(2)}_{kl} = \eps_{ll'} 
\rd_{l'} A_{k}.
$$
Therefore, 
we have
$$
-\Del \left( \bom - \bm{w}_0 - \bm{w}_1 \right)_i = 
\eps_{ijk} 
\left(
\eps_{ll'} \rd_ln_j \rd_{l'} A_{k}
\right).
$$
Therefore, by the compensation properties (see, for example, \cite{Wen80}, 
Lemma A.1 in \cite{BreCor84}, Theorem 3.1.3 in \cite{Helein})
we obtain $\bom- \bm{w}_0 -\bm{w}_1 \in H^1(B:\R^3)\cap C(B:\R^3)$, yielding $\bom\in C(B:\R^3)$. 
If we take a smaller ball $B'\subset B$, we may suppose $\sup_{x,y\in B'} |\bom (x) -\bom(y)|$ arbitrary small. 
Then the results in Section 8 in \cite{LadUra} can be applied for \eqref{E3.6}, concluding $\bom\in C^\infty(\R^2:\R^3)$.
\end{proof}

\subsection{Linear estimates}
\label{S4.2}
%
%
In this section, we give a proof of Proposition \ref{Pc2.1}. 
Let us recall the setting. 
Given $\xi_*\in X$ and $\be>0$, we consider the operator 
\begin{equation}\label{EZC1}
F^*F + \frac 12 \ovl{\xi}_*^2 + \be^2,
\end{equation}
where 
$$
F= \rd_\rho + \frac{\cos \te}{\rho} = \rd_\rho + \frac{\rho^2 -4r^2}{\rho(\rho^2+4r^2)},
\qquad
\ovl{\xi_*} = \xi_* + \frac{4r}{(\rho^2+4r^2)^{1/2}}.
$$
%
%
%
Obviously, \eqref{EZC1} is a positive operator, and thus invertible. 
Now we introduce the norm
$$
\nor{f}{X_{\xi_*}} = \left( \nor{f}{X}^2 + \frac 12 \nor{\ovl{\xi_*}f}{L^2}^2 \right)^{1/2}.
$$
Then, the main claim of this section is the following:
\begin{prop}\label{PZC1}
There exist constants $C,\del_0,\be_*>0$ such that for $\xi_*\in X$ with $\nor{\xi_*}{X} \le \del_0$, the following estimates hold:
\begin{equation}\label{EZC2}
\nor{(F^*F + \frac 12 \ovl{\xi_*}^2 + \be^2)^{-1} f}{X_{\xi_*}} \le \frac{C}{\be^{1-s}} 
\nor{\rho^s f}{L^2},\qquad 0\le s\le 1,
\end{equation}
\begin{equation}\label{EZC3}
\begin{aligned}
\nor{(F^*F + \frac 12 \ovl{\xi_*}^2 + \be^2)^{-1} f}{X_{\xi_*}} \le \frac{C}{\be^{1+s}} 
\left(
\nor{\rho^{1-s} \rd_\rho f}{L^2} +
\nor{\rho^{-s} f}{L^2} 
+
\nor{\rho^{1-s} \ovl{\xi_*} f}{L^2}
\right),\\
0<s \le 1
\end{aligned}
\end{equation}
for $f:(0,\infty) \to \R$ and $0<\be<\be_*$.
\end{prop}
\begin{remark}
Given $\xi_*\in X$ with $\nor{\xi_*}{X}\le \del$, \eqref{EZC2} and \eqref{EZC3} hold for all $\be>0$ with a constant $C$ possibly dependent on $\xi_*$. Proposition \ref{PZC1} asserts the uniformity of $C$ in the regime of small $\be$. 
\end{remark}
In what follows, we abbreviate $\xi_*$, $\ovl{\xi}_*$ as $\xi$, $\ovl{\xi}$, respectively. 
The proof of Proposition \ref{PZC1} proceeds as follows. 
We rewrite the inverse of \eqref{EZC1} as
$$
\left(F^*F + \frac 12 \ovl{\xi}^2 + \be^2\right)^{-1} 
= 
(I+R_\xi(\be)W)^{-1} R_{\xi}(\be)
$$
where
$$
R_\xi (\be) = (-\Del_\xi + \be^2)^{-1},\qquad 
-\Del_\xi = -\rd_\rho^2-\frac 1\rho \rd_\rho + \frac{1}{\rho^2} + \frac 12 \ovl{\xi}^2,\qquad 
W= -\frac{32r^2}{(\rho^2+4r^2)^2}.
$$
Then we will see that $\nor{R_\xi(\be)}{X_\xi}$ is estimated by the same bound as in \eqref{EZC2} and \eqref{EZC3} (see Proposition \ref{P4.1} below). 
Thus it suffices to prove $X_\xi\to X_\xi$-estimate for $(I+R_\xi(\be)W)^{-1}$, uniform with respect to $\xi$ and $\be$. 
This outline of proof is the same as that in \cite{GusWan21}, where similar estimates are proved without an extra potential $\frac12 \ovl{\xi}^2$. 
However, 
in our situation, 
we will require additional argument to handle this new factor.
First, unlike the case in \cite{GusWan21}, the limit of $R_\xi(\be)$ as $\be\to 0$ is no longer explicit due to the presence of the implicit potential $\frac 12\ovl{\xi}^2$. 
Therefore, we need to replace the argument in an abstract manner by means of functional analytic techniques. 
Second, the uniformity of the estimate for $(I+R_\xi(\be)W)^{-1}$ in $\xi$ is a new factor compared to \cite{GusWan21}, which is proved by contradiction argument (see Lemma \ref{LZC1}).  
%
%
Before going to the detail, 
we introduce some notation. 
For nonnegative integer $k$, 
define
$$
\nor{f}{X_\xi^k} = \inp{f}{(-\Del_\xi)^k f}_{L^2}^{1/2}.
$$
Notice that $X_\xi^0= L^2$ and $X_\xi^1 = X_\xi$. 
Also, we denote by $X_\xi^*$ the dual space of $X_\xi$. 
%
We identify a function $f:(0,\infty)\to \R$ with the map $g\mapsto \int_0^\infty fg \rho d\rho$ if it is bounded from $X_\xi$ to $\R$. 
For this reason, we denote the coupling of $f\in X_\xi^*$ and $g\in X_\xi$ by $\inp{f}{g}_{L^2}$ if there is no risk of confusion. 
Then, it holds that
\begin{equation}\label{E4.1}
  \nor{f}{X_\xi^*} \le \nor{\rho f}{L^2},\qquad \nor{\rd_\rho f}{X_\xi^*}\le C\nor{f}{L^2}.
\end{equation}
Indeed, for $g\in X$, we have
\begin{gather*}
\left|\inp{f}{g}_{L^2}\right| \le \nor{\rho f}{L^2} \nor{\frac g\rho}{L^2} 
\le
\nor{\rho f}{L^2} \nor{g}{X}
\le \nor{\rho f}{L^2} \nor{g}{X_\xi},
\\
\left|\inp{\rd_\rho f}{g}_{L^2} \right|
= \left|\inp{f}{-\rd_\rho g -\frac 1\rho g}_{L^2}\right|
\le C \nor{f}{L^2}\nor{g}{X_\xi},
\end{gather*}
and thus \eqref{E4.1} follows by definition.

\subsubsection{Resolvent estimates}

\begin{prop}\label{P4.1}
There exists $C>0$ 
such that for all $\be>0$ and $\xi\in X$, we have
\begin{equation}\label{E4.2}
\nor{R_\xi (\be) f}{X_\xi} \le 
\frac{C}{\be^{1-s}} \nor{\rho^s f}{L^2},\qquad 0\le s\le 1,
\end{equation}
\begin{equation}\label{E4.3}
\nor{R_\xi (\be) f}{X_\xi} \le \frac{C}{\be^{1+s}} \left(\nor{\rho^{1-s} \rd_\rho f}{L^2} 
+ \nor{\rho^{-s}f}{L^2}
+ \nor{\rho^{1-s} \ovl\xi f}{L^2}
\right),\qquad 0<s\le 1
\end{equation}
for $\be>0$, $f:(0,\infty)\to\R$. 
\end{prop}
\begin{proof}
We divide the proof into several steps.\par
\textit{Step 1.} 
For integer $k\ge 0$, we first prove 
\begin{align}\label{EZD1}
\nor{R_\xi(\be) f}{X_\xi^k} \le \frac{1}{\be^2} \nor{f}{X^k_\xi}
,&&
\nor{R_\xi(\be) f}{X_\xi^{k+1}} \le \frac{1}{\be} \nor{f}{X^k_\xi}
.
\end{align}
In particular, the first inequality in \eqref{EZD1} with $k=1$ implies \eqref{E4.3} with $s=1$, while 
the second inequality in \eqref{EZD1} with $k=0$ yields \eqref{E4.2} with $s=0$. 
For $f\in X_\xi^k$, let $g= R_\xi(\be)f$. Then
\begin{equation}\label{EZD2}
\inp{(-\Del_\xi)^k g}{f}_{L^2} 
=
\inp{(-\Del_\xi)^k g}{(-\Del_\xi+\be^2) g}_{L^2}
=
\nor{g}{X_\xi^{k+1}}^2 + \be^2 \nor{g}{X_\xi^k}^2 
\ge \be^2 \nor{g}{X_\xi^k}^2.
\end{equation}
Hence
$$
\be^2 \nor{g}{X_\xi^k}^2 \le  \left|\inp{(-\Del_\xi)^k g}{f}_{L^2}\right| \le \nor{g}{X_\xi^k} \nor{f}{X_\xi^k}, 
$$
yielding the first inequality in \eqref{EZD1}. 
If we use \eqref{EZD2} again, it follows that
$$
\nor{g}{X_\xi^{k+1}}^2
\le \nor{g}{X_\xi^k} \nor{f}{X_\xi^k} 
\le \frac{1}{\be^2} 
\nor{f}{X_\xi^k}^2,
$$
implying the second inequality in \eqref{EZD1} as desired.\par
\textit{Step 2.} \eqref{E4.2} with $s=1$. We show the stronger bound
\begin{equation}\label{E4.5}
  \nor{R_\xi(\be) f}{X_\xi} \le C \nor{f}{X_\xi^*},
\end{equation}
which combined with \eqref{E4.1},  implies the desired estimate.  
For $f\in X_\xi$, let $g=R_\xi(\be) f$. Then, \eqref{E4.5} immediately follows since
$$
      \nor{g}{X_\xi}^2 \le
      \nor{g}{X_\xi}^2 + \be^2 \nor{g}{L^2}^2
      =\inp{g}{(-\Del_\xi + \be^2) g}_{L^2}
      = \inp{g}{f}_{L^2}
      \le \nor{g}{X_\xi} \nor{f}{X_\xi^*}.
    $$
\textit{Step 3.} \eqref{E4.2} and \eqref{E4.3} for $0<s<1$. 
We prove this via interpolation using dyadic decomposition in space variable. 
Let $\chi \in C_0^\infty(\R)$ be a cut-off function with
\begin{align*}
\chi(\rho)=1 \quad (|\rho|\le 1),&&
\chi(\rho)=0 \quad (|\rho|\ge 2),&&
0\le \chi \le1\quad \text{on } \R.
\end{align*}
Define $\chi_k(\rho) = \chi(\frac{\rho}{2^{k}}) - \chi(\frac{\rho}{2^{k+1}})$. 
Then we have
$\sum_{k\in\Z} \chi_k=1$ on $\R\setminus \{0\}$, and $\supp \chi_k \subset \{2^{k-1}\le |\rho|\le 2^{k+1} \}$. 
We also define $\tilde{\chi}_k = \chi_{k-1}+\chi_k + \chi_{k+1}$. 
Using these, we decompose as
\begin{align*}
&
\nor{\rho^{1-s} \rd_\rho f}{L^2}^2 + 
\nor{\rho^{-s} f}{L^2}^2 
+
\nor{\rho^{1-s}\ovl\xi f}{L^2}^2\\
&
\sim
\sum_{n\in \Z} \left( 2^{2(1-s)k} \nor{\tilde{\chi}_k 
\rd_\rho f}{L^2}^2
+ 2^{-2sk} \nor{\tilde{\chi}_k 
f}{L^2}^2 
+
2^{2(1-s)k}
\nor{\tilde\chi_k \ovl\xi f}{L^2}^2
\right)
\end{align*}
for $f\in X_\xi$. 
We estimate by 
using \eqref{E4.2} with $s=0$ and \eqref{E4.3} with $s=1$ as
\begin{align*}
&\nor{R_\xi (\be)f}{X_\xi}^2 
\sim \sum_{k\in\Z} \nor{R_\xi (\be) \chi_k 
f}{X_\xi}^2 \\
&\le 
\sum_{k\in\Z} 
\nor{R_\xi (\be) \chi_k 
f}{X_\xi}^{2(1-s)} 
\nor{R_\xi (\be) \chi_k 
f}{X_\xi}^{2s} \\
&\le C
\sum_{k\in\Z} 
\left( \be^{-1} \nor{\chi_k f}{L^2}
\right)^{2(1-s)}
\be^{-4s}
\left( \nor{ \rd_\rho( \chi_k 
f) }{L^2} + 2^{-k} \nor{\chi_k f}{L^2} 
+ \nor{\chi_k \ovl\xi f}{L^2}  \right)^{2s} \\
&\le C 
\be^{-2-2s} \sum_{k\in\Z}  \nor{\chi_k f}{L^2}^{2(1-s)}
\left( \nor{ \chi_k 
\rd_\rho f }{L^2}^{2s} + 2^{-2sk} \nor{\tilde\chi_k f}{L^2}^{2s} 
+ \nor{\chi_k \ovl\xi f}{L^2}^{2s} 
\right) \\
&\le C 
\be^{-2-2s} \sum_{k\in\Z} 
\left( 2^{2(1-s)k} \nor{ \tilde{\chi}_k 
 \rd_\rho f }{L^2}^{2} + 2^{-2sk} \nor{\tilde{\chi}_k f}{L^2}^2
+ 2^{2(1-s)k} \nor{\chi_k \ovl\xi f}{L^2}^{2} 
\right) \\
&\le C 
\be^{-2-2s} \left(
\nor{\rho^{1-s} \rd_\rho f}{L^2}^2 + \nor{\rho^{-s} f}{L^2}^2 
+ \nor{\rho^{1-s} \ovl\xi f}{L^2}^2
\right), 
\end{align*}
which concludes \eqref{E4.3}. 
The proof of \eqref{E4.2} goes similarly as above. 
\end{proof}

\subsubsection{Uniform estimate for $(I+R_\xi (\be)W)^{-1}$}

\begin{prop}\label{P4.2}
There exist $C>0$, $\del_0>0$, $\be_*>0$ 
such that 
\begin{equation}\label{EZE1}
\nor{(I+R_\xi (\be) W)^{-1}}{X_\xi \to X_\xi} \le C
\end{equation}
for all $\xi\in B_X(\del_0)$ 
and $0<\be<\be_*$. 
\end{prop}

We give the proof by piling up several small lemmas.  
First of all, we write $W$ as 
$$
W =W_1 W_2,\qquad W_1 := \frac{1}{\rho^2},\qquad
W_2 := - \frac{32r^2 \rho^2}{ (\rho^2+4r^2)^{2}}.
$$

\begin{lem}\label{LA1}
Let $T_1:\rho L^2\to X_\xi$ be the bounded operator defined by
$$
\inp{T_1 f}{g}_{X_\xi} = 
\inp{W_1 f}{g}_{L^2},\qquad (f\in \rho L^2,\ g\in X_\xi).
$$
Then $R_\xi (\be) W_1 f \wto T_1 f$ in $X_\xi$ for any $f\in \rho L^2$. 
(Later we will see that this convergence will be upgraded into normed sense.)
\end{lem}

\begin{proof}
For $f\in \rho L^2$ and $g\in X$, we observe
\begin{align*}
\inp{R_\xi (\be) W_1 f}{g}_{X_\xi} 
&= \inp{-\Del_\xi (-\Del_\xi + \be^2  )^{-1} W_1 f}{g}_{L^2} \\
&= \inp{W_1 f}{g}_{L^2} - \be^2 \inp{R_\xi (\be) W_1 f}{g}_{L^2}.
\end{align*}
Accordingly, we define $T_1\in \boB(\rho L^2, X_\xi)$ by
$$
\inp{T_1 f}{g}_{X_\xi} = \inp{W_1 f}{g}_{L^2},
$$
which is well-defined 
since 
\begin{align*}
|\inp{W_1 f}{g}_{L^2}| 
\le \int_0^\infty \frac{1}{\rho^2} |fg| \rho d\rho
\le C \nor{\frac{f}{\rho}}{L^2} \nor{g}{X_\xi}. 
\end{align*}
Then we have
\begin{equation}\label{E4.6}
\inp{(R_\xi(\be) W_1 -T_1) f}{g}_{X_\xi} = -\be^2 \inp{R_\xi(\be) W_1 f}{g}_{L^2}.
\end{equation}
Now it suffices to show that \eqref{E4.6} converges to $0$ as $\be\to 0$ for any $f\in \rho L^2$ and $g\in X_\xi$. 
For this, we may assume $g\in C_0^\infty (0,\infty)$, since if we take $\{g_n\}\subset 
C_0^\infty (0,\infty)$ with $g_n\to g$ in $X_\xi$,
\begin{align*}
\be^2 |\inp{R_\xi (\be) W_1 f}{g-g_n}_{L^2}|
&\le \be^2 \nor{R_\xi(\be)W_1 f}{X_\xi^*} \nor{g-g_n}{X_\xi} \\
&\le \nor{W_1 f }{X_\xi^*} \nor{g-g_n}{X_\xi} \le \nor{\frac f\rho}{L^2} \nor{g-g_n}{X_\xi} \xrightarrow{n\to\infty} 0
\end{align*}
uniformly in $\be$. 
Here, we have
\begin{align*}
|\inp{R_\xi(\be) W_1 f}{g}_{L^2}| &
\le \nor{ \rho W_1 f}{L^2} \nor{\frac 1\rho R_\xi(\be) g}{L^2} \le C \nor{\rho W_1 f}{L^2} \nor{\rho g}{L^2}, 
\end{align*}
where the implicit constant is uniform in $\be$. Hence the convergence follows.
\end{proof}

Let $T_1$ be as in Lemma \ref{LA1}. Then it holds that
\begin{equation}\label{E4.7}
R_\xi (\be) W_1 -T_1 = \be^2 R_\xi(\be) T_1.
\end{equation}
Indeed, 
by the resolvent identity,
$$
R_\xi(\be) W_1 - R_\xi(\gam) W_1 = R_\xi(\be) (\be^2 -\gam^2) R_\xi(\gam) W_1,
$$
which converges weakly to \eqref{E4.7} as $\gam\to 0$. 
%
Next, we define
$$
T \equiv T_\xi := T_1 W_2 = \rm{w-}\lim_{\be\to 0} R_\xi (\be) W.
$$
By this definition, we only know that the image of $T$ is included in $X_\xi$. However, the following result shows decay property of $Tf$ for $f\in \rho L^2$.
\begin{lem}
For $f\in \rho L^2$, it holds that  $\rho^{1-s} \rd_\rho T f$, $\rho^{-s} T f 
,\rho^{1-s} \ovl{\xi} f
\in L^2$ for $s\in (0,1)$. Moreover, there exists $C>0$ independent of $\xi$ such that
\begin{equation}\label{E4.8}
\nor{\rho^{1-s} \rd_\rho T f}{L^2} + \nor{\rho^{-s} T f}{L^2} 
+\frac 12 \nor{\rho^{1-s} \ovl\xi T f}{L^2}
\le C \nor{\rho^{2-s} W f}{L^2} \le C \nor{\frac f\rho}{L^2} .
\end{equation}
\end{lem}
\begin{proof}
First suppose that $\rho^{1-s} \rd_\rho Tf$, $\rho^{-s}Tf$, $\rho^{1-s} \ovl{\xi} Tf\in L^2$. 
Noting that
$$
\inp{Tf}{\rho^{2-2s}Tf}_{X_\xi} 
=
\inp{\rd_\rho Tf}{\rd_\rho 
(\rho^{2-2s}Tf)}_{L^2} 
+ 
\nor{\rho^{-s}Tf}{L^2}^2
+
\frac12 \nor{\rho^{1-s} \ovl{\xi} Tf}{L^2}^2,
$$
we have
\begin{equation}\label{EZD4}
\begin{aligned}
&\nor{\rho^{1-s} \rd_\rho T f}{L^2}^2 + \nor{\rho^{-s} T f}{L^2}^2 + \frac 12
\nor{\rho^{1-s}\ovl\xi Tf}{L^2}^2 \\
&= \nor{\rho^{1-s} \rd_\rho T f}{L^2}^2 + 
\inp{T f}{\rho^{2-2s} Tf}_{X_\xi} - \inp{\rd_\rho T f}{\rd_\rho (\rho^{2-2s} T f)}_{L^2}\\
&= \inp{T f}{\rho^{2-2s} Tf}_{X_\xi} - 
(2-2s) \inp{\rd_\rho T f}{\rho^{1-2s} Tf}_{L^2} \\
&= \inp{W f}{\rho^{2-2s} Tf}_{L^2} - (2-2s) \inp{\rd_\rho T f}{\rho^{1-2s} Tf}_{L^2}\\
&\le \nor{\rho^{2-s} W f}{L^2} \nor{\rho^{-s} T f}{L^2} 
+ (1-s) \left(
\nor{\rho^{1-s} \rd_\rho T f}{L^2}^2 + \nor{\rho^{-s} Tf}{L^2}^2
\right) ,
\end{aligned}
\end{equation}
and hence
\begin{align*}
s \left(
\nor{\rho^{1-s} \rd_\rho T f}{L^2}^2 + \nor{\rho^{-s} Tf}{L^2}^2
\right) 
+ \frac 12 \nor{\rho^{1-s}\tilde\xi Tf}{L^2}^2
&\le 
\frac 2s \nor{\rho^{2-s}W f}{L^2}^2 +
\frac s2 \nor{\rho^{-s} Tf}{L^2}^2,
\end{align*}
implying \eqref{E4.8} as desired. \par
To verify the above argument, we finally check 
$\rho^{1-s} \rd_\rho Tf$, $\rho^{-s}Tf$, $\rho^{1-s} \ovl{\xi} Tf\in L^2$. 
Noting that $Tf$ is the weak limit of $R_\xi(\gam) Wf$ as $\gam\to 0$, 
we prove this by obtaining uniform bounds in $\gam$. 
First of all, we claim
\begin{equation}\label{EZD3}
\rho^{1-s} \rd_\rho R_\xi(\gam) Wf,\ 
\rho^{-s}  R_\xi(\gam) Wf,\ 
\rho^{1-s} \ovl{\xi} R_\xi(\gam) Wf \in L^2 \qquad 
\text{for } \gam>0,\ 
0\le s\le 1.
\end{equation}
Indeed, by $f\in \rho L^2$, Proposition \ref{P4.1} implies \eqref{EZD3} with $s=0$. 
On the other hand, we have
$$
\rho R_\xi(\gam) Wf = 
R_\xi (\gam) \left(
\rho Wf - 2\rd_\rho R_\xi(\gam) Wf - \frac 1\rho R_\xi(\gam) Wf
\right),
$$
yielding \eqref{EZD3} with $s=1$. The rest case of $s$ is followed by interpolation. 
Next, the same estimate as \eqref{EZD4} with $Tf$ replaced by $R_\xi (\gam) f$ gives 
\begin{align*}
&\nor{\rho^{1-s} \rd_\rho R_\xi (\gam) W f}{L^2} + \nor{\rho^{-s} R_\xi (\gam) W f}{L^2} 
+ \nor{\rho^{1-s} \ovl\xi R_\xi (\gam) W f}{L^2} \le C \nor{\rho^{2-s} Wf}{L^2} .
\end{align*}
Since the bound is independent of $\gam$, 
there exists $\{\gam_n\}_{n=1}^\infty$ with $\gam_n\to 0$ as $n\to\infty$ such that 
$\rho^{1-s} \rd_\rho R_\xi(\gam_n) Wf$, 
$\rho^{-s}  R_\xi(\gam_n) Wf$, and $\rho^{1-s} \ovl{\xi} R_\xi (\gam_n)Wf$ are weakly convergent in $L^2$. 
However, 
by Lemma \ref{LA1}, 
the limits should coincide 
with $\rho^{1-s} \rd_\rho Tf$, $\rho^{-s} Tf$, 
$\rho^{1-s} \ovl\xi  Tf$, respectively. 
Hence we obtain the desired conclusion.
\end{proof}

\begin{lem}\label{LA4}
$R_\xi (\be) W \xrightarrow{\be\to 0} T$ in $\boB (X_\xi)$, where the 
convergence is uniform in $\xi$.
\end{lem}
\begin{proof}
Let $f\in \rho L^2$, and take $s\in (0,1)$. By \eqref{E4.3}, \eqref{E4.7} and \eqref{E4.8},
\begin{align*}
\nor{R_\xi (\be) W f - Tf}{X_\xi} 
&=
\be^2 \nor{R_\xi (\be) T f}{X_\xi} \\
&\le C \be^2 \be^{-1-s} 
\left(\nor{\rho^{1-s} \rd_\rho T f}{L^2} + 
\nor{\rho^{-s} T f}{L^2} + \nor{\rho^{1-s} \ovl\xi Tf}{L^2}  \right) \\
&\le C \be^{1-s} \nor{\frac f\rho}{L^2} \le C \be^{1-s} \nor{f}{X_\xi} \xrightarrow{\be\to 0} 0.
\end{align*}
\end{proof}

\begin{lem}\label{LA2}
$W_2:X_\xi \to \rho L^2 $ is a compact operator. 
In particular, $T$ and $R_\xi(\be) W \in \boB(X_\xi)$ are compact.
\end{lem}

\begin{proof}
Let $\{f_n\}$ be a bounded sequence in $X_\xi$. 
Then by taking a subsequence, $\{ f_n\}$ has a weak limit $f$ in $X_\xi$. 
It suffices to show $\frac 1\rho W_2 f_n \to \frac 1\rho W_2f$ in $L^2$ by taking further subsequence. 
For any $0<\rho_0<\rho_1<\infty$, we claim that $\{f_n\}_{n=1}^\infty$ satisfies the conditions of the Ascoli-Arzel\`a theorem. 
The uniform boundedness follows immediately by $\nor{f_n}{L^\infty} \le C \nor{f_n}{X} =C$. On the other hand,
$$
|f_n(\tau_0) -f_n(\tau_1)| = 
\left|
\int_{\tau_0}^{\tau_1} \rd_\rho f_n (\rho) d\rho
\right|
\le 
\nor{f_n}{X} 
\left(\int_{\tau_0}^{\tau_1} \frac 1{\rho^2} \rho d\rho \right)^{\frac 12}
= \nor{f_n}{X} \left(\log \left( \frac{\tau_1}{\tau_0} \right)\right)^{\frac 12},
$$
implying equicontinuity of $\{ f_n\}$ on $[\rho_0,\rho_1]$. 
Therefore, by the Ascoli-Arzel\`a theorem, up to a subsequence, $\{f_n\}$ converges uniformly on $[\rho_0,\rho_1]$. 
By the diagonal argument, it follows that there exist $f:(0,\infty)\to\R$ and a subsequence of $\{f_n\}$ such that 
\begin{equation}\label{EZF2}
f_n \xrightarrow{n\to\infty} f \qquad \text{uniformly in any closed 
interval } I\subset (0,\infty). 
\end{equation}
On the other hand, since $\{\frac 1\rho W_2 f_n\}$ is bounded in $L^2$, we may suppose that 
$\{\frac 1\rho W_2 f_n\}$ is weakly convergent, and its limit should be $\frac 1\rho W_2 f$. 
Next, 
we observe that 
$\{\frac 1\rho W_2 f_n\}$ is tight, in the sense that 
for any $\eps>0$, there exists $\rho_0, \rho_1>0$ such that
\begin{equation}\label{E4.9}
\int_{(0,\rho_0) \cap (\rho_1,\infty)}  
\frac{|W_2 f_n|^2}{\rho^2}
 \rho d\rho
< \eps
\end{equation}
for all $n\in\N$. Indeed, 
\begin{gather*}
\int_0^{\rho_0} \frac{|W_2f_n|^2}{\rho^2} \rho d\rho 
= 
\int_0^{\rho_0} \frac{2^{10} r^4\rho^4 }{(\rho^2 +4r^2)^4} \frac{|f_n|^2}{\rho^2} \rho d\rho 
\le 
C\rho_0^4 \sup_{n\in\N} \nor{f_n}{X}^2,
\\
\int_{\rho_1}^\infty \frac{|W_2f_n|^2}{\rho^2} \rho d\rho 
= 
\int_{\rho_1}^\infty \frac{2^{10} r^4\rho^4 }{(\rho^2 +4r^2)^4} \frac{|f_n|^2}{\rho^2} \rho d\rho 
\le 
\frac{C}{(\rho_1^2 +1)^2} \sup_{n\in\N} \nor{f_n}{X}^2,
\end{gather*}
and hence \eqref{E4.9} follows. 
Similarly, we have
\begin{equation}\label{EZF3}
\int_{(0,\rho_0) \cap (\rho_1,\infty)}  
\frac{|W_2 f|^2}{\rho^2}
 \rho d\rho
< \eps,
\end{equation}
using the same $\rho_0$, $\rho_1$ as above. Therefore, combining \eqref{EZF2}, \eqref{E4.9} and \eqref{EZF3} leads to $\frac 1\rho W_2 f_n \to \frac 1\rho W_2 f$ strongly in $L^2$ as desired.
\end{proof}

\begin{lem}\label{LA3}
There exists $\del_0>0$ such that if $\nor{\xi}{X} \le \del_0$, then 
$$
\Ker (I+T) =0, \quad 
\text{and}\quad \Ker (I+R_\xi (\be) W)=0.
$$
\end{lem}

\begin{proof}
Let $f\in \Ker (I+T)$. 
We take $\{f_n\}\subset C_0^\infty(0,\infty)$ with $f_n \xrightarrow{n\to\infty} f$ in $X_\xi$. 
Then for $\vph\in C_0^\infty (0,\infty)$, we have
\begin{equation}\label{E4.12}
\begin{aligned}
\inp{(I+T)f_n}{\vph}_{X_\xi} &= 
\lim_{\be\to 0} \inp{(I+R_\xi(\be) W)f_n}{\vph}_{X_\xi}
= \lim_{\be\to 0} \inp{
-\Del_\xi (I+R_\xi(\be) W)f_n}{\vph}_{L^2} \\
&= \lim_{\be\to 0} \left( \inp{(-\Del_\xi + W )f_n}{\vph}_{L^2} - \be^2 \inp{R_\xi (\be) Wf_n}{\vph}_{L^2} \right) \\
&= \inp{(-\Del_\xi + W)f_n}{\vph}_{L^2}\\ 
&= \left\langle
f_n, \left(-\rd_\rho^2 -\frac 1\rho \rd_\rho + \frac 1{\rho^2}\right)  \vph\right\rangle_{L^2} +\left\langle\left(\frac 12 \ovl\xi^2 + W\right)f_n,\vph\right\rangle_{L^2}.
\end{aligned}
\end{equation}
Taking $n\to\infty$, we have
\begin{equation}
0 = \left\langle f, 
\left(
-\rd_\rho^2 -\frac 1\rho \rd_\rho 
 + \frac 1{\rho^2}\right)\vph\right\rangle_{L^2} + 
\left\langle
\left(\frac 12 \ovl\xi^2 + W\right)f,\vph\right\rangle_{L^2},
\end{equation}
which means
\begin{equation}\label{E4.10}
(
-\Del_\xi + W)f =0 \quad \text{in } 
\boD'(0,\infty).
\end{equation}
Now we claim that $f\in C^1(0,\infty)$. 
For this, we use elliptic regularity theory. 
First, let $\Om_R = \{ x\in \R^2|  \frac 1R< |x| <R\}$ for fixed $R>2$, and 
define $\tilde{f}(x)= e^{i\psi} f(\rho)$ with $x=\rho e^{i\psi}$. Then
$$
\Del \tilde{f} 
= \frac 12 \ovl\xi^2 \tilde{f} + W \tilde{f} 
\quad \text{in } \boD'(\Om_R).
$$
In particular, we have $\tilde{f}$, $\Del \tilde{f}\in H^1(\Om_R)$. 
Thus, the elliptic regularity theory implies 
$\tilde{f}\in H^3 (\tilde\Om_R)$ where 
$\tilde\Om_R = \{ x\in \R^2|  \frac 2R< |x| < \frac R2 \}$. 
(See, for example, \cite{GilTru}.)
Especially, we have $\tilde{f}\in C^1(\tilde\Om_R)$ and thus $f\in C^1(\frac 2R ,\frac R2 )$. 
Since $R$ is arbitrary, we have $f\in C^1(0,\infty).$\par
%
%
Next, by using the same sequence $\{f_n\}_{n=1}^\infty$ as above, we have
\begin{align*}
\inp{(I+T)f_n}{f_n}_{X_\xi} &= 
\lim_{\be\to 0} \inp{(I+R_\xi(\be) W)f_n}{f_n}_{X_\xi}
= \lim_{\be\to 0} \inp{-\Del_\xi (I+R_\xi(\be) W)f_n}{f_n}_{L^2} \\
&= \lim_{\be\to 0} \left( \inp{(-\Del_\xi + W )f_n}{f_n}_{L^2} - \be^2 \inp{R_\xi (\be) Wf_n}{f_n}_{L^2} \right) \\
&= \inp{(-\Del_\xi + W)f_n}{f_n}_{L^2}\\
&= \nor{F f_n}{L^2}^2
+ \int_0^\infty \frac 12 \ovl\xi^2 |f_n|^2 \rho d\rho 
\ge \int_{0}^{\infty} \frac 12 \ovl\xi^2 |f_n|^2 \rho d\rho.
\end{align*} 
Taking $n\to\infty$, we obtain
\begin{equation}\label{E4.11}
0 \ge \int_{0}^{\infty} \frac 12 \ovl\xi^2 |f|^2 \rho d\rho.
\end{equation}
Here, by the assumption $\nor{\xi}{X} \le \del_0 \ll 1$, there is an interval $[\rho_0,\rho_1]$ such that $\ovl{\xi} >0$ on $[\rho_0,\rho_1]$. 
Therefore, \eqref{E4.11} implies 
$f(\rho) =0 $ on $[\rho_0,\rho_1]$. 
Since $f$ is a $C^1$-solution to \eqref{E4.10}, 
we have 
$f\equiv 0$ on $(0,\infty)$ by the uniqueness of the solutions.\par
Next, let $f\in \Ker (I+R_\xi (\be)W)$. We first note that $f\in L^2$, since the dual estimate of \eqref{E4.2} with $s=0$ implies
$$
\nor{f}{L^2} =\nor{R_\xi (\be)Wf}{L^2} \le \frac C\be \nor{Wf}{X_\xi^*} \le 
\frac C\be \nor{\rho Wf}{L^2} \le \frac{C}{\be} \nor{\frac{f}{\rho}}{L^2} <\infty.
$$
Here we take a sequence $\{f_n\}\subset C_0^\infty(0,\infty)$ with 
$f_n\xrightarrow{n\to\infty} f$ in $X_\xi$. Then
\begin{align*}
\inp{(-\Del_\xi +\be^2)(I+R_\xi(\be) W)f_n}{f_n}_{L^2} 
&= 
\inp{(-\Del_\xi + \be^2 +W)f_n}{f_n}_{L^2} \\
&=
\nor{Ff_n}{L^2}^2 + \be^2 \nor{f_n}{L^2}^2 + 
\int_0^\infty \frac 12 \ovl\xi^2 |f_n|^2 \rho d\rho. 
\end{align*}
On the other hand
\begin{align*}
&\inp{(-\Del_\xi +\be^2)(I+R_\xi (\be) W)f_n}{f_n}_{L^2} \\
&= 
\inp{(I+R_\xi (\be) W)f_n}{f_n}_{X_\xi } 
+ \be^2 \inp{(I+R_\xi (\be) W)f_n}{f_n}_{L^2}  \\
& =
\inp{(I+R_\xi (\be) W)f_n}{f_n}_{X_\xi } 
+ \be^2 \nor{f_n}{L^2}^2
+ \be^2 \inp{R_\xi (\be) Wf_n}{f_n}_{L^2}.
\end{align*}
Hence
$$
\inp{(I+R_\xi (\be) W)f_n}{f_n}_{X_\xi } 
+ \be^2 \inp{R_\xi (\be) Wf_n}{f_n}_{L^2}
\ge \int_{0}^{\infty} \frac 12 \ovl\xi^2 |f_n|^2 \rho d\rho,
$$
As $n\to\infty$, we have
$$
\inp{(I+R_\xi (\be) W)f}{f}_{X_\xi } 
+ \be^2 \inp{R_\xi (\be) Wf}{f}_{L^2}
\ge \int_{0}^{\infty} \frac 12 \ovl\xi^2 |f|^2 \rho d\rho,
$$
where the convergence is justified by
\begin{align*}
&\left|\inp{R_\xi (\be) W f_n}{f_n}_{X_\xi} - \inp{R_\xi (\be)W f}{f}_{X_\xi}
\right| \\
&\le \left|\inp{R_\xi (\be) W(f_n-f)}{f_n}_{X_\xi}\right| + 
\left|\inp{R_\xi(\be)W f}{f-f_n}_{X_\xi}\right| \\
&\le 
 \nor{R_\xi(\be) W(f_n-f)}{X_\xi} \nor{f_n}{X_\xi}
+ \nor{R_\xi(\be) Wf}{X_\xi} \nor{f_n-f}{X_\xi} \\
&\le C \nor{f-f_n}{X_\xi } 
(\nor{f_n}{X_\xi} 
+ \nor{f}{X_\xi}
)\xrightarrow{n\to\infty} 0,\\
&\left|\inp{R_\xi (\be) W f_n}{f_n}_{L^2} - \inp{R_\xi(\be)W f}{f}_{L^2}
\right| \\
&\le \left|\inp{R_\xi(\be) W(f_n-f)}{f_n}_{L^2}\right| + 
\left|\inp{R_\xi (\be)W f}{f-f_n}_{L^2}\right| \\
&=
\left|\inp{ W(f_n-f)}{R_\xi(\be)f_n}_{L^2}\right| + 
\left|\inp{W f}{R_\xi (\be)(f-f_n)}_{L^2}\right| \\
&\le 
 \nor{ W(f_n-f)}{X_\xi^*} \nor{R_\xi (\be) f_n}{X_\xi}
+ \nor{ Wf}{X_\xi^*} \nor{R_\xi (\be)(f_n-f)}{X_\xi} \\
&\le 
\frac{C}{\be^2} \nor{f-f_n}{X_\xi} (\nor{f_n}{X_\xi} +\nor{f}{X_\xi}) \xrightarrow{n\to\infty} 0.
\end{align*}
Hence, by assumption of $f$, we have
\begin{align*}
0 &= 
\inp{(I+R_\xi (\be) W)f}{f}_{X_\xi}  + \be^2 \inp{(I+R_\xi(\be) W)f}{f}_{L^2} \\
&=
\inp{(I+R_\xi (\be) W)f}{f}_{X_\xi} 
+ \be^2 \nor{f}{L^2}^2
+ \be^2 \inp{R_\xi (\be) Wf}{f}_{L^2} \\
&\ge 
\be^2 \nor{f}{L^2}^2
+\int_0^\infty \frac 12 \ovl\xi^2 |f|^2 \rho d\rho,
\end{align*}
which implies $f=0$.
\end{proof}
%
\begin{lem}\label{LZC1}
There exists $\del_0>0$ such that $I+T_\xi$ is invertible for $\xi\in B_X (\del_0)$, and
\begin{equation}\label{EFeb5}
\sup_{\nor{\xi}{X} \le \del_0 } \nor{(I+T_\xi)^{-1}}{X_\xi\to X_\xi} <\infty.
\end{equation}
\end{lem}
\begin{proof}
Since $T_\xi$ is compact by Lemma \ref{LA2}, the Fredholm alternative implies that  
$I+T_\xi$ is invertible if and only if 
its kernel is $0$, which is true by Lemma \ref{LA3}.\par
Suppose that \eqref{EFeb5} is false. 
Then there exists $\{\xi_n\}_{n=1}^\infty$ with $\nor{\xi_n}{X} \to 0$ such that
$$
\nor{(I+T_{\xi_n})^{-1}}{X_{\xi_n}\to X_{\xi_n}} \ge n.
$$
In particular, there is $f_n\in X_{\xi_n}$ such that 
$$
\nor{f_n}{X_{\xi_n}} \le \frac 2n,\qquad \text{and}\qquad 
\nor{(I+T_{\xi_n})^{-1} f_n}{X_{\xi_n}} =1.
$$
Write $g_n := (I+T_{\xi_n})^{-1} f_n$. 
Since 
$
\nor{g_n}{X} \le \nor{g_n}{X_{\xi_n}} =1
$, 
$g_n$ has a weak limit $g$ in $X$ if we take a subsequence. 
Now we claim that $g\equiv 0$. 
Let $\vph\in C_0^\infty (0,\infty)$. 
By the same computation as \eqref{E4.12}, we obtain
$$
\inp{(I+T_{\xi_n}) g_n}{\vph}_{X_{\xi_n}} = 
\inp{g_n}{(-\rd_\rho^2 - \frac 1\rho \rd_\rho + \frac{1}{\rho^2}) \vph}_{L^2} + \inp{(\frac 12 \ovl\xi_n^2 +W) g_n}{\vph}_{L^2}.
$$
Here, we observe that
\begin{gather*}
|\inp{(I+T_{\xi_n}) g_n}{\vph}_{X_{\xi_n}}| = 
|\inp{f_n}{\vph}_{X_{\xi_n}}| \le \nor{f_n}{X_{\xi_n}} \nor{\vph}{X_{\xi_n}}
\le \frac{C}{n} \xrightarrow{n\to\infty} 0,
\\
|\inp{ (\ovl\xi_n^2 - \frac{16r^2}{\rho^2+4r^2})  g_n}{\vph}_{L^2}| 
\le C \nor{\xi_n}{L^\infty} \nor{g_n}{L^\infty} 
\xrightarrow{n\to\infty} 0,
\\
\inp{(\frac{8r^2}{\rho^2+4r^2} +W) g_n}{\vph}_{L^2}
\xrightarrow{n\to\infty} \inp{(\frac{8r^2}{\rho^2+4r^2} +W) g}{\vph}_{L^2}.
\end{gather*}
Hence we have
\begin{equation}\label{E4.13}
\begin{aligned}
0 
&= \inp{g}{(-\rd_\rho^2-\frac 1\rho \rd_\rho + \frac 1{\rho^2}) \vph}_{L^2} + \inp{(\frac{8r^2}{\rho^2+4r^2} +W) g}{\vph}_{L^2} \\
&= \inp{Fg}{F\vph}_{L^2} + \inp{ \frac{8r^2}{\rho^2+4r^2} g}{\vph}_{L^2}.
\end{aligned}
\end{equation}
Now we take a sequence $\{\vph_n\}_{n=1}^\infty\subset C_0^\infty (0,\infty)$ such that $\vph_n\to g$ in $X$ as $n\to\infty$. If we apply \eqref{E4.13} for $\vph=\vph_n$, and take $n\to\infty$, we obtain
$$
0= \nor{Fg}{L^2}^2 +
\int_0^\infty \frac 12 \sin^2 \te |g|^2 \rho d\rho,
$$
which yields $g=0$.\par
Here, by definition of $T_{\xi_n}$,
\begin{equation}\label{E4.14}
1= \nor{g_n}{X_{\xi_n}}^2 = \inp{f_n - T_{\xi_n} g_n}{g_n}_{X_{\xi_n}}
= \inp{f_n }{g_n}_{X_{\xi_n}} - \inp{Wg_n}{g_n}_{L^2}.
\end{equation}
Now we claim that the right hand side of \eqref{E4.14} converges to $0$ as $n\to\infty$, which leads to a contradiction. 
First, we have
$$
|\inp{f_n }{g_n}_{X_{\xi_n}}| \le 
\nor{f_n}{X_{\xi_n}} \nor{g_n}{X_{\xi_n}} \le \frac 2n \xrightarrow{n\to\infty} 0.
$$
On the other hand, we apply Lemma \ref{LA2} with $\ovl\xi=0$, which implies that $\frac 1\rho W_2 g_n \to 0$ strongly in $L^2$, if we take subsequence if necessary. Therefore, 
$$
\inp{Wg_n}{g_n}_{L^2} 
\le \nor{\frac 1\rho W_2 g_n}{L^2} 
\nor{g_n}{X_{\xi_n}} 
\xrightarrow{n\to\infty} 0.
$$
Hence we complete the proof.
\end{proof}

\begin{lem}
\label{LZC2}
$I+R_\xi (\be)W$ is invertible, and 
$(I+R_\xi (\be)W)^{-1} \xrightarrow{\be\to 0} (I+T)^{-1}$ in $\boB(X_\xi)$, 
where the convergence is uniform in $\xi \in B_X (\del_0)$. 
\end{lem}
\begin{proof}
The invertibility of $I+R_\xi(\be)W$ follows in the same way as that of $I+T$. 
Note that
$$
(I+R_\xi (\be)W)^{-1} = (I+ T- T +R_\xi(\be)W )^{-1} 
=  \left[ I - (I+T)^{-1} (T-R_\xi (\be)W ) \right]^{-1} (1+T)^{-1} .
$$
Thus, we have
\begin{align*}
&(I+R_\xi (\be)W)^{-1} - (1+T)^{-1}\\
&=
\sum_{n=1}^\infty \left( (I+T)^{-1} (T-R_\xi (\be)W) \right)^n
(1+T)^{-1}.
\end{align*}
By Lemma \ref{LA4}, 
we have $\nor{R_\xi(\be)W -T}{X_\xi\to X_\xi} \to 0$ as $\be\to 0$ uniformly in $\xi$. 
Now, let $\eps>0$ be arbitrary number. Thanks to Lemma \ref{LZC1}, 
there exists $\be_\eps>0$ such that if $\be < \be_\eps$, then
$$
\nor{(I+T)^{-1} (R_\xi(\be)W -T)}{X_\xi\to X_\xi} 
\le 
\nor{(1+T)^{-1}}{X_\xi\to X_\xi}
\nor{(R_\xi(\be)W -T)^{-1}}{X_\xi\to X_\xi} 
\le \eps
$$
uniformly in $\xi\in B_X (\del_0)$. 
Therefore, 
$$
\nor{(I+R_\xi (\be)W)^{-1} - (1+T)^{-1}}{X_\xi\to X_\xi} \le C \nor{(I+T)^{-1}}{X_\xi\to X_\xi} \eps,
$$
which, together with Lemma \ref{LZC1},  proves the uniform convergence as claimed.
\end{proof}
\begin{proof}[Proof of Proposition \ref{P4.2}]
By Lemma \ref{LZC2}, $I+R_\xi(\be)W$ is invertible for all $\be>0$, and there exists $\be_*>0$, independent of $\xi\in B_X (\del_0)$, such that 
$$
\sup_{0<\be\le\be_*} 
\nor{(I+R_\xi(\be)W)^{-1}}{X_\xi \to X_\xi} 
\le 2 \nor{(I+T)^{-1}}{X_\xi\to X_\xi}\le C.
$$
By Lemma \ref{LZC1}, this bound is
uniform in $\xi\in B_\xi(\del_0)$, which completes the proof. 
\end{proof}

\textbf{Acknowledgements} 
The authors wish to thank Stephen Gustafson for an insightful discussion which led us recognize the $F^*F$ structure in the linearized operator. 
S.~Ibrahim is supported by the NSERC grant No. 371637-2025. I.~Shimizu is supported by JSPS KAKENHI Grant Number 23KJ1416. \medskip\par

\textbf{Statements and Declarations} The authors have no relevant financial or non-financial interests to disclose.


\begin{thebibliography}{99}

\bibitem{Aka19} 
Akahori,~T., Ibrahim~S., Ikoma,~N., Kikuchi,~H., Nawa,~H.: 
Uniqueness and nondegeenracy of ground states to nonlinear scalar field equations involving the Sobolev critical exponent in their nonlinearities for high frequencies, Calc. Var. \textbf{58}(4), 120 (2019)

\bibitem{BarSinRosSch20} 
Barton-Singer,~B., Ross,~C., Schroers,~B.~J.: 
Magnetic skyrmions at critical coupling, 
Comm. Math. Phys. \textbf{375}, 2259--2280 (2020)

\bibitem{BerLio83} 
Berestycki,~H., Lions,~P.-L.: 
Nonlinear scalar 
field equations, I: Existence of a ground state
, Arch. Rational Mech. Anal. \textbf{82}(4), 315--345 (1983)

\bibitem{BerMurSim20} 
Bernand-Mantel,~A., Muratov,~C.~B., Simon,~T.~M.: 
Unraveling the role of dipolar versus Dzyaloshinskii-Moriya interactions in stabilizing compact magnetic skyrmions, 
Phys. Rev. B \textbf{101}, 045416 (2020)

\bibitem{BerMurSim21} 
Bernand-Mantel,~A., Muratov,~C.~B., Simon,~T.~M.: 
A quantitative description of skyrmions in ultrathin ferromagnetic films and rigidity of degree $\pm 1$ harmonic maps from $\R^2$ to $\mathbb{S}^2$, 
Arch. Rat. Mech. Anal. \textbf{239}(1), 219--299 (2021)

\bibitem{BogHub94}
Bogdanov,~A.~N., Hubert,~A.: 
Thermodynamically stable magnetic vortex states in magnetic crystals, 
J. Magn. Magn. Mater. \textbf{138}, 255--269 (1994)

\bibitem{BogHub99}
Bogdanov,~A.~N., Hubert,~A.: 
The stability of vortex-like structures in uniaxial ferromagnets, 
J. Magn. Magn. Mater. \textbf{195}(1), 182--192 (1999)

\bibitem{BogYab89} Bogdanov,~A.~N., Yablonskii,~D.~A.: 
Thermodynamically stable ``vortices'' in magnetically ordered crystals. The mixed state of magnets, Zh. Eksp. Teor. Fiz. \textbf{95}, 178--182 (1989)

\bibitem{BreCor84} 
Brezis,~H., Coron~J.-M.: Multiple solutions of $H$-systems and Rellich's conjecture, 
Comm. Pure Appl. Math. \textbf{37}(2), 149--187 (1984)

\bibitem{ButLemBea18} 
B\"uttner,~F., Lemesh,~I., Beach,~S.~D.~G.: 
Theory of isolated magnetic skyrmions: from fundamentals to room temperture applications, 
Sci. Rep. \textbf{8}, 4464 (2018)

\bibitem{ColGus20} 
Coles,~M., Gustafson,~S.: 
Solitary waves and dynamics for subcritical perturbations of energy critical NLS, 
Publ. Res. Inst. Math. Sci. \textbf{56}(4), 647--699 (2020)

\bibitem{DorMel17} 
D\"oring,~L., Melcher,~C.: 
Compactness results for static and dynamic chiral skyrmions near the conformal limit, 
Calc. Var. \textbf{56}(3). 60 (2017)

\bibitem{Felsager} 
Felsager,~B.: 
Geometry, Particles, and Fields, 
Springer-Verlag, 
1998.

\bibitem{FerCroSam13} 
Fert,~A., Cros,~V., Sampaio,~J.: 
Skyrmions on the track, 
Nature. Nanotechnology \textbf{8}, 160--161 (2013)

\bibitem{GilTru} 
Gilbarg,~D., Trudinger,~N. S.: 
Elliptic Partial Differential Equations of Second Order, 
Springer-Verlag, Berlin, Heiderberg, 2001.

\bibitem{GusWan21} 
Gustafson,~S., Wang,~L.: 
Co-rotational chiral magnetic skyrmions near harmonic maps, J. Funct. Anal. \textbf{280}(4), 108867 (2021)

\bibitem{Helein} 
H\'elein,~F.: 
Harmonic Maps, Conservation Laws, and Moving Frames, 
Cambridge University Press, New York, 2002. 

\bibitem{HubSch} 
Hubert,~A., Sch\"afer,~R.: Magnetic Domains: The Analysis of Magnetic Microstructures, Springer-Verlag, 1998. 

\bibitem{IbrShi23}
Ibrahim,~S., Shimizu,~I.: 
Phase transition threshold and stability of magnetic skyrmions, 
Commun. Math. Phys. \textbf{402}, 2627--2640 (2023)

\bibitem{KomMelVen20} 
Komineas,~S., Melcher,~C., Venakides,~S.: 
The profile of chiral skyrmions of small radius
, Nonlinearity \textbf{33}(7), 3395--3408 (2020)

\bibitem{KomMelVen21} 
Komineas,~S., Melcher,~C., Venakides,~S.:  
Chiral skyrmions of large radius
, Phys. D \textbf{418}, 132842 (2021)

\bibitem{KomMelVen23} 
Komineas,~S., Melcher,~C., Venakides,~S.:
Chiral magnetic skyrmions across length scales
, New J. Phys. \textbf{25}, 023013 (2023)


\bibitem{LadUra} 
Ladyzhenskaya,~O.~A., Ural'tseva,~N.~N.: 
Linear and Quasilinear Elliptic Equations, Academic Press, 1968. 

\bibitem{LanLif35} 
Landau,~L., Lifshitz,~E.: 
On the theory of the dispersion of magnetic permeability in ferromagnetic bodies, Phys. Zeitsch. der Sw. \textbf{8}, 153--169 (1935)

\bibitem{Leo16}
Leonov,~A.~O., Monchesky,~T.~L., Romming, N., Kubetzka,~A., Bogdanov,~A.~N., Wiesendanger,~R.: 
The properties of isolated chiral skyrmions in thin magnetic films, New J. Phys., 065003 (2016)

\bibitem{LiMel18} 
Li, X., Melcher, C.: Stability of axisymmetric chiral skyrmions, J. Funct. Anal. \textbf{275}(10), 2817--2844 (2018)

\bibitem{Mel14}
Melcher,~C.: 
Chiral skyrmions in the plane, Proc. R. Soc. A \textbf{470}(2172), 20140394 (2014)

\bibitem{MurSimSla} 
Muratov,~C.~B., Simon,~T. M., Slastikov,~V.~V.: 
Existence of higher degree minimizers in the magnetic skyrmion problem, Arch. Ration. Mech. Anal. \textbf{249}(5), Paper No. 61, 31 (2025).

\bibitem{NagTok13} 
Nagaosa,~N., Tokura,~Y.: Topological properties and dynamics of magnetic skyrmions, Nature Nanotech. \textbf{8}, 899--911 (2013)

\bibitem{Sch19} 
Schroers,~B.~J.: 
Gauged sigma models and magnetic skyrmions, SciPost Phys. \textbf{7}, 030 (2019)

\bibitem{Struwe} 
Struwe,~M.: 
Variational Methods: Applications to Nonlinear Partial Differential Equations and Hamiltonian Systems, 2nd ed., Springer-Verlag, 1996.

\bibitem{Sohr} 
Sohr,~H.: 
The Navier-Stokes Equations: An Elementary Functional Analytic Approach, Springer Basel AG, 2001.

\bibitem{Wen80} 
Wente,~H.-C.: 
Large solutions to the volume constrained Plateau problem, Arch. Rational Mech. Anal. \textbf{75}(1), 59--77 (1980)

\end{thebibliography}
\end{document}